\documentclass{amsart}
\usepackage{amsfonts,amssymb,amsmath,amsthm}
\usepackage{url}
\usepackage{enumerate}

        \usepackage[french,english]{babel}
\usepackage[T1]{fontenc}
        \usepackage{amstext}
        \usepackage{mathtools}
       \usepackage{dsfont}

\newtheorem{theo}{Th\'eor\`eme}[section]
\newtheorem{lem}[theo]{Lemme}
\newtheorem{prop}[theo]{Proposition}
\newtheorem{cor}[theo]{Corollaire}
\newtheorem*{rem*}{Remarque}
\newtheorem{rem}[theo]{Remarque}

\newcommand{\beqyn}{\begin{eqnarray}}
\newcommand{\enqyn}{\end{eqnarray}}
\newcommand{\blem}{\begin{lem}}
\newcommand{\elem}{\end{lem}}
\newcommand{\brop}{\begin{prop}}
\newcommand{\erop}{\end{prop}}
\newcommand{\bcor}{\begin{cor}}
\newcommand{\ecor}{\end{cor}}
\newcommand{\brem}{\begin{rem}}
\newcommand{\erem}{\end{rem}}
\newcommand{\brems}{\begin{rem*}}
\newcommand{\erems}{\end{rem*}}
\newcommand{\benum}{\begin{enumerate}}
\newcommand{\enum}{\end{enumerate}}
\newcommand{\nident}{\noindent}

\newcommand{\btheo}{\begin{theo}}
\newcommand{\etheo}{\end{theo}}
\newcommand{\bdem}{\begin{proof}}
\newcommand{\edem}{\end{proof}}

\newcommand{\rmA}{\mathrm{A}}

\newcommand{\rmE}{\mathrm{E}}
\newcommand{\rmF}{\mathrm{F}}

\newcommand{\rmI}{\mathrm{I}}

\newcommand{\rmK}{\mathrm{K}}

\newcommand{\Hom}{\mathrm{Hom}}
\newcommand{\Gal}{\mathrm{Gal}}
\newcommand{\Lie}{\mathrm{Lie}}
\newcommand{\Id}{\mathrm{Id}}

\newcommand{\Ad}{\mathrm{Ad}}

\newcommand{\Tr}{\mathrm{Tr}}
\newcommand{\End}{\mathrm{End}}

\newcommand{\diag}{\mathrm{diag}}
\newcommand{\Gl}{\mathrm{GL}}

\newcommand{\Res}{\mathrm{Res}}

\newcommand{\vol}{\mathrm{vol}}
\newcommand{\Rel}{\mathrm{Re}}

\newcommand{\ind}{\mathds{1}}

\newcommand{\upla}{\underline{\rho}}
\newcommand{\uphi}{\upphi}
\newcommand{\al}{\alpha}

\newcommand{\la}{\lambda}
\newcommand{\La}{\Lambda}

\newcommand{\brtau}{\bar \tau}

\newcommand{\inv}{^{-1}}

\newcommand{\bilif}{\langle \cdot,\cdot \rangle}

\newcommand{\matx}[4]{\begin{pmatrix} #1 & #2 \\ #3 & #4 \\ \end{pmatrix}}
\newcommand{\dsl}{\displaystyle \left(}
\newcommand{\rb}{\right)}

\newcommand{\bsl}{\backslash}

\newcommand{\rar}{\rightarrow}
\newcommand{\irar}{\xrightarrow{\sim}}
\newcommand{\hrar}{\hookrightarrow}
\newcommand{\smin}{\smallsetminus}
\newcommand{\sbs}{\subseteq}
\newcommand{\sbn}{\subsetneq}

\newcommand{\sps}{\supseteq}

\newcommand{\ps}[1]{\prescript{#1}{}}

\newcommand{\psl}[1]{\prescript{}{#1}}

\newcommand{\gll}{\mathfrak{gl}}

\newcommand{\all}{\mathfrak{a}}

\newcommand{\kl}{\mathfrak{k}}
\newcommand{\fl}{\mathfrak{l}}
\newcommand{\gl}{\mathfrak{g}}

\newcommand{\ml}{\mathfrak{m}}
\newcommand{\nl}{\mathfrak{n}}

\newcommand{\tl}{\mathfrak{t}}
\newcommand{\ul}{\mathfrak{u}}

\newcommand{\ol}{\mathfrak{o}}

\newcommand{\N}{\mathbb{N}}
\newcommand{\Z}{\mathbb{Z}}
\newcommand{\A}{\mathbb{A}}
\newcommand{\C}{\mathbb{C}}

\newcommand{\Q}{\mathbb{Q}}
\newcommand{\R}{\mathbb{R}}

\newcommand{\Gm}{\mathbb{G}_m}
\newcommand{\Ga}{\mathbb{G}_a}

\newcommand{\calF}{\mathcal{F}}

\newcommand{\calI}{\mathcal{I}}
\newcommand{\calJ}{\mathcal{J}}
\newcommand{\calK}{\mathcal{K}}

\newcommand{\calO}{\mathcal{O}}
\newcommand{\calP}{\mathcal{P}}
\newcommand{\calR}{\mathcal{R}}
\newcommand{\calS}{\mathcal{S}}

\newcommand{\calV}{\mathcal{V}}

\newcommand{\cad}{c'est-\`a-dire }

\newcommand{\nin}{\notin}

\newcommand{\tlPhi}{\tilde \Phi}

\newcommand{\tlP}{\widetilde{P}}
\newcommand{\tlR}{\widetilde{R}}
\newcommand{\tlQ}{\widetilde{Q}}
\newcommand{\tlS}{\widetilde{S}}

\newcommand{\tlU}{\widetilde{U}}
\newcommand{\tlV}{\widetilde{V}}

\newcommand{\tlul}{\widetilde{\ul}}

\newcommand{\tltl}{\widetilde{\tl}}

\newcommand{\tlzero}{\widetilde{0}}
\newcommand{\tlone}{\widetilde{1}}
\newcommand{\tltwo}{\widetilde{2}}

\newcommand{\htau}{\hat \tau}
\newcommand{\hf}{\hat f}

\newcommand{\hDelta}{\widehat \Delta}

\newcommand{\huphi}{\hat \upphi}

\newcommand{\brP}{\overline{P}}

\newcommand{\brLa}{\overline{\La}}

\usepackage{upgreek}
\usepackage{mathabx}

\numberwithin{equation}{section}
\author{Micha\l \ Zydor}
\address{
Université Paris Diderot
Institut de Mathématiques de Jussieu-Paris Rive Gauche
UMR7586
Bâtiment Sophie Germain
Case 7012
75205 PARIS Cedex 13
France}
\email{michal.zydor@imj-prg.fr}


\begin{document}

\title[La formule des traces de Jacquet-Rallis]{La variante infinitésimale 
de la formule des traces de Jacquet-Rallis pour
les groupes unitaires}

\selectlanguage{english}
\begin{abstract} 
We establish an infinitesimal version of the 
Jacquet-Rallis trace formula for unitary groups. 
Our formula is obtained by integrating a 
truncated kernel à la Arthur. 
It has a geometric side which is a 
sum of distributions $J_{\mathfrak{o}}$ indexed by classes of elements 
of the Lie algebra of $U(n+1)$ stable by $U(n)$-conjugation
as well as the "spectral side" 
consisting of the Fourier transforms 
of the aforementioned distributions. 
 We prove that the distributions $J_{\mathfrak{o}}$ 
are invariant and depend only on the choice of 
the Haar measure on $U(n)(\mathbb{A})$. 
For regular semi-simple classes $\mathfrak{o}$, $J_{\mathfrak{o}}$ is 
a relative orbital integral of Jacquet-Rallis. 
For classes $\mathfrak{o}$ called relatively regular semi-simple, we express $J_{\mathfrak{o}}$
in terms of relative orbital integrals regularised by means of zêta functions.
\end{abstract}
\selectlanguage{french}
\begin{abstract}
Nous établissons une variante infinitésimale 
de la formule des traces de Jacquet-Rallis 
pour les groupes unitaires. 
Notre formule s'obtient par 
intégration d'un noyau tronqué 
à la Arthur. Elle possède un côté géométrique 
qui est une somme de distributions $J_{\mathfrak{o}}$ indexée 
par les classes d'éléments de l'algèbre de Lie de $U(n+1)$ 
stables par $U(n)$-conjugaison
ainsi qu'un "côté spectral" 
formé des transformées de Fourier 
des distributions précédentes. 
On démontre que les distributions $J_{\mathfrak{o}}$ 
sont invariantes et ne dépendent que du choix 
de la mesure de Haar sur $U(n)(\mathbb{A})$. 
Pour des classes $\mathfrak{o}$ semi-simples régulières $J_{\mathfrak{o}}$ 
est une intégrale
orbitale relative de Jacquet-Rallis.
Pour les classes $\mathfrak{o}$ dites relativement semi-simples régulières, 
on exprime $J_{\mathfrak{o}}$
en terme des intégrales orbitales relatives régularisées à l'aide des fonctions zêta.
\end{abstract}
\maketitle

\section{Introduction}

\subsection{Contexte}
Dans \cite{jacqrall} Jacquet et Rallis proposent une approche 
\`a la conjecture globale de Gan-Gross-Prasad pour les groupes unitaires 
via une formule des traces relative. 
Soient $\rmE/\rmF$ 
une extension quadratique de corps globaux, 
$\sigma \in \Gal(\rmE/\rmF)$ non-trivial, 
 $\A$ l'anneau des adèles de $\rmF$ 
et $\eta : \rmF^{*}\backslash\A^{*} \rightarrow \C^{*}$ le 
caract\`ere quadratique associé à 
l'extension $\rmE/\rmF$ 
par la théorie du corps de classes. 
Fixons $V$ et $W$ des $\rmE$-espaces hermitiens 
non-dégénérés,
de dimensions $n$ et $n+1$ respectivement, 
munis de l'inclusion $V \hrar W$
et notons $U = U(V)$ et $\tlU = U(W)$ les groupes unitaires associés. 
L'inclusion $V \hrar W$ réalise $U$ comme un sous-groupe de 
$\tlU$.
Les formules des traces qu'ils suggèrent les mènent à 
postuler l'égalité entre les intégrales formelles suivantes:
\begin{equation}\label{JRRTF}
\int_{U(\rmF)\backslash U(\A)}
\sum_{\gamma \in \tlU(\rmF)}f(x^{-1}\gamma x)
 dx =  \\
\int_{\Gl_{n}(\rmF)\backslash\Gl_{n}(\A)}
\sum_{\gamma' \in S_{n+1}(\rmF)}g(x^{-1}\gamma' x)
\eta(\det x)dx
\end{equation}
où  $S_{n+1}(\rmF) = \{x \in \Gl_{n+1}(\rmE): x\sigma(x) =1\}$, 
$f$ et $g$ sont des fonctions lisses \`a support compact 
sur $U(\A)$ et $S_{n+1}(\A)$ respectivement 
qui vérifient certaines conditions de compatibilité des 
intégrales orbitales locales. 

Ces formules des traces sont liées à la conjecture 
de Gan-Gross-Prasad de la façon suivante. 
Chacune des deux intégrales peut être vue
 comme une période d'un noyau automorphe. 
 Au moins formellement, chacune admet 
 deux développements: l'un, qui exploite la décomposition 
 spectrale du noyau automorphe, 
 est la somme des périodes des représentations 
 du spectre automorphe, 
 l'autre, qui utilise l'action de 
 $U$ sur $\tlU$ et de $\Gl_{n}$ sur $S_{n+1}$, 
 est une somme des intégrales orbitales relatives. 
 Pour $\Gl_{n}$, les périodes automorphes 
 sont reliés à des valeurs spéciales de fonctions 
 $L$ par la théorie de Rankin-Selberg. 
 L'égalité (\ref{JRRTF}) relie  
 les périodes des groupes unitaires à 
 celles de $\Gl_{n}$ et donc aux 
 valeurs spéciales de fonctions $L$. 
 Cette égalité devrait résulter de l'égalité 
 des décompositions géométriques. 
 On renvoie 
à \cite{ggp} pour un énoncé précis 
de la conjecture de Gan-Gross-Prasad.

Zhang \cite{zhang2}, proposition 2.11,
établit 
l'égalité (\ref{JRRTF}) pour des fonctions 
vérifiant certaines contraintes locales. 
Il s'en est servi pour démontrer une partie 
substantielle de la conjecture de Gan-Gross-Prasad. 
Cette même formule 
des traces a ensuite permis à Zhang
\cite{zhang1}
d'obtenir certains cas 
du raffinement de la conjecture de Gan-Gross-Prasad 
dû à Ichino et Ikeda \cite{ichinoIkeda} pour les groupes orthogonaux et
adapté aux groupes unitaires par
N. Harris \cite{harris}.

\subsection{Nos résultats}\label{par:nosResIntr}

Les deux intégrales dans (\ref{JRRTF}) ne 
sont pas convergentes 
en général et elles nécessitent d'être tronquées. 
Dans cet article 
nous nous proposons de donner une version infinitésimale
de la formule des traces relative de Jacquet-Rallis dans 
le cas des groupes unitaires par un processus de troncature 
\`a la Arthur. Le cas des groupes linéaires 
est traité dans \cite{leMoi2}. 

Plus précisément, la formule qu'on tronque 
est
\begin{equation}\label{JRRTFalgebres}
\int_{U(\rmF)\backslash U(\A)}
\sum_{\xi \in \tlul(\rmF)}f(x^{-1}\xi x)dx = 
\int_{U(\rmF)\backslash U(\A)}
k_{f}(x)dx, \quad f \in \calS(\tlul(\A))
\end{equation}
o\`u $\tlul$ est l'alg\`ebre de Lie du groupe unitaire $\tlU$ 
et $\calS(\tlul(\A))$ est l'espace des fonction 
de classe Bruhat-Schwartz sur $\tlul(\A)$ 
(voir le paragraphe \ref{bschwartz}). 

Décrivons notre troncature brièvement. Fixons $P_{0}$ 
un $\rmF$-sous-groupe parabolique minimal de $U$ ainsi qu'une décomposition de Levi $P_{0} = M_{0}N_{0}$ avec $M_{0}$ une partie 
de Levi de $P_{0}$ et $N_{0}$ son radical unipotent. Tout sous-groupe 
parabolique standard $P$ (i.e. $P \sps P_{0}$) admet alors 
une unique décomposition de Levi $P = M_{P}N_{P}$ où $M_{0} \sbs M_{P}$.
Tout sous-groupe parabolique standard de $U$ est défini comme 
le stabilisateur d'un drapeau isotrope dans $V$. 
On associe à tout sous-groupe parabolique standard $P$ 
le sous-groupe parabolique $\tlP$ de $\tlU$ défini 
comme le stabilisateur du même drapeau que $P$ mais vu dans 
l'espace hermitien $W$. On a alors une unique décomposition de 
Levi $M_{\tlP}N_{\tlP}$ de $\tlP$ telle que $M_{\tlP} \sps M_{P}$.
Dans le paragraphe 
\ref{lesInvs} on décrit une décomposition de $\tlul(\rmF)$ 
en classes stables par $U(\rmF)$-conjugaison, on note 
$\calO$ l'ensemble de ces classes. 

Pour un sous-groupe parabolique standard $P$, $\ol \in \calO$ 
et $f \in \calS(\tlul(\A))$
on pose
\[
k_{f,P,\mathfrak{o}}(x) =  \sum_{\xi \in 
\Lie(M_{\tlP})(\rmF) \cap \mathfrak{o}}
\int_{\Lie(N_{\tlP})(\A)}f(x\inv (\xi + U)x)dU, \ 
x \in M_{P}(\rmF)N_{P}(\A)\backslash U(\A).
\]
Le noyau tronqué est défini alors comme
\[
k^{T}_{f,\mathfrak{o}}(x) = 
\sum_{P\supseteq P_{0}} 
(-1)^{d_{P}}\sum_{\delta \in P(\rmF)\backslash U(\rmF)}
\htau_{P}(H_{P}(\delta x) - T)k_{f,P,\mathfrak{o}}(\delta x), \ 
x \in U(\rmF)\backslash U(\A)
\]
où, si l'on pose $\all_{P} := 
\Hom_{\Z}(\Hom_{\rmF}(M_{P}, \Gm), \R)$, alors 
$H_{P} : U(\A) \rar \all_{P}$ 
est l'application de Harish-Chandra qui 
dépend du choix d'un sous-groupe compact maximal dans $U(\A)$, 
$\htau_{P}$ 
est la fonction caractéristique
d'un cône obtus dans $\all_{P}$, 
$d_{P} = \dim_{\R} \all_{P}$ et 
$T$ est un paramètre dans $\all_{0} := \all_{P_{0}}$
(voir le paragraphe \ref{par:prelimstrace}). 
Remarquons que 
Ichino et Yamana \cite{ichYam2}, en s'inspirant de \cite{jlr}, définissent
des périodes régularisés qui apparaissent dans le contexte de 
la conjecture globale de Gan-Gross-Prasad pour les groupes unitaires 
et leur opérateur de troncature $\La_{m}^{T}$
utilise les mêmes ensembles des sous-groupes paraboliques de $U$ 
et $\tlU$ que nous. 
Notre premier résultat, démontré dans la section \ref{sec:convergence}, 
est alors:
\begin{theo}[cf. \ref{thm:MainConv}]
Pour tout
 $T \in T_{+} + \mathfrak{a}_{0}^{+}$ on a
\begin{displaymath}
\sum_{\mathfrak{o} \in \mathcal{O}}
\int_{U(\rmF) \backslash U(\A)}|k_{f,\mathfrak{o}}^{T}(x)|dx < \infty.
\end{displaymath}
\end{theo}
\noindent
Ici $\all_{0}^{+}$ c'est le
cône aigu engendré par les combinaisons linéaires à coefficients 
positifs de copoids dans $\all_{0}$ et 
$T_{+} \in \all_{0}^{+}$ est précisé dans 
le paragraphe \ref{par:prelimstrace}.

Ensuite, on s'intéresse au comportement de l'application $T \mapsto \int_{U(\rmF) \backslash U(\A)}k_{f,\mathfrak{o}}^{T}(x)dx$. 
Dans le le paragraphe \ref{par:asymptChapt} on obtient:
\begin{theo}[cf. \ref{mainQualitThm}]
La fonction
\[
T \mapsto\  J^{T}_{\ol}(f) := 
\int_{U(\rmF) \bsl U(\A)}k_{f,\ol}^{T}(x)dx
\] 
où $T$
parcourt $T_{+} + \mathfrak{a}_{0}^{+}$ 
est un polynôme-exponentielle (voir paragraphe \ref{par:fonsPolExp})
dont la partie purement polynomiale est constante.  
\end{theo}
\noindent
En fait, si l'on remplace 
l'intégrale dans (\ref{JRRTFalgebres}) par
$
\int_{U(\rmF)\backslash U(\A)}
F^{U}(x,T)k_{f}(x)dx
$
où la fonction $F^{U}$ est 
à support compact (voir le paragraphe \ref{par:prelimstrace}) 
alors 
on obtient une expression convergente. 
D'après les résultats de 
\cite{levy}, cette expression est asymptotiquement 
un polynôme-exponentielle 
 en $T$. Notre distribution $J_{\ol}^{T}$ égale alors
 ce polynôme-exponentielle (cf. corollaire \ref{almostPolCor}).

On note $J_{\ol}(f)$ la partie constante du 
polynôme-exponentielle $J_{\ol}^{T}(f)$. 
Il s'avère que la distribution 
$J_{\ol}$ a des propriétés remarquables. On obtient
 \begin{itemize}
\item (cf. Théorème \ref{invarianceTheo}) 
La distribution $J_{\ol}$ est invariante 
par conjugaison par $U(\A)$.
\item (cf. Paragraphe \ref{par:noChoixMade})
La distribution $J_{\ol}$ ne dépend que du choix de la mesure de Haar sur $U(\A)$. 
 \end{itemize}

Dans le paragraphe \ref{par:regOrbsChap}, on constate que 
les distributions $J_{\mathfrak{o}}^{T}$ 
pour les classes 
semi-simples régulières, \cad les classes composées d'une seule $U(\rmF)$-orbite 
de stabilisateur trivial,  ne dépendent pas 
du param\`etre $T$, donc $J_{\ol}^{T} = J_{\ol}$, 
et qu'elles s'expriment comme 
des intégrales orbitales relatives qui 
apparaissent déjà dans \cite{jacqrall}. 

Dans la section
\ref{FourierTransChap} on obtient 
la formule des traces de Jacquet-Rallis infinitésimale pour 
les groupes unitaires.
\begin{theo}[cf. \ref{thm:RTFJRI}] 
Pour tout $f \in \calS(\tlul(\A))$ on a
\begin{displaymath}
\sum_{\mathfrak{o} \in \mathcal{O}}J_{\mathfrak{o}}(f) = 
\sum_{\mathfrak{o} \in \mathcal{O}}J_{\mathfrak{o}}(\hat f).
\end{displaymath}
\end{theo}
\noindent
Ici $\hat f$ 
c'est une transformée de Fourier (on en considère plusieurs)
de $f$. 

La section \ref{sec:orbRrssS} est consacrée à donner une expression plus explicite 
des distributions $J_{\ol}$ pour des $\ol$ dites relativement semi-simples régulières. 
Ces classes sont des réunions finies d'orbites et ses éléments admettent des tores 
pour centralisateurs. 
Soient $\ol$ une telle classe et $P$ un sous-groupe parabolique standard de $U$
tel que l'unique orbite fermée dans $\ol$ intersecte $\Lie (M_{\tlP})(\rmF)$ 
non-trivialement, minimal pour cette propriété.
On choisit des représentants $X_{\calI}$ des orbites de dimension maximale dans $\ol$
de façon qu'ils aient le même centralisateur, noté $T_{0}$, 
qui vérifie $T_{0} \sbs M_{P}$. Le $\rmF$-tore $T_{0}$ est alors 
anisotrope et le quotient $T_{0}(\rmF) \bsl T_{0}(\A)$ est compact. 
Pour $f \in \calS(\tlul(\A))$ et $X_{\calI}$ l'intégrale
\begin{equation*}
\zeta_{\calI}(f)(\la) = 
\vol (T_{0}(\rmF) \bsl T_{0}(\A))
\int\limits_{\mathclap{T_{0}(\A)\bsl U(\A)}}
f(\Ad(x^{-1})X_{\calI})e^{\la(H_{P}(x))}dx, \quad \la \in \Hom_{\R}(\all_{P},\C)
\end{equation*}
converge sur un ouvert
de $\Hom_{\R}(\all_{P},\C)$ qui dépend de $X_{\calI}$, 
et admet un prolongement méromorphe 
à $\Hom_{\R}(\all_{P},\C)$ noté aussi $\zeta_{\calI}(f)$. 
La fonction $\zeta_{\calI}(f)$ peut avoir un pôle en $\la = 0$ mais 
la somme de $\zeta_{\calI}(f)$ sur toutes les orbites de dimension maximale dans $\ol$ 
n'en a pas. En fait, la valeur en $\la = 0$ de cette somme est liée à notre distribution de façon suivante:

\begin{theo}[cf. \ref{thm:theThmOrbs}]\label{thm:theThmOrbs0} 
On a
\[
J_{\ol}(f) = \dsl \sum_{X_{\calI}} \zeta_{\calI}(f)\rb (0)
\]
où la somme porte sur toutes les orbites de dimension maximale dans $\ol$.
\end{theo}
Dans \cite{leMoi2} on définit aussi la notion 
d'une classe relativement semi-simple régulière dans le contexte du groupe linéaire
et on obtient une formule complètement analogue.

La troncature que nous proposons et ses 
propriétés se transposent bien
dans le cas des groupes et donnent
le développement géométrique 
de la formule des traces de Jacquet-Rallis 
pour les groupes unitaires. 
Ceci sera établi dans \cite{leMoi3} 
où on donnera les formules 
des traces relatives de Jacquet-Rallis grossières.
Par ailleurs, notre formule a un intérêt propre. 
Elle peut être un outil puissant pour des questions 
d'analyse harmonique locale. Par 
exemple dans le cas de l'endoscopie classique, 
un analogue simple pour les algèbres de Lie 
de la formule 
des traces d'Arthur joue un rôle central 
dans la preuve du transfert par Waldspurger \cite{wald1}. Notons 
aussi que Zhang \cite{zhang2} passe par les algèbres de Lie
pour démontrer le transfert
dans la formule des traces simple qu'il utilise.

%
\textbf{Remerciements}. Je remercie chaleureusement mon directeur de thèse Pierre-Henri Chaudouard pour ses multiples conseils. 
Je remercie aussi Atsushi Ichino et Shunsuke Yamana pour la possibilité de consulter leur article 
\cite{ichYam2} avant qu'il ne soit publié. Merci à Jacek Jendrej pour son aide dans la preuve du lemme \ref{lem:signCombi}. 
Je remercie finalement le rapporteur pour sa lecture attentive. 
Ce travail a été partiellement soutenu par le projet Ferplay ANR-13-BS01-0012.

\section{Prolégom\`enes}\label{prolegolego}

\subsection{Préliminaires pour la formule des traces}\label{par:prelimstrace}

Soient $\rmF$ un corps de nombres et 
$U$ un $\rmF$-groupe algébrique 
réductif. Pour tout $\rmF$-sous-groupe de Levi $M$ 
de $U$ (\cad un facteur de Levi d'un $\rmF$-sous-groupe parabolique 
de $U$) soit $\calF(M)$ l'ensemble des $\rmF$-sous-groupes 
paraboliques de $U$ contenant $M$ et $\calP(M)$ le sous-ensemble 
de $\calF(M)$ composé des sous-groupes paraboliques 
admettant $M$ comme facteur de Levi.
On fixe 
un sous-groupe de Levi minimal $M_{0}$ de $U$. 
Fixons aussi un $P_{0} \in \calP(M_{0})$. 
On appelle les éléments de $\calF(M_{0})$ 
les sous-groupes paraboliques semi-standards 
et les éléments de $\calF(M_{0})$ 
contenant $P_{0}$ les sous-groupes paraboliques standards. 
Tout $\rmF$-sous-groupe parabolique de $U$ 
est $U(\rmF)$-conjugué à un unique sous-groupe parabolique standard.
On utilisera toujours le symbole $P$, avec des indices éventuellement, 
pour noter un sous-groupe parabolique semi-standard.
Pour tout $P \in \calF(M_{0})$ soit $N_{P}$ le radical unipotent 
de $P$ et $M_{P}$ le facteur de Levi de $P$ contenant $M_{0}$. 
On a alors $P = M_{P}N_{P}$. On note $A_{P}$ 
le tore central de $M_{P}$ déployé sur $\rmF$ 
et maximal pour cette propriété. Les objets associés au sous-groupe 
$P_{0}$ seront noté avec l'indice $0$, on écrit donc 
$N_{0}$, $A_{0}$ au lieu de $N_{P_{0}}$, $A_{P_{0}}$ etc.
De même, quand il n'y aura pas d'ambiguïté on écrit 
$N_{1}$ au lieu de $N_{P_{1}}$ etc, pour un $P_{1} \in \calF(M_{0})$.

Soit  $P \in \calF(M_{0})$.
On définit le $\R$-espace vectoriel 
$\mathfrak{a}_{P} := \Hom_{\Z}(\Hom_{\rmF}(M_{P}, \Gm), \R)$, 
isomorphe à $\Hom_{\Z}(\Hom_{\rmF}(A_{P}, \Gm), \R)$ grâce à l'inclusion $A_{P} \hrar M_{P}$, 
ainsi que son espace dual $\all_{P}^{*} = \Hom_{\rmF}(M_{P}, \Gm) \otimes_{\Z} \R$ et on pose
 \begin{equation}\label{eq:dPDef}
 d_{P} = \dim_{\R} \all_{P}, \quad d_{Q}^{P} = d_{Q} - d_{P}, \ 
 Q \sbs P.
\end{equation}

Si $P_{1}\subseteq P_{2}$, on a un
homomorphisme injectif canonique
$\mathfrak{a}_2^{*} \hookrightarrow \mathfrak{a}_1^{*}$ 
qui donne la projection 
$\mathfrak{a}_1  \twoheadrightarrow \mathfrak{a}_2$, 
dont on note 
$\mathfrak{a}_{1}^{2} =\mathfrak{a}_{P_{1}}^{P_{2}}$ le noyau.
On a aussi l'inclusion 
$\mathfrak{a}_{{2}}\hookrightarrow \mathfrak{a}_{{1}}$, 
qui est une section de
$\mathfrak{a}_1 \twoheadrightarrow \mathfrak{a}_2$, 
grâce à la restriction des caractères de $A_{{1}}$ à
$A_{{2}}$. 
Il s'ensuit  que si 
$P_{1} \subseteq P_{2}$ alors
\begin{equation}\label{eq:decomp}
\mathfrak{a}_1 = \mathfrak{a}_1^2\oplus \mathfrak{a}_2.
\end{equation}
Conformément à cette décomposition, on pose aussi 
$(\all_{1}^{2})^{*} = 
\{\la \in \all_{1}^{*}| \la(H) = 0 \ \forall H \in \all_{2}\}$. 
Grâce à la décomposition (\ref{eq:decomp}) ci-dessus, on considère les espaces
$\all_{1}^{*}$ et $(\all_{1}^{2})^{*}$ comme des sous-espaces de $\all_{0}^{*}$.

Notons $\Delta_{P}^{U} = \Delta_{P}$ 
l'ensemble de racines 
simples pour l'action de $A_{P}$ sur $N_{P}$. 
Il y a une correspondance bijective entre les sous-groupes paraboliques $P_{2}$
contenant $P_{1}$ et les sous-ensembles 
$\Delta_{1}^{2} = \Delta_{P_{1}}^{P_{2}}$ de 
$\Delta_1 = \Delta_{P_{1}}$. En fait, 
$\Delta_1^2$
est l'ensemble de racines simples 
pour l'action de $A_1$ sur $N_1 \cap M_2$ et 
l'on a
\begin{displaymath}
\mathfrak{a}_2 = \{H \in \mathfrak{a}_1| \al(H) = 0 
\ \forall \al \in \Delta_1^2\}. 
\end{displaymath}
En plus $\Delta_{1}^{2}$ (les restrictions de ses éléments 
à $\all_{1}^{2}$) est une base de $(\all_{1}^{2})^{*}$.

Fixons $P_{1} \subseteq P_{2}$ et soit $B \in \calP(M_{0})$ 
contenu dans $P_{1}$. On a alors l'ensemble 
$\Delta_{B}^{\vee} = \{\al^{\vee} \in \all_{0}| \al  \in \Delta_{B}\}$ 
de coracines simples associées aux racines simples $\Delta_{B}$.
À une racine $\al \in \Delta_1^2$ on associe 
de mani\`ere univoque une coracine $\al^{\vee} \in \all_{1}^{2}$ 
de façon suivante: $\al$ est une restriction 
d'une unique $\al_{0} \in \Delta_{B}^{2} \smin \Delta_{B}^{1}$ 
à $\all_{1}^{2}$ et on définit $\al^{\vee}$ comme la projection de 
$\al_{0}^{\vee}$ à $\all_{1}^{2}$. Cela ne dépend pas du choix de $B$.
Notons également $\hDelta_1^2$ et 
$(\hDelta_1^2)^{\vee}$ 
les bases de $(\mathfrak{a}_{{1}}^{{2}})^{*}$ et 
$\mathfrak{a}_{{1}}^{{2}}$ duales \`a 
$(\Delta_{{1}}^{{2}})^{\vee}$ et
$\Delta_{{1}}^{{2}}$ 
 respectivement. 
Si $P_{2} = U$ on note simplement $\Delta_{{1}}, 
\Delta_{{1}}^{\vee}$ etc. 

Soient $P, P_{1}, P_{2} \in \calF(M_{0})$, on note
\begin{displaymath}
\mathfrak{a}_{P}^{+} = 
\{H \in \mathfrak{a}_{P}| \alpha(H) > 0 \
\forall \alpha \in \Delta_{P}\}
\end{displaymath}
et si $P_{1} \subseteq P_{2} $ 
notons $\tau_{1}^{2}$, $\htau_1^2$
les fonction caractéristiques de
\begin{equation*}
\{H \in \mathfrak{a}_1| \alpha(H) > 0 \
\forall \alpha \in \Delta_1^2\}, \quad  
\{H \in \mathfrak{a}_{1}| \varpi(H) > 0 \
\forall \varpi\ \in \hDelta_{{1}}^{{2}}\}
\end{equation*}
respectivement. 
On note $\tau_{P}$ pour $\tau_{P}^{U}$ et 
$\htau_{P}$ pour $\htau_{P}^{U}$. 

Soit $\A = \A_{\rmF}$ l'anneau des adèles de $\rmF$ 
et soit $|\cdot |_{\A}$ la valeur absolue standard sur le groupe 
des idèles $\A^{*}$.
 Pour tout $P \in \calF(M_{0})$, posons 
 $H_{P} : M_{P}(\A) \rightarrow \mathfrak{a}_{P}$
défini comme
\begin{equation}\label{harishChandra}
\langle H_{P}(m),\chi \rangle = \log (|\chi (m)|_{\A}), \quad \chi \in 
\Hom_{\rmF}(M_{P}, \Gm), \ m \in M_{P}(\A).
\end{equation}
C'est un homomorphisme continu et surjectif, 
donc si l'on note $M_{P}(\A_{})^{1}$ 
son noyau, on obtient la suite exacte suivante
\begin{displaymath}
1 \rightarrow M_{P}(\A_{})^{1} \rightarrow M_{P}(\A_{}) 
\rightarrow \mathfrak{a}_{P} \rightarrow 0.
\end{displaymath}
Soit $A_{P}^{\infty}$ la composante neutre du groupe des \(\R\)-points du 
tore déployé et défini sur $\Q$ maximal pour cette propriété dans 
le \(\Q\)-tore 
\(\Res_{\rmF/\Q}A_{P}\). Alors, comme $\rmF \otimes_{\Q} \R$ 
s'injecte dans $\A_{}$, on a 
naturellement $A_{P}^{\infty}\hookrightarrow A_{P}(\A_{}) \hookrightarrow M_{P}(\A_{})$. En plus, 
la restriction de $H_{P}$ \`a $A_{P}^{\infty}$ est 
un isomorphisme donc 
$M_{P}(\A_{})$ est un produit direct de 
$M_{P}(\A_{})^{1}$ et $A_{P}^{\infty}$.

Fixons $K$ un sous-groupe compact maximal admissible de 
$U(\A)$ par rapport à $M_{0}$. La notion d'admissibilité 
par rapport à un sous-groupe de Levi minimal est définie 
dans le paragraphe 1 de \cite{arthur2}. 
On a donc, que pour tout sous-groupe parabolique semi-standard $P$, 
$K \cap M_{P}(\A)$ est admissible dans $M_{P}(\A)$ 
et on obtient aussi la décomposition d'Iwasawa
$U(\A) = P(\A)K = N_{P}(\A)M_{P}(\A) K$ ce 
qui nous permet d'étendre $H_{P}$ \`a $U(\A)$ en posant 
$H_{P}(x) = H_{P}(m)$ où $x = nmk$ avec 
$m \in M_{P}(\A), n \in N_{P}(\A), k \in K$. Dans ce cas $H_{P}(x)$ 
ne dépend pas du choix de $m$.

On fixe un $T_{1} \in -\mathfrak{a}_{0}^{+}$ et 
un compact $\omega \subseteq N_0(\A)M_0(\A)^{1}$ pour 
qu'on puisse définir $F^{P}(x,T)$ pour tout 
$T \in T_{1} + \mathfrak{a}_0^{+}$ comme le fait Arthur dans 
le paragraphe 6 de
\cite{arthur3}.
Pour tous sous-groupes 
paraboliques $P_{1}\subseteq P_{2}$ et 
$T \in T_{1} + \mathfrak{a}_{0}^{+}$ soit 
\begin{equation}\label{eq:A12T}
A_{1,2}^{\infty}(T) =  \{a \in A_{1}^{\infty} \cap M_{{2}}(\A)^{1}| 
\alpha(H{_{1}}(a)-T_{1}) > 0, 
\varpi(H_{{1}}(a) - T) < 0 \ \forall \al \in \Delta_{{1}}^{{2}}, 
\varpi \in \hDelta_{{1}}^{{2}}\}.
\end{equation}
Pour un $T \in T_{1} + \mathfrak{a}_{0}^{+}$ on 
définit la fonction $F^{P}(x,T)$ o\`u $x \in U(\A)$ comme 
la fonction caractéristique de la projection 
de  $\omega A_{0,P}^{\infty}(T)K$ 
sur 
$A_{P}^{\infty}N_{P}(\A)M_{P}(\rmF)\backslash U(\A)$. 
On a alors la propriété que si $x = nmak$ 
o\`u $n \in N_{P}(\A)$, $m \in M_{P}(\A)^{1}$, 
$a \in A_{P}^{\infty}$ et $k \in K$ alors 
$F^{P}(x,T) = F^{P}(m,T)$. Fixons 
aussi un $T_{+} \in \mathfrak{a}_{{0}}^{+}$ tel 
que pour tout $T \in T_{+} + \mathfrak{a}_{{0}}^{+}$ et 
pour tout sous-groupe parabolique $P$, les fonctions 
$F^{P}(\cdot,T)$ sont définies et vérifient le lemme 
6.4 de \cite{arthur3}. En particulier, on peut choisir 
le compact $\omega$, et on le fait, de façon que la 
fonction 
$F^{P}$ restreint \`a $M_{P}(\A)$ soit égale \`a 
$F^{M_{P}}$ relativement aux compacts $\omega \cap M_{P}(\A)$ 
et $K \cap M_{P}(\A)$. 

On note $\Omega$ 
le groupe de Weyl de $(U,A_{{0}})$. 
Celui-ci op\`ere naturellement sur
les espaces $\mathfrak{a}_{0}$ et $\all_{0}^{*}$. 
Pour tout $s \in \Omega$, on choisit un représentant $w_{s}$
dans l'intersection de $U(\rmF)$ 
avec le normalisateur de $A_{{0}}$. 

Pour tous $P,Q \in \calF(M_{0})$
notons 
$\Omega(\mathfrak{a}_{P},\mathfrak{a}_{Q})$ 
l'ensemble des isomorphismes distincts de 
$\mathfrak{a}_{P}$ sur $\mathfrak{a}_{Q}$ obtenus 
par restriction d'un élément de $\Omega$ \`a 
$\mathfrak{a}_{P}$. Pour tout 
$s_{0} \in \Omega(\mathfrak{a}_{P},\mathfrak{a}_{Q})$ 
il existe un unique $s \in \Omega$ tel que la restriction 
de $s$ \`a $\mathfrak{a}_{P}$ soit $s_{0}$ et 
$s^{-1}\al$ est une combinaison linéaire \`a coefficients positifs d'éléments de
$\Delta_{{0}}$ pour 
tout $\al \in \Delta_{{0}}^{Q}$. Ainsi, 
on peut voir $\Omega(\mathfrak{a}_{P},\mathfrak{a}_{Q})$ comme 
un sous-ensemble de $\Omega$. 

Fixons maintenant un sous-groupe parabolique standard 
$P_{1}$. Pour un sous-groupe parabolique standard $P$
notons $\Omega(\mathfrak{a}_{1},P)$ l'ensemble 
de $s \in \bigcup_{Q\supseteq P_{0}}
\Omega(\mathfrak{a}_{1},\mathfrak{a}_{Q})$
tels que $\mathfrak{a}_{P} \subseteq s\mathfrak{a}_{1}$ et 
$s^{-1}\al$ est une combinaison linéaire 
\`a coefficients positifs d'éléments de 
$\Delta_{1}$ pour 
tout $\al \in \Delta_{Q}^{P}$. 
Si l'on pose alors 
\begin{equation}\label{eq:defFMP}
\calF(M_{1}, P) = \{ Q \in \calF (M_{1})| \exists \gamma \in U(\rmF)| 
\gamma Q \gamma^{-1} = P\}
\end{equation}
on a la décomposition en parties disjointes 
$\calF(M_{1}) = \coprod_{P \sps P_{0}} \calF(M_{1},P)$ ainsi que 
 la bijection entre $\Omega(\mathfrak{a}_{1},P)$ et $ \calF(M_{1},P)$ 
donnée par 
\begin{equation}\label{eq:semiStConjBij}
\Omega(\mathfrak{a}_{1},P) \ni s \mapsto s^{-1}P \in \calF(M_{1},P)
\end{equation}
où l'on note $s^{-1}P := w_{s}^{-1}Pw_{s}$.

Pour tous $P \in \calF(M_{0})$, 
$s \in \Omega$, $T \in \all_{0}$ et $x \in U(\A)$ 
notons la formule:
\begin{equation}\label{eq:htauSemist}
\htau_{P}(H_{P}(w_{s}x) - T) = 
\htau_{s^{-1}P}(H_{s^{-1}P}(x) - s^{-1}T - H_{s^{-1}P}( w_s^{-1})).
\end{equation}

Notons finalement, que parfois, pour économiser l'espace,
 on utilisera la notation $[G]$ pour noter $G(\rmF) \bsl G(\A)$.

\subsection{Généralités sur le groupe unitaire}\label{generalities}

Soit $\rmE$ une extension quadratique de $\rmF$ 
avec $\sigma$ le générateur du groupe de Galois $\Gal(\rmE/\rmF)$. 
Soit
$W$ un $\rmE$-espace vectoriel de dimension finie $n+1$, 
où $n \in \N$, 
muni d'une forme hermitienne non-dégénérée $\tlPhi$, 
et notons $\tlU = U(W,\tlPhi)$ le groupe unitaire associé. 
Choisissons $e_{0} \in W$ tel que $\nu_{0} := \tlPhi(e_{0},e_{0}) \neq 0$, 
notons $D_{0}$ la droite engendrée par $e_{0}$ 
et soit $V$ le sous-espace de $W$ orthogonal à $D_{0}$. 
La restriction de 
$\tlPhi$ à $V$, notée $\Phi$, est alors non-dégénérée et 
l'on pose $U = U(V,\Phi)$. Les groupes $\tlU$ et $U$ sont des $\rmF$-groupes 
algébriques et l'on voit $U$ comme un sous-groupe de 
$\tlU$ grâce à l'inclusion 
$V \hrar W$.


Un sous-espace $V' \subseteq V$ est appelé \textit{isotrope}
si pour tous $v_{1},v_{2} \in V'$ on a 
$\Phi(v_{1},v_{2}) = 0$. 
Soit $d_{0}$ l'indice de Witt de $\Phi$ \cad le 
maximum des dimensions de sous-espaces isotropes de $V$. 
Fixons donc $V_{d_{0}}$ et $V_{d_{0}}^{\sharp}$ deux sous-espaces isotropes de $V$ de 
dimension $d_{0}$ qui sont 
en dualité grâce à $\Phi$ et 
choisissons $d$ droites 
$D_{1},\ldots ,D_{d_{0}}$ (resp. $D_{1}^{\sharp},\ldots ,D_{d_{0}}^{\sharp}$) 
engendrant 
$V_{d_{0}}$ (resp. $V_{d_{0}}^{\sharp}$)
de façon que
$\Phi(D_{i}, D_{j}^{\sharp}) \neq \{0\}$ si et seulement si $i = j$
pour tous $1 \le i,j \le d_{0}$. 
Pour $i = 0 ,\ldots, d_{0}$ posons 
$V_{i} = \bigoplus_{k=1}^{i}D_{i}$, 
$V_{i}^{\sharp} = \bigoplus_{k=1}^{i}D_{i}^{\sharp}$ 
et $Z_{i}$ le complément orthogonal de 
$V_{i} \oplus V_{i}^{\sharp}$ dans $V$. Notons $Z = Z_{d_{0}}$ et remarquons que
la restriction de $\Phi$ à $Z$ est anisotrope. 
On a alors 
\[
Z_{i} = Z \bigobot_{k = i +1}^{d_{0}}(D_{k} \oplus D_{k}^{\sharp}), \quad 
V = V_{i}^{\bot} \oplus V_{i}^{\sharp}=
V_{i} \oplus Z_{i} \oplus V^{\sharp}_{i}, \quad i = 0, 1,\ldots d_{0}.
\]

Par \textit{drapeau isotrope} dans $V$ 
on entend une suite $V_{0}', \ldots, V_{j}'$
de sous-espaces isotropes de $V$ telle que 
$V_{0}' = \{0\}$ et $V_{i-1}'\subsetneq V_{i}'$ 
pour $i = 1,\ldots, j$. 
Les stabilisateurs des drapeaux isotropes dans 
$U$ sont précisément les sous-groupes paraboliques de $U$ définis 
sur $\rmF$. Soit $M_{0}$ le sous-groupe 
de $U$ défini par 
le fait qu'il fixe les droites $D_{i},D_{j}^{\sharp}$ pour 
tous $1 \le i,j \le d_{0}$ et $M_{0}|_{Z} = U(\Phi_{|Z})$. 
Alors $M_{0}$ est un sous-groupe de Levi minimal de $U$.
Les éléments de $\calF(M_{0})$ sont en bijection 
avec des drapeaux isotropes de type
\begin{equation}\label{drapParab}
  \{0\}  = V_{0}' 
\sbn V_{1}'  \sbn \cdots \sbn V_{j}'
\end{equation}
où $V_{i}'$ est engendré par certaines 
des droites $D_{1}, \ldots, D_{d_{0}}, D_{1}^{\sharp}, \ldots, D_{d_{0}}^{\sharp}$ 
pour $i = 0, \ldots, j$ (de façon que $V_{i}'$ soit isotrope). 
Si $P \in \calF(M_{0})$ 
stabilise le drapeau comme ci-dessus (\ref{drapParab}) 
on pose
\[
V_{P} := V_{j}', \quad Z_{P} = Z \ \bigobot\limits_{
\mathclap{
\begin{subarray}{c}
1 \le i \le d_{0} \\
(D_{i} \oplus D_{i}^{\sharp}) \cap V_{P} = \{0\}
\end{subarray}}} \ 
(D_{i} \oplus D_{i}^{\sharp}), \quad 
V_{P}^{\sharp} = 
\bigoplus\limits_{
\mathclap{
\begin{subarray}{c}
1 \le i \le d _{0}\\
D_{i} \sbs V_{P}
\end{subarray}}} D_{i}^{\sharp} 
\oplus 
\bigoplus\limits_{
\mathclap{
\begin{subarray}{c}
1 \le i \le d_{0} \\
D_{i}^{\sharp} \sbs V_{P}
\end{subarray}}} D_{i}.
\] 
On a alors $V= V_{P} \oplus Z_{P} \oplus V_{P}^{\sharp}$.

Pour se placer dans le cadre du paragraphe \ref{par:prelimstrace} 
on choisit le sous-groupe parabolique minimal $P_{0} \in \calF(M_{0})$ 
défini comme le stabilisateur du drapeau isotrope
\begin{equation}\label{drapIsoFix}
  \{0\}  = V_{0} 
\subsetneq V_{1} \subsetneq \cdots \subsetneq V_{d_{0}-1} \subsetneq V_{d_{0}}.
\end{equation}
Alors, les sous-groupes paraboliques standards de $U$  correspondent  
aux $l+1$-uplets $(i_{0},i_{1},\ldots, i_{l})$ où
$0 \le l \le d_{0}$ et $0 = i_{0} < i_{1} \ldots < i_{l}  \le d_{0}$. 
Le $l$-uplet $(i_{0},i_{1},\ldots, i_{l})$
correspond au sous-groupe de $U$ fixant le drapeau
\begin{equation}\label{drapParabStand}
\{0\}  = V_{i_{0}} 
\subseteq  \cdots \subseteq V_{i_{l}}.
\end{equation}

%

Pour un sous-groupe parabolique semi-standard $P$ dans $U$ défini 
par un drapeau comme dans (\ref{drapParab}), notons 
$\tlP$ le sous-groupe parabolique de $\tlU$ fixant 
le m\^eme drapeau sauf que 
l'on regarde les sous-espaces $V_i'$ comme 
des sous-espaces de l'espace hermitien 
$W$.

Soit $M_{\tlzero} = M_{\tlP_{0}}$ le facteur de Levi 
de $\tlP_{0}$ qui stabilise les droites 
$D_{i}, D_{j}^{\sharp}$ 
pour $1 \le i,j \le d_{0}$. 
Pour $P \in \calF(M_{0})$ notons 
$M_{\tlP}$ le facteur de Levi de $\tlP$ 
contenant $M_{\tlzero}$. 
L'application $P \mapsto \tlP$ est  une bijection 
entre $\calF(M_{0})$ et $\calF(M_{\tlzero})$. 
En plus, sa restriction à $\calP(M_{Q})$, où 
$Q \in \calF(M_{0})$ est fixé, est une bijection 
entre $\calP(M_{Q})$ et $\calP(M_{\tlQ})$. On a pourtant:

\brem\label{tlPisNotMiniRem}
 En général, le $\rmF$-sous-groupe parabolique 
$\tlP_{0}\subseteq \tlU$ 
n'est pas minimal. 
Il l'est si et seulement s'il 
n'y a pas de $v \in Z$ tel 
que $\Phi(v,v)=- \nu_{0}$.
\erem

Posons $\tlul = \Lie(\tlU)$ et 
$\ul = \Lie(U) \hrar \tlul$.
On définit la forme bilinéaire  
$\langle \cdot \ , \cdot\rangle$ 
sur $\tlul$ invariante par adjonction donnée par
\begin{equation}\label{eq:kiliform}
\langle \cdot \ , \cdot\rangle 
: \tlul \times \tlul \rightarrow \Ga, 
\quad \langle X, Y \rangle = \Tr(XY).
\end{equation}

Soit $P \in \calF(M_{0})$. 
Notons $\ml_{P} = \Lie (M_{P})$, $\ml_{\tlP} = \Lie(M_{\tlP})$, 
$\nl_{P} = \Lie (N_{P})$ et $\nl_{\tlP} = \Lie(N_{\tilde P})$.
Il existe un unique $\brP \in \calP(M_{P})$, appelé le 
sous-groupe parabolique opposé à $P$, 
tel que $\brP \cap P = M_{P}$.
Notons alors 
$\bar \nl_{P} := \nl_{\brP}$ et 
$\bar \nl_{\tlP} := \nl_{\widetilde{\brP}}$. 
La restriction de la forme $\bilif$ à 
$\bar \nl_{\tlP} \times \nl_{\tlP}$ 
(resp. $\bar \nl_{P} \times \nl_{P}$) 
est non-dégénérée donc 
l'espace $\bar \nl_{\tlP}$ (resp. $\bar \nl_{P}$)
s'identifie à l'espace dual à $\nl_{\tlP}$ 
(resp. $\nl_{P}$) grâce à cette forme.

Pour des sous-groupes paraboliques semi-standards
$P_{1} \subseteq P_{2}$ on définit aussi 
$N_{P_{1}}^{P_{2}} := M_{P_{2}} \cap N_{P_{1}}$, 
$\nl_{P_{1}}^{P_{2}} := \Lie N_{P_{1}}^{P_{2}}$, 
$\nl_{\tlP_{1}}^{\tlP_{2}} := \Lie (M_{\tlP_{2}} \cap N_{\tlP_{1}})$,  
$\bar \nl_{P_{1}}^{P_{2}} := 
\ml_{P_{2}} \cap \bar \nl_{P_{1}}$ et 
$\bar \nl_{\tlP_{1}}^{\tlP_{2}} := 
\ml_{\tlP_{2}} \cap \bar \nl_{\tlP_{1}}$.
 On a alors 
$\ml_{P_{2}} = \bar \nl_{P_{1}}^{P_{2}} \oplus \ml_{P_{1}} 
\oplus \nl_{P_{1}}^{P_{2}}$, etc.

Lorsqu'il n'y aura pas d'ambigu\"ité
 on écrira les objets associés à $\tlU$ avec le tilde: 
donc par exemple, on notera $\ml_{\tlone}$ au lieu 
de $\ml_{\tlP_{1}}$, $\nl_{\tlone}^{\tltwo}$ au 
lieu de $\nl_{\tlP_{1}}^{\tlP_{2}}$ etc. En particulier les objets 
associé au groupe $\tlP_{0}$ seront noté avec l'indice 
$\tlzero$.

\subsection{Les invariants}\label{lesInvs}

Pour un $\rmE$-espace vectoriel $\calV$ soit 
$\gll_{\rmE/\rmF}(\calV)$ le $\rmF$-groupe algébrique $\Res_{\rmE/\rmF}(\End_{\rmE}(\calV))$.
Le décomposition $W = V \oplus D_{0}$ 
induit l'inclusion $\Gl(V) \hrar \Gl(W)$. 
Soit $X \in \gll_{\rmE/\rmF}(W)$, on a
\begin{equation}\label{xisamatrix}
X = 
\begin{pmatrix}
B & u \\
v & d \\
\end{pmatrix} 
\end{equation}
où $B \in \gll_{\rmE/\rmF}(V)$, $u \in \Res_{\rmE/\rmF}(\Hom_{\rmE}(D_{0},V))$, 
$v \in \Res_{\rmE/\rmF}(\Hom_{\rmE}(V,D_{0}))$ et 
$d \in \gll_{\rmE/\rmF}(D_{0})$. 
On identifie $u$ avec $u(e_{0}) \in \Res_{\rmE/\rmF}(V)$ et
$v$ avec l'élément de $\Res_{\rmE/\rmF}(V^{*})$ défini 
par $x \mapsto \Phi(e_{0},v(x))$. 
Notons $\ul_{1}$ l'algèbre de Lie du groupe unitaire 
$U(D_{0}, \tlPhi|_{D_{0}})$.
On a 
$U(D_{0}, \tlPhi|_{D_{0}}) \hrar \Res_{\rmE/\rmF}(\Gl(D_{0}))$,
$U \hrar \Res_{\rmE/\rmF}(\Gl(V))$ ainsi que 
$\tlU \hrar \Res_{\rmE/\rmF}(\Gl(W))$
d'où 
$\ul_{1} \hrar \gll_{\rmE/\rmF}(D_{0})$, 
$\ul \hrar \gll_{\rmE/\rmF}(V)$ et
$\tlul \hrar \gll_{\rmE/\rmF}(W)$.
Donc, on a 
que $X \in \gll_{\rmE/\rmF}(W)$, décomposé comme (\ref{xisamatrix}), 
est dans 
 $\tlul$ si et seulement si
\begin{equation}\label{vIsPhi}
B \in \ul \text{ et }
v = u^{\sharp} := -\Phi(u,\cdot)
\text{ et } d \in \ul_{1}.
\end{equation}

Le groupe $\Res_{\rmE/\rmF}\Gl(V)$ agit sur $\gll_{\rmE/\rmF}(W)$ 
par conjugaison.
On introduit alors les invariants suivants de cette action.
Soit $X \in \gll_{\rmE/\rmF}(W)$ décomposé comme dans 
(\ref{xisamatrix}). On pose 
$A_{0}(X) = \Phi(e_{0},d(e_{0}))$ et 
$A_{i}(X) = vB^{i-1}u  $ pour $i = 1,2,\ldots ,n$. 
Notons aussi $B_{i}(X)$ les coefficients 
de polynôme caractéristique 
de $B$:
\begin{equation*}
\det(T-B) = T^{n} - B_{1}(X)T^{n-1} + \cdots + (-1)^{n}B_{n}(X).
\end{equation*}
On a alors un $\rmF$-morphisme $\Res_{\rmE/\rmF}(\Gl(V))$-invariant
$Q : \gll_{\rmE/\rmF}(W) \rightarrow 
\Res_{\rmE/\rmF}(\A_{\rmE}^{2n+1})$, où
$\A_{\rmE}$ c'est la $\rmE$-droite affine,  donné par
\begin{displaymath}
Q(X) = (A_{0}(X), A_{1}(X),\ldots, A_{n}(X), B_{1}(X), \ldots, B_{n}(X)).
\end{displaymath} 

On dit que $X \in \gll_{\rmE/\rmF}(W)$ 
est \textit{semi-simple régulier} s'il vérifie les conditions 
de la proposition suivante, due à \cite{rallSchiff}, théorème 6.1 et 
proposition 6.3. 

\brop\label{prop:relelt}
Soit $X = \matx{B}{u}{v}{d} \in \gll_{\rmE/\rmF}(W)$,
alors les conditions 
suivantes sont équivalentes
\begin{enumerate}[1)]
\item $\det(a_{ij}) \neq 0$ o\`u $a_{ij} = vB^{i+j}u$, 
$0 \le i,j \le n-1$.
\item Le stabilisateur de $X$ dans 
$\Res_{\rmE/\rmF}(\Gl(V))$ 
est trivial et 
l'orbite de $X$ dans $\gll_{\rmE/\rmF}(W)$ 
pour l'action de $\Res_{\rmE/\rmF}(\Gl(V))$ 
est fermée pour 
la topologie de Zariski.
\end{enumerate}
\erop

Par restriction, on a l'application 
$U$-invariante $Q : \tlul \rightarrow \Res_{\rmE/\rmF}(\A_{\rmE}^{2n+1})$. 
On dit qu'un élément de $\tlul$
est semi-simple régulier s'il l'est comme un 
élément de $\gll_{\rmE/\rmF}(W)$. 
En vertu de la proposition \ref{prop:relelt}, l'ensemble des éléments semi-simples réguliers dans 
$\tlul$ est un ouvert pour la topologie de Zariski. 
La proposition suivante, démontrée dans \cite{rallSchiff}, 
proposition 17.2, accentue le rôle privilégié 
des éléments semi-simples réguliers.

\brop\label{prop:orbRegConj} 
 Soient $X $ et $ Y $ deux éléments semi-simples réguliers de 
$\tlul(\rmF)$. Alors $Q(X) = Q(Y)$ si et seulement s'ils 
sont conjugués par $U(\rmF)$. 
\erop

 On introduit une relation d'équivalence sur 
 $\tlul(\rmF)$: $X_{1} =\matx{B_{1}}{u_{1}}{u_{1}^{\sharp}}{d_{1}} ,
 X_{2}  = \matx{B_{2}}{u_{2}}{u_{2}^{\sharp}}{d_{2}} 
 \in \tlul(\rmF)$ sont 
 équivalents si et seulement s'ils sont dans la même fibre 
 de $Q$ 
 et si les parties semi-simples de 
 $B_{1}$ et $B_{2}$ sont conjugués sous $U(\rmF)$. 
 Notons 
 $\mathcal{O}$ l'ensemble de classes d'équivalence pour 
 cette relation. 
Grâce à la proposition \ref{prop:orbRegConj} ci-dessus, si $\ol \in \calO$ est 
 telle qu'il existe un $X \in \ol$ semi-simple régulier alors 
 tous les éléments de $\ol$ le sont 
 et $\ol$ est une $U(\rmF)$-classe de conjugaison dans $\tlul(\rmF)$. 

On étudiera maintenant 
les intersection des classes $\ol \in \calO$ 
avec les algèbres de 
Lie
de sous-groupes paraboliques relativement standards.

\blem\label{shortLemmeOrb}
 Soient $X = \matx{B}{u}{u^{\sharp}}{d} \in\tlul(\rmF)$ et  
 $P \in \calF(M_{0})$. Alors
\begin{enumerate}[a)]
\item 
$X \in \ml_{\tlP}(\rmF)$ si et seulement si 
$B \in \ml_{P}(\rmF)$ et 
$u \in Z_{P}$,
\item 
$X \in \nl_{\tlP}(\rmF)$ si et seulement si 
$B \in \nl_{P}(\rmF)$ et $u \in V_{P}$.
\end{enumerate}

\bdem 
\textit{a)}. 
Supposons $X \in \ml_{\tlP}(\rmF)$. 
Si l'on note $\tlV_{P}^{\bot}$ l'espace orthogonal \`a $V_{P}$ 
dans $W$ on a $\tlV_{P}^{\bot} = V_{P}^{\bot}\oplus D_{0}$. 
Donc $Xe_{0} = u+de_{0} \in Z_{P} \oplus D_{0}$ 
d'où $u \in Z_{P}$. 
On voit alors que la restriction de $X$ à $V_{P}$ 
et à $V_{P}^{\sharp}$
égale $B$, 
d'où 
$B \in \ml_{P}(\rmF)$.
La réciproque est évidente. 
L'analyse analogue démontre \textit{b)}. 
\edem
\elem

\brop\label{classIntUnip}
 Soit $P \in \calF(M_{0})$. 
Alors, pour tous $X \in \ml_{\tlP}(\rmF)$ et $N \in \nl_{\tlP}(\rmF)$ on a
\begin{displaymath}
Q(X) = Q(X+N).
\end{displaymath}

\bdem 
Soient alors 
$X = \matx{B}{u}{u^{\sharp}}{d} \in \ml_{\tlP}(\rmF)$ 
et $N = \matx{B_{N}}{u_{N}}{u_{N}^{\sharp}}{d_{N}} \in \nl_{\tlP}(\rmF)$.
D'après le lemme \ref{shortLemmeOrb} on a
$B \in \ml_{P}(\rmF), B_{N} \in \nl_{P}(\rmF), 
u \in Z_{P}, u_{N} \in 
V_{P}$. 
Soit $V_{P}^{\bot} = V_{P} \oplus Z_{P}$ l'espace orthogonal à $V_{P}$ dans $V$.
En utilisant le fait que 
$[B,B_{N}] = BB_{N} - B_{N}B \in \nl_{P}(\rmF)$ on voit 
qu'on a 
\begin{equation*}
((B+B_{N})^{i} - B^{i})V_{P}^{\bot}
\subseteq V_{P}, \ \forall \ i \ge 0.
\end{equation*}
En utilisant cela, le fait que 
$u + u_{N} \in V_{P}^{\bot}$ ainsi que 
la relation (\ref{vIsPhi}), 
on obtient
\begin{multline*}
A_{i+1}(X+N) = \Phi(u+u_{N}, (B+B_{N})^{i}(u+u_{N})) = 
\Phi(u+u_{N}, ((B+B_{N})^{i}-B^{i})(u+u_{N})) \\
+\Phi(u+u_{N}, B^{i}(u+u_{N})) = 0 + 
\Phi(u_{N}, B^{i}u_{N}) + 
\Phi(u,B^{i}u_{N}) + 
\Phi(u_{N}, B^{i}u) + 
\Phi(u, B^{i}u) = \\
=
\Phi(u,B^{i}u) = 
A_{i+1}(X)
\end{multline*} 
pour $i \ge 0$. L'égalité $A_{0}(X+N) = A_{0}(X)$ est claire.

Quant aux invariants $B_{i}$ , 
l'identité $B_{i}(X+N) = B_{i}(X)$ est un résultat
classique de l'algèbre linéaire.
\edem
\erop

\bcor\label{invsIntPar}
 Soient $P \in \calF(M_{0})$ et 
$\mathfrak{o} \in \mathcal{O}$.
Pour tous $R \subseteq \ml_{\tlP}(\rmF)$ et 
$S \subseteq \nl_{\tlP}(\rmF)$ on a
\begin{displaymath}
\mathfrak{o} \cap (R \oplus S) = (\mathfrak{o} \cap R) \oplus S.
\end{displaymath}
\ecor

\subsection{Racines simples}\label{racSimples}

Pour tout 
choix de vecteurs non nuls dans $D_{i},D_{j}^{\sharp}$ et pour toute base de 
$Z$, l'ensemble des $\rmF$-points de $A_{0} = A_{P_{0}}$ égale
\begin{multline}
\{\diag(x_{1},\ldots, x_{d_{0}}) :=  
\begin{pmatrix}
x_{1} &  &  &  &  & &  \\
 & \ddots &  &  &  & &  \\
 &  & x_{d_{0}} &  &  &  &  \\
 & & & \Id_{Z} & & &  \\
& & & & x_{d_{0}}^{-1} &  & \\
& & & & & \ddots &  \\
& &  & & & & x_{1}^{-1} \\
\end{pmatrix} | \  x_{1}, \ldots, x_{n} \in \rmF^{*}\}.
\end{multline}

Soit $P \in \calF(M_{0})$. Remarquons que
$A_{P} = A_{\tlP}$ (on voit toujours $U$ comme un sous-groupe de 
$\tlU$).
On va alors écrire $A_{P}$ tout court pour 
dénoter l'un d'eux.

Soient $P_{1}, P_{2} \in \calF(M_{0})$ standards. 
On s'intéresse 
\`a la relation entre $\Delta_{1}^{2}$ 
et $\Delta_{\tlone}^{\tltwo}$. 
Supposons donc d'abord que $P_{1} = P_{0}$ et 
$P_{2} = U$. On voit 
alors que $\Delta_{0}$ 
égale $\{\al_{1}, \al_{2}, \ldots, \al_{d_{0}}\}$ où
$\al_{i}(\diag(x_{1},\ldots, x_{d_{0}})) = x_{i}/x_{i+1}$ pour 
$1 \le i < d_{0}$ et 
\begin{displaymath}
\al_{d}(\diag(x_{1},\ldots, x_{d_{0}})) = 
\begin{cases} x_{d_{0}}^{2} & \text{si }2d_{0} = \dim V, \\
x_{d_{0}} & \text{sinon.} \\
\end{cases}
\end{displaymath}
Par définition de l'inclusion 
$U \hookrightarrow \tilde U$ on voit que 
$ \Delta_{\tlzero} = 
\{\tilde \al_{1}, \tilde \al_{2}, \ldots, \tilde \al_{d_{0}}\}$ 
o\`u $\tilde \al_{i} = \al_{i}$ pour $1 \le i < d_{0}$ et 
$\tilde \al_{d_{0}} = 
\displaystyle \frac{\al_{d_{0}}}{2}$ 
si $2d_{0} = \dim V$ et 
$\tilde \al_{d_{0}} = \al_{d_{0}}$ sinon.
On obtient donc le lemme suivant:

\blem\label{lem:DeltaSontLesMemes}
 Pour tous sous-groupes
paraboliques standards $P_{1}\subseteq P_{2}$, 
les éléments de $\Delta_{\tlone}^{\tltwo}$ sont 
égaux aux éléments de $\Delta_{1}^{2}$ 
\`a la multiplication par $1/2$ pr\`es.

\bdem 
Cela découle du fait que pour un sous-groupe 
parabolique $P$, les éléments de 
$\Delta_{P}$ (resp. $\tilde \Delta_{P}$)
sont des restrictions de $\Delta_{0}\smallsetminus 
\Delta_{0}^{P}$ 
(resp. $\Delta_{\tlzero} \smallsetminus 
\Delta_{\tlzero}^{\tlP}$) \`a $A_{P}$. 
\edem
\elem

\subsection{Les mesures de Haar}\label{par:haarMesSubs}

Fixons $dx$ une mesure de Haar sur $U(\A)$. Soit 
$P = M_{P}N_{P}$ un sous-groupe parabolique semi-standard de 
$U$. On fixe alors, pour tout sous-groupe connexe $V$ de $N_{P}$ 
(resp. toute sous-alg\`ebre $\mathfrak{h}$ de $\nl_{\tlP}$)
l'unique mesure de Haar sur $V(\A)$ (resp. $\mathfrak{h}(\A)$)
pour laquelle le volume de $V(\rmF)\backslash V(\A)$ 
(resp. $\mathfrak{h}(\rmF)\backslash \mathfrak{h}(\A)$)
soit $1$. Choisissons la mesure de Haar sur $\rmK$ 
normalisé de m\^eme façon. On fixe aussi une norme 
euclidienne $\|\cdot \|$ $\Omega$-invariante 
sur $\mathfrak{a}_{0}$
et des
mesures de Haar sur tous les sous-espaces 
de $\mathfrak{a}_{0}$ 
compatibles avec cette norme. On en déduit
la mesure de Haar sur $A_{P}^{\infty}$ 
grâce à l'isomorphisme
$H_{P}$.  

Soient $dp$ la mesure de Haar sur $P(\A)$ invariante 
\`a gauche normalisé de façon que $dx = dpdk$ 
(grâce à la décomposition d'Iwasawa).
Notons $\rho_{P}$ (resp. $\rho_{\tlP}$)
l'élément de $\all_{P}^{*}$ 
tel que $d(\Ad(m)n) = e^{2\rho_{P}(H_{P}(m))}dn$ 
pour $m \in M_{P}(\A)$ et $n \in N_{P}(\A)$ 
(resp. $d(\Ad(m)U_{\tlP}) = e^{2\rho_{\tlP}(H_{P}(m))}dU_{\tlP}$ 
pour $m \in M_{P}(\A)$ et $U_{\tlP} \in \nl_{\tlP}(\A)$)
Il s'ensuit 
qu'il existe une unique 
mesure de Haar $dm$ sur 
$M_{P}(\A)$ telle 
que si l'on écrit 
$p = nm$ 
où $p \in P(\A)$, $n \in N_{P}(A)$ et
$m \in M_{P}(\A)$ 
alors $dp = e^{-2\rho_{P}(H_{P}(m))}dndm$. 
On pose aussi 
\[
\upla_{P} := 2(\rho_{\tlP} - \rho_{P}).
\] 
Puisque $\rho_{P}$ (resp. $\rho_{\tlP}$) 
est la demi-somme des poids (avec multiplicités) pour l'action de $A_{P}$ sur $\nl_{P}$
(resp. $\nl_{\tlP}$), on a,
grâce à la partie \textit{b)} du lemme \ref{shortLemmeOrb}, que 
$\upla_{P}$ 
c'est juste la somme de poids (avec multiplicités) 
pour l'action du tore $A_{P}$ 
sur $V_{P}$.

\subsection{Fonctions de Bruhat-Schwartz}\label{bschwartz}

On note  
$\A^{\infty}$ l'anneau des ad\`eles finis de $\rmF$ et 
$\A_{\infty}$ le produit de
complétions de $\rmF$ aux
places infinies, de façon qu'on a 
$\A = \A_{\infty}\times \A^{\infty}$. On fixe 
une norme $\| \cdot \|$ sur 
$\tlul(\A_{\infty})$. 
Si $X \in \tlul(\A)$, par $\| X \|$ on entend 
$\|\cdot \|$ appliqué \`a la projection de $X$ sur sa partie 
infinie grâce \`a la décomposition canonique 
$\tlul(\A) = \tlul(\A_{\infty})\times \tlul(\A^{\infty})$. 

Soient $\calV_{i}(\rmF) \subseteq \tlul(\rmF)$ 
des sous-$\rmF$-espaces en somme directe 
et notons 
$\calV_{i}(\A) = \calV_{i}(\rmF)\otimes_{\rmF}\A$
pour $1 \le i \le k$.
Notons 
$\calV(\rmF) = \bigoplus_{1 \le i \le k} V_{i}(\rmF)$ 
et $\calV(\A) = \calV(\rmF)\otimes_{\rmF}\A$. 
Donc $\calV(\A)$ et les $\calV_{i}(\A)$ sont 
naturellement des sous-$\A$-modules de 
$\tlul(\A)$.
Pour $X \in \calV(\A)$ notons $X_{i}$ sa projection 
sur $\calV_{i}(\A)$. On a alors le 
lemme suivant qui se démontre en utilisant 
l'équivalence des normes sur 
$\tlul(\A_{\infty})$. 

\blem\label{AGinequality}
 Il existe un réel positif $c_{\calV}$ tel que 
pour tout $X \in \calV(\A)$ on a
\begin{displaymath}
\| X\| \ge c_{\calV} \| X_{1}\|^{\frac{1}{k}} \cdots \| X_{k}\|^{\frac{1}{k}}.
\end{displaymath}
\elem

Notons $\mathcal{S}(\tlul(\A))$ l'ensemble des 
fonctions des Bruhat-Schwartz sur $\tlul(\A)$ 
\textit{i.e.} l'espace de fonctions sur 
$\tlul(\A)$ engendré par des fonctions du type 
$f_{\infty}\otimes \chi^{\infty}$ o\`u 
$f_{\infty}$ est une fonction de la classe de 
Schwartz sur $\tlul(\A_{\infty})$ et 
$\chi^{\infty}$ est une fonction caractéristique d'un 
compact ouvert de $\tlul(\A^{\infty})$. 
Un opérateur différentiel $\partial$ sur 
$\tlul(\A_{\infty})$ s'étend
sur $\tlul(\A)$ en posant 
$\partial (f_{\infty}\otimes \chi^{\infty}) = 
\partial f_{\infty} \otimes \chi^{\infty}$. 
On a donc que pour tout 
$f \in \mathcal{S}(\tlul(\A))$, 
tout opérateur différentiel $\partial$ sur 
$\tlul(\A_{\infty})$ et
tout $n \in \N$ il existe une constante 
$C_{f,\partial, n}$ telle 
que
\begin{displaymath}
\sup_{X \in \tlul(\A)}\|X\|^{n}|\partial f(X)| \le 
C_{f,\partial,n}.
\end{displaymath}
\section{Convergence du noyau}\label{sec:convergence}

À partir de cette section, jusqu'à la 
section \ref{sec:orbRrssS}, par un sous-groupe parabolique de $U$ 
on entend un sous-groupe parabolique standard.

Soit $f \in \mathcal{S}(\tlul(\A))$. Pour 
un sous-groupe parabolique $P$ de $U$ et 
une classe $\mathfrak{o} \in \mathcal{O}$ notons pour 
$x \in M_{P}(\rmF)N_{P}(\A)\backslash U(\A)$
\begin{equation}\label{eq:kPolDef}
k_{P,\mathfrak{o}}(x) = k_{f, P,\mathfrak{o}}(x)=  \sum_{\xi \in \ml_{\tlP}(\rmF) \cap \mathfrak{o}}
\int_{\nl_{\tlP}(\A)}f(x\inv (\xi + U_{P})x)dU_{P}
\end{equation}
et pour $x \in U(\rmF)\backslash U(\A)$
\begin{displaymath}
k^{T}_{\mathfrak{o}}(x) = k^{T}_{f, \mathfrak{o}}(x) =  
\sum_{P\supseteq P_{0}} 
(-1)^{d_{P}}\sum_{\delta \in P(\rmF)\backslash U(\rmF)}
\htau_{P}(H_{P}(\delta x) - T)k_{P,\mathfrak{o}}(\delta x)
\end{displaymath}
où $d_{P}$ est défini par (\ref{eq:dPDef}). N.B. la somme sur $\delta$ 
dans $P(\rmF)\backslash U(\rmF)$ est finie en vertu du lemme 5.1 de \cite{arthur3}.

Le but de ce chapitre est de démontrer le résultat suivant:

\begin{theo}\label{thm:MainConv} Soit $T_{+} \in \all_{0}^{+}$ 
introduit dans le paragraphe \ref{par:prelimstrace}. Alors, 
pour tout
 $T \in T_{+} + \mathfrak{a}_{0}^{+}$ on a
\begin{displaymath}
\sum_{\mathfrak{o} \in \mathcal{O}}
\int_{U(\rmF) \backslash U(\A)}|k_{\mathfrak{o}}^{T}(x)|dx < \infty.
\end{displaymath}

\bdem
Suivant Arthur \cite{arthur3} paragraphe 6, 
pour $P_{1}\subseteq P_{2}$
on définie la fonction $\sigma_{1}^{2}$ par
\begin{equation*}
\sigma_1^2(H) = \sigma_{P_{1}}^{P_{2}}(H)  = 
(-1)^{d_{2}}
\sum_{Q \supseteq P_{2}}
(-1)^{d_{Q}}\tau_{P_{1}}^{Q}(H)\htau_{Q}(H),\ 
H \in \mathfrak{a}_{{1}}.
 \end{equation*}
 C'est la fonction caractéristique d'un sous-ensemble de $\all_{1}$. 
 On pose aussi
\begin{align}
&\chi^{T}_{P_{1},{P_2}}(x)  = \chi^{T}_{{1},{2}}(x) = F^{{1}}(x,T)
\sigma_{{1}}^{{2}}(H_{{1}}(x)-T) \label{chipp}, \ 
x \in P_{1}(\rmF)\backslash U(\A),
\\
&k_{P_{1},{P_2},\mathfrak{o}}(x) =k_{{1},{2},
\mathfrak{o}}(x) = \sum_{P_{1} \subseteq P \subseteq P_{2}}
(-1)^{d_{P}}k_{P,\mathfrak{o}}(x), \ 
x \in P_{1}(\rmF)\backslash U(\A).
\end{align}
Alors, en utilisant le lemme 6.4 de \cite{arthur3} et 
l'invariance de $k_{P,\mathfrak{o}}$ et $H_{P}$ par rapport 
\`a $P(\rmF)$ on 
s'aperçoit qu'on a l'identité
\begin{equation}\label{chgmtordre}
k_{\mathfrak{o}}^{T}(x) = \sum_{P_{1} \subseteq P_{2}} 
\sum_{\delta \in P_{1}(\rmF) \backslash U(\rmF)}
\chi_{{1},{2}}^{T}(\delta x)k_{{1},{2},\mathfrak{o}}(\delta x).
\end{equation}
On a donc
\begin{displaymath}
\sum_{\mathfrak{o} \in \mathcal{O}}
\int_{U(\rmF)\backslash U(\A)}|k_{\mathfrak{o}}^{T}(x)|dx \le 
\sum_{\mathfrak{o} \in \mathcal{O}}
\sum_{P_{1} \subseteq P_{2}} 
\int_{P_{1}(\rmF)\backslash U_{}(\A)}\chi_{1,2}^{T}(x)
|k_{1,2,\mathfrak{o}}(x)|dx.
\end{displaymath}
Il suffit donc de  montrer que
pour tous sous-groupes paraboliques $P_{1} \subseteq P_{2}$
\begin{equation}\label{comesdownto}
\sum_{\mathfrak{o} \in \mathcal{O}}
\int_{P_{1}(\rmF)\backslash U(\A)}
\chi_{{1},{2}}^{T}(x)|k_{{1},{2},\mathfrak{o}}(x)|dx 
< \infty.
\end{equation}
Il résulte du lemme 6.1 de \cite{arthur3} 
que si
$P_{1} = P_{2} \neq U$ on a 
$\sigma_{{1}}^{{2}} \equiv 0$ donc l'intégrale 
(\ref{comesdownto}) vaut $0$ dans ce cas. D'autre part, quand 
$P_{1} = P_{2} = U$ on a que $F^{U}$ est la 
fonction caractéristique d'un compact dans $U(\rmF) \bsl U(\A)$ donc 
l'intégrale est bien finie. On peut supposer alors que 
$P_{1} \neq P_{2}$. Dans ce cas on montrera quelque chose de 
plus fort que (\ref{comesdownto}).

\begin{theo}\label{thm:MainConv2} 
Soient $f \in \mathcal{S}(\tlul(\A))$ et 
$P_{1}, P_{2}$ deux sous-groupes 
paraboliques 
tels que $P_{1}\subsetneq P_{2}$. Alors, pour tout réel $\varepsilon_{0} > 0$ et tout $N \in \N$ 
il existe une 
constante $C$ qui ne dépend que de $N$, $f$ et 
$\varepsilon_{0}$ telle que
\begin{displaymath}
\sum_{\mathfrak{o} \in \mathcal{O}}
\int_{P_{1}(\rmF)\backslash U(\A)}
\chi_{{1},{2}}^{T}(x)|k_{{1},{2},\mathfrak{o}}(x)|dx 
< Ce^{-N\|T\|}
\end{displaymath}
pour tout $T \in T_{+} + \mathfrak{a}_{{0}}^{+}$ 
tel que $\al(T) > 
\varepsilon_{0}\|T\|$ pour tout $\al \in \Delta_{0}$.
\end{theo}
\nident
Notons au passage le corollaire immédiat de ce théorème: 
\bcor\label{almostPolCor}
 Soit $k(x) = \sum_{\mathfrak{o} \in 
\mathcal{O}}k_{U,\mathfrak{o}}(x)$ et 
$k^{T}(x) = \sum_{\mathfrak{o} \in \mathcal{O}}k^{T}_{\mathfrak{o}}(x)$.
Pour tout réel $\varepsilon_{0} > 0$ et tout $N \in \N$ 
il existe une 
constante $C$ qui ne dépend que de $N$, $f$ et 
$\varepsilon_{0}$ telle que
\begin{equation}
\int_{U(\rmF)\backslash U(\A)}
|F^{U}(x,T)k(x) - k^{T}(x)|dx 
< Ce^{-N\|T\|}
\end{equation}
pour tout $T \in T_{+} + \mathfrak{a}_{{0}}^{+}$ 
tel que $\al(T) > 
\varepsilon_{0}\|T\|$ pour tout $\al \in \Delta_{0}$.
\ecor

On introduit d'abord quelques notations.
Pour tous sous-groupes paraboliques $Q \subseteq  S$ posons
\begin{displaymath}
(\bar\nl_{\tlQ}^{\tlS})' = \bar \nl_{\tlQ}^{\tlS}
\smallsetminus \bigcup_{ Q \subseteq R \subsetneq S}
\bar \nl_{\tlQ}^{\tlR}.
\end{displaymath}
Alors $(\bar\nl_{\tlQ}^{\tlS})'$ est ouvert dans 
$\bar \nl_{\tlQ}^{\tlS}$ et on a
une décomposition en parties localement fermées
\begin{equation}\label{locFermLin}
\bar \nl_{\tlQ}^{\tlS} = \coprod_{Q\subseteq R\subseteq S}
(\bar\nl_{\tlQ}^{\tlR})'.
\end{equation}
En supposant toujours que $Q \subseteq S$ notons aussi
\begin{displaymath}
\ml_{\tlS,\tlQ}' = \ml_{\tlS} \smallsetminus 
\bigcup_{ Q \subseteq R \subsetneq S} \Lie(\tlR) =
(\bar \nl_{\tlQ}^{\tlS})' \oplus \ml_{\tlQ} \oplus \nl_{\tlQ}^{\tlS}.
\end{displaymath}
Fixons les sous-groupes paraboliques 
$P_{1} \subsetneq P_{2}$. 
Soit $P$ un sous-groupe parabolique tel que 
$P_{1}\subseteq P \subseteq P_{2}$.  
On a donc
\begin{displaymath}
\ml_{\tlP} = \coprod_{P_{1} \subseteq S \subseteq P} 
(\ml_{\tlS,\tlone}'\oplus \nl_{\tlS}^{\tlP})
\end{displaymath}
et, en utilisant le corollaire \ref{invsIntPar}, on obtient, 
pour tout $\mathfrak{o} \in \mathcal{O}$
\begin{equation*}
\mathfrak{o} \cap \ml_{\tlP}(\rmF) = 
\coprod_{P_{1} \subseteq S \subseteq P} 
\mathfrak{o} \cap 
(\ml_{\tlS,\tlone}'\oplus \nl_{\tlS}^{\tlP})(\rmF) = 
\coprod_{P_{1} \subseteq S \subseteq P} 
(\mathfrak{o} \cap \ml_{\tlS,\tlone}'(\rmF))
\oplus \nl_{\tlS}^{\tlP}(\rmF).
\end{equation*}
 Grâce à cela, 
on peut réécrire $k_{P,\mathfrak{o}}(x)$ comme
\begin{displaymath}
\sum_{P_{1} \subseteq S \subseteq P}
\sum_{\eta \in \nl_{\tlS}^{\tlP}(\rmF)}
\sum_{\zeta \in \ml_{\tlS,\tlone}'(\rmF) \cap \mathfrak{o}}
\int_{\nl_{\tlP}(\A)}f(x\inv (\eta + \zeta+ U_{P})x)dU_{P}.
\end{displaymath}
Fixons un caractère additif non-trivial $\psi$ sur 
$\rmF \bsl \A$. 
En appliquant la formule sommatoire 
de Poisson pour la somme portant sur
$\eta \in \nl_{\tlS}^{\tlP}(\rmF)$ 
de la fonction
\begin{displaymath}
\nl_{\tlS}^{\tlP}(\A) \ni Y \longmapsto
\int_{\nl_{\tlP}(\A)}
f(x\inv (Y + \zeta  + U_{P})x)dU_{P}, 
\end{displaymath}
pour tout $\zeta \in \ml_{\tlS,\tlone}'(\rmF) \cap \mathfrak{o}$,
on obtient
\begin{displaymath}
k_{P,\mathfrak{o}}(x) = 
\sum_{P_{1} \subseteq S \subseteq P}
\sum_{\eta \in \bar \nl_{\tlS}^{\tlP}(\rmF)}
\sum_{\zeta \in \ml_{\tlS,\tlone}'(\rmF) \cap \mathfrak{o}}
\Phi_{S}(x,\zeta,\eta),
\end{displaymath}
o\`u
\begin{displaymath}
\Phi_S(x,X,Y) =
\int\limits_{\mathrlap{\nl_{\tlS}(\A)}} \ 
f(x\inv (X + U_{S})x)
\psi (\langle U_{S}, Y\rangle)dU_{S}, \  x \in U(\A), \ 
X \in \ml_{\tlS}(\A), \
Y \in \bar \nl_{\tlS}^{\tltwo}(\A).
\end{displaymath}
En utilisant l'égalité (\ref{locFermLin}) 
on peut écrire $k_{P,\mathfrak{o}}(x)$ aussi comme
\begin{displaymath}
\sum_{P_{1} \subseteq S \subseteq R \subseteq  P}
\sum_{\zeta \in \ml_{\tlS,\tlone}'(\rmF) \cap \mathfrak{o}}
\sum_{\eta \in (\bar \nl_{\tlS}^{\tlR})'(\rmF)}
\Phi_{S}(x,\zeta,\eta).
\end{displaymath}
Grâce à cette formule, on a pour tout $\mathfrak{o} \in \mathcal{O}$
\begin{equation}\label{doubleNoyauLong}
\begin{split}
k_{1,2,\mathfrak{o}}(x) & =  
\sum_{P_{1} \subseteq P \subseteq P_{2}}
(-1)^{d_{P}}
k_{P,\mathfrak{o}}(x) \\ 
& = 
\sum_{P_{1} \subseteq S \subseteq R \subseteq   P \subseteq P_{2}}
(-1)^{d_{P}} 
\sum_{\zeta \in \ml_{\tlS,\tlone}'(\rmF) \cap \mathfrak{o}}
\sum_{\eta \in (\bar\nl_{\tlS}^{\tlR})'(\rmF)} 
\Phi_{S}(x,\zeta,\eta) \\ & = 
\sum_{P_{1} \subseteq S \subseteq R \subseteq   P_{2}}
\sum_{\zeta \in \ml_{\tlS,\tlone}'(\rmF) \cap \mathfrak{o}}
\sum_{\eta \in( \bar\nl_{\tlS}^{\tlR})'(\rmF)}
\Phi_{S}(x,\zeta,\eta)
\sum_{R \subseteq P \subseteq P_{2}}
(-1)^{d_{P}}.
\end{split}
\end{equation}
On invoque maintenant l'identité due \`a Arthur \cite{arthur3}, 
proposition 1.1
\begin{equation}\label{basicidentity}
\sum_{\{ P | R \subseteq P \subseteq P_{2}\}}
(-1)^{d_{P}-d_{2}} = 
\begin{cases}
0 & \text{si } R \neq P_{2}, \\
1 & \text{sinon}. \\
\end{cases}
\end{equation}
On en déduit que la somme (\ref{doubleNoyauLong}) décrivant 
$k_{1,2,\mathfrak{o}}(x)$ se réduit \`a
\begin{displaymath}
(-1)^{d_{2}}
\sum_{P_{1} \subseteq S \subseteq P_{2}} 
\sum_{\eta \in (\bar\nl_{\tlS}^{\tltwo})'(\rmF)}
\sum_{\zeta \in \ml_{\tlS,\tlone}'(\rmF) \cap \mathfrak{o}}
\Phi_{S}(x,\zeta,\eta).
\end{displaymath}
Remarquons que pour un $S$ entre $P_{1}$ et $P_{2}$ fixé, le terme correspondant
dans la somme ci-dessus est $P_{1}(\rmF)$-invariant. 
Ainsi, pour démontrer le théorème \ref{thm:MainConv2} il suffit 
de majorer
\begin{equation}\label{mainThmConv3}
\int_{P_{1}(\rmF)\backslash U(\A)}
\chi_{1,2}^{T}(x)
\sum_{\eta \in (\bar\nl_{\tlS}^{\tltwo})'(\rmF)}
\sum_{\zeta \in \ml_{\tlS,\tlone}'(\rmF)}
|\Phi_{S}(x,\zeta,\eta)|dx
\end{equation}
o\`u $P_{1} \subseteq S\subseteq P_{2}$ sont fixés.  
Remarquons que la double somme sur $\mathfrak{o} \in \mathcal{O}$ 
et $\zeta \in 
 \ml_{\tlS,\tlone}'(\rmF) \cap \mathfrak{o}$ 
s'est réduite à la somme sur tout 
$\zeta \in  \ml_{\tlS,\tlone}'(\rmF)$.

 On a la décomposition  
\begin{displaymath}
P_{1}(\rmF)\backslash  U(\A) = 
N_{1}(\rmF)\backslash N_{1}(\A)  \times
A_{1}^{\infty} \times
M_{1}(\rmF)\backslash M_{1}(\A)^{1}
\times K.
\end{displaymath}
Suivant cette décomposition $x = n_{1}a_{1}m_{1}k$ et 
$dx = e^{-2\rho_{1}(H_{0}(a_{1}))}dn_{1}da_{1}dm_{1}dk$.
Remarquons quand m\^eme que pour qu'on puisse intégrer 
sur les quotients $N_{1}(\rmF)\backslash N_{1}(\A) $ et 
$M_{1}(\rmF)\backslash M_{1}(\A)^{1}$ il faut 
d'abord intégrer sur $N_{1}(\rmF)\backslash N_{1}(\A)$ 
et puis sur $M_{1}(\rmF)\backslash M_{1}(\A)^{1}$. 
On va remplacer maintenant les intégrales sur 
$M_{1}(\rmF)\backslash M_{1}(\A)^{1}$, 
$N_{1}(\rmF)\backslash N_{1}(\A)$ et $K$ par 
un supremum sur un ensemble convenable. 
Puisque on suppose que $F^{P_{1}}(x,T) = F^1(m_{1},T)=1$
on peut restreindre l'intégrale sur 
$M_{1}(\rmF)\backslash M_{1}(\A)^{1}$ 
aux éléments $m_1$ 
ayant un représentant dans
\begin{displaymath}
\omega_{1} A_{0,1}^{\infty}(T) (K \cap M_{1}(\A)^{1}), \quad 
\text{ où }
\omega_{1} := \omega \cap M_{1}(\A)^{1},
\end{displaymath}
où le compact $\omega \subseteq N_{0}(\A)M_{0}(\A)^{1}$ est 
défini dans le paragraphe \ref{par:prelimstrace} et $A_{0,1}^{\infty}(T)$ 
est défini par (\ref{eq:A12T}) dans le même paragraphe.
%
En utilisant les faits que $A_{1}$ commute avec $M_{1}$ et 
que le volume de $M_{1}(\rmF)\backslash M_{1}(\A)^{1}$ est fini, 
on a pour 
tout $a_{1} \in A_{1}^{\infty}$
\begin{multline}\label{firstStepMajorisation}
\int\limits_{K}
\int\limits_{M_{1}(\rmF)\backslash M_{1}(\A)^{1}}
\int\limits_{[N_{1}]}
F^{1}(m_{1},T)
\sum_{\eta \in (\bar\nl_{\tlS}^{\tltwo})'(\rmF)}
\sum_{\zeta \in \ml_{\tlS,\tlone}'(\rmF)}
|\Phi_{S}(n_{1}a_{1}m_{1}k,\zeta ,\eta)|
dn_{1}dm_{1}dk
\le \\
\sup_{
\begin{subarray}{c}
k_{0} \in \omega_{1}, a_{0}^{1} 
\in A_{0,1}^{\infty}(T),\\
k \in K 
\end{subarray}}
\int\limits_{[N_{1}]}
\sum_{\eta \in (\bar\nl_{\tlS}^{\tltwo})'(\rmF)}
\sum_{\zeta \in \ml_{\tlS,\tlone}'(\rmF)}
|\Phi_{S}(n_{1}k_{0}a_{0}^{1}a_{1}k,\zeta ,\eta)|dn_{1}
\end{multline}
\`a une constante multiplicative indépendante de $T$ pr\`es.

Remarquons que la fonction $|\Phi_S(x,\cdot,\cdot)|$ 
est invariante \`a gauche 
par $N_{2}(\A)$. En effet, soient $n_{2} \in N_{2}(\A)$, 
$X \in \ml_{\tlS}(\rmA)$ et 
$Y \in \bar \nl_{\tlS}^{\tltwo}(\A)$. 
Il existe un $U_{2} \in \nl_{\tltwo}(\A)$ tel
que $n_{2}^{-1}Xn_{2} = X + U_{2}$ et on a 
\begin{multline*}
|\Phi_{S}(n_2x,X,Y)| = 
|\int_{\nl_{\tlS}(\A)}f(x^{-1}(X + n_{2}^{-1}U_{S}n_{2}+U_{2})x)
\psi (\langle U_{S}, Y\rangle)dU_{S}| = \\
|\int_{\nl_{\tlS}(\A)}
f(x^{-1}(X + n_{2}^{-1}U_{S}n_{2}+U_{2})x)
\psi (\langle 
n_{2}^{-1}U_{S}n_{2}+U_{2}, Y\rangle)dU_{S}|, \ x \in U(\A).
\end{multline*}
Dans la dernière égalité on utilise le fait 
que la forme bilinéaire $\langle \cdot, \cdot \rangle$ 
est invariante par conjugaison, que $n_{2}^{-1}Yn_{2} - Y \in \nl_{\tltwo}(\A)$ 
et $\langle \cdot, \cdot \rangle$ est triviale sur $\nl_{\tlS} \times \nl_{\tltwo}$
 et que $|\psi| \equiv 1$. 
Finalement, le changement de variable
$n_{2}^{-1}U_{S}n_{2}+U_{2} \mapsto U_{S}$ ne change pas 
la mesure de Haar sur $\nl_{\tlS}(\A)$ ce qui 
démontre que 
$|\Phi_{S}(n_2x,X,Y)| = |\Phi_{S}(x,X,Y)|$.

Fixons alors un compact $K_{1}'$ de 
$N_{1}^{2}(\A)$ qui se surjecte sur  
$N_{1}^{2}(\rmF)\backslash N_{1}^{2}(\A)$ 
et posons $K_{1} = K_{1}' \omega_{1}$. 
C'est un compact de $N_{0}^{2}(\A)M_{0}(\A)^{1}$ 
et on a
\begin{equation}\label{eq:noN2neq}
\begin{split}
\int\limits_{[N_{1}]}
\sum_{\eta \in (\bar\nl_{\tlS}^{\tltwo})'(\rmF)}
\sum_{\zeta \in \ml_{\tlS,\tlone}'(\rmF)}
|\Phi_{S}(n_{1}k_{0}a_{0}^{1}a_{1}k,\zeta ,\eta)|dn_{1} & = \\
 \int\limits_{[N_{1}^{2}]}
\sum_{\eta \in (\bar\nl_{\tlS}^{\tltwo})'(\rmF)}
\sum_{\zeta \in \ml_{\tlS,\tlone}'(\rmF)}
|\Phi_{S}(n_{1}^{2}k_{0}a_{0}^{1}a_{1}k,\zeta ,\eta)|dn_{1}^{2} & \le \\ 
\sup_{k_{1} \in K_{1}}
\sum_{\eta \in (\bar\nl_{\tlS}^{\tltwo})'(\rmF)}
\sum_{\zeta \in \ml_{\tlS,\tlone}'(\rmF)}
|\Phi_{S}(k_{1}a_{0}^{1}a_{1}k,\zeta ,\eta)|.
\end{split}
\end{equation}

Ensuite on a $\Phi_{S}(k_{1}a_{0}^{1}a_{1}k,\zeta ,\eta) = 
\Phi_{S}(a_{0}^{1}a_{1}(a_{0}^{1}a_{1})^{-1}k_{1}(a_{0}^{1}a_{1})k,\zeta ,\eta )$.
 On prétend que, tout comme $k_{1}$, l'expression 
$(a_{0}^{1}a_{1})^{-1}k_{1}(a_{0}^{1}a_{1})$  reste dans un compact 
fixé de $N_{0}^{2}(\A)M_{0}(\A)^{1}$
pour tout $a_{1} \in A_{1}^{\infty}$ tel que 
$\sigma_{1}^{2}(H_{1}(a_{1})-T)=1$ et tout $a_{0}^{1} \in A_{0,1}^{\infty}(T)$.  
On aura besoin du lemme suivant:

\blem\label{toughLinLem}
Pour tout $a_{0}^{1} \in A_{0,1}^{\infty}(T)$ et tout
$a_{1} \in A_{1}^{\infty}$ tel que 
$\sigma_{1}^{2}(H_{1}(a_{1}) - T) = 1$ on a
$\tau_{0}^{2}(H_{0}(a_{0}^{1}a_{1}) - T_{1}) = 1$. 

\bdem 
La preuve se trouve dans \cite{chaud}, paragraphe 4. 
\edem
\elem

\nident
Donc, comme $\Delta_{0}^{2}$ c'est l'ensemble de racines simples 
pour l'action du  tore $A_{0}^{\infty}$ 
sur $N_{0}^{2}(\A)$ on obtient bien que 
pour tout $a \in A_{0}^{\infty}$ tel que $\tau_{0}^{2}(H_{0}(a) - T_{1}) = 1$
l'ensemble
$a^{-1}K_{1}aK$ 
est
contenu dans un compact de $U(\A)$ 
indépendant de $T$ et $a$ qu'on 
va noter $\Gamma$. 

On voit donc que pour tous
$a_{0}^{1} \in A_{0,1}^{\infty}(T)$, 
$a_{1} \in A_{1}^{\infty}$,
$k \in K$ et $k_{1} \in K_{1}$
on a
\[
\sum_{
\begin{subarray}{c}
\eta \in (\bar\nl_{\tlS}^{\tltwo})'(\rmF) \\
\zeta \in \ml_{\tlS,\tlone}'(\rmF)
\end{subarray}}
|\Phi_{S}(a_{0}^{1}a_{1}(a_{0}^{1}a_{1})^{-1}k_{1}(a_{0}^{1}a_{1})k,\zeta ,\eta)|
\le \sup_{y \in \Gamma}
\sum_{
\begin{subarray}{c}
\eta \in (\bar\nl_{\tlS}^{\tltwo})'(\rmF) \\
\zeta \in \ml_{\tlS,\tlone}'(\rmF)
\end{subarray}}
|\Phi_{S}(a_{0}^{1}a_{1}y,\zeta ,\eta)|.
\]
Ensuite, en faisant le changement de variable 
$(a_{0}^{1}a_{1})^{-1}U_{S}(a_{0}^{1}a_{1}) \mapsto U_{S}$  
dans la définition de la fonction 
$\Phi_{S}$ on obtient
\begin{equation*}
\Phi_{S}(a_{0}^{1}a_{1}y,\zeta,\eta) = e^{2\rho_{\tlS}(H_{0}(a_{0}^{1}a_{1}))}
\Phi_{S}(y,\Ad(a_{0}^{1}a_{1})^{-1}\zeta,\Ad(a_{0}^{1}a_{1})^{-1}\eta).
\end{equation*}
En prenant compte de ceci ainsi que 
des majorations 
(\ref{firstStepMajorisation}) et
(\ref{eq:noN2neq})
on obtient que
l'expression 
(\ref{mainThmConv3}) est 
majorée par 
\begin{multline}\label{mainThmConv4}
\int\limits_{A_{1}^{\infty}} \ \,
\sup_{\mathmakebox[0.8cm]{\begin{subarray}{c}
a_{0}^{1} \in A_{0,1}^{\infty}(T) \\
y \in \Gamma
\end{subarray}
}}
e^{2(\rho_{\tlS}-\rho_{1})(H_{0}(a_{1}a_{0}^{1}))}
 \sigma_{1}^{2}(H_{1}(a_{1})-T) \\
\sum_{\eta \in (\bar\nl_{\tlS}^{\tltwo})'(\rmF)}
\sum_{\zeta \in \ml_{\tlS,\tlone}'(\rmF)}
|\Phi_{S}(y,\Ad((a_{0}^{1}a_{1})^{-1})\zeta,\Ad((a_{0}^{1}a_{1})^{-1})\eta)|
da_{1}.
\end{multline}

Fixons une $\Q$-base $\{e_{i}\}$
de $\tilde \ul(\rmF)$ composée 
de vecteurs propres pour l'action de $A_{0}$ 
sur $\tilde \ul$. Les éléments de
l'espace $\mathcal{S}(\tlul(\A))$ sont, 
par définition,  les sommes finies des fonctions de type 
$f_{\infty}\otimes \chi^{\infty}$ o\`u 
$f_{\infty}$ est une fonction de la classe de 
Schwartz sur $\tlul(\A_{\infty})$ et
$\chi^{\infty}$ est une fonction 
caractéristique d'un compact ouvert dans $\tlul(\A^{\infty})$. 
Alors, comme $\Gamma$ est compact, 
il existe un compact
$\tltl \sbs \tlul(\A^{^{\infty}})$ 
tel que $\tltl = \prod_{v}\tltl_{v}$ où $\tltl_{v}$ 
est un compact ouvert dans $\tlul(\rmF_{v})$ où $\rmF_{v}$ est 
la complétion 
de $\rmF$ par rapport à une place finie $v$, 
tel que pour 
tous
$y \in \Gamma$ et 
$X \in \ml_{\tlS}(\A)$ 
(resp. $Y \in \bar\nl_{\tlS}^{\tltwo}(\A)$)
on a $\Phi_{S}(y, X, \cdot) = 0$
(resp. $\Phi_{S}(y, \cdot, Y) =0$) 
si la projection de $X$ sur $\ml_{\tlS}(\A^{\infty})$ 
(resp. de $Y$ sur $\bar \nl_{\tlS}^{\tltwo}(\A^{\infty})$)
n'est pas dans 
$\tltl$. 
Or, les éléments de $A_{0}^{\infty}$ n'agissent 
que sur la partie infinie d'un $\zeta \in \tlul(\rmF)$ 
dans la formule $\Ad(a_{0})\zeta$ où $a_{0} \in A_{0}^{\infty}$. 
Alors, si l'on décompose $\zeta = \sum b_{i}e_{i}$  avec 
$b_{i} \in \Q$, on a que la composante correspondante 
à une place finie $v$ de $\Ad(a_{0})\zeta$ c'est juste
$\sum b_{i}e_{i}$. 
Donc, pour que cela 
soit dans un compact $\tltl_{v}$, il faut que les valuations 
$p$-adiques des coefficients $b_{i}$ soient bornées 
où $v | p$. Cela démontre que, quitte à changer la 
base $\{e_{i}\}$ par une autre de type 
$\{c_{i}e_{i}\}$ où $c_{i} \in \Q^{*}$, on peut 
supposer que les sommes dans 
(\ref{mainThmConv4})
portent en effet sur des éléments du
$\Z$-réseau engendré par les $\{e_{i}\}$ 
que l'on note $\mathcal{R}$, 
qui appartiennent également \`a 
$(\bar\nl_{\tlS}^{\tlP_{2}})'(\rmF)$ et 
$\ml_{\tlS,\tlone}'(\rmF) $.
On pose
\[
\calR_{1} = \calR \cap 
(\bar\nl_{\tlS}^{\tltwo})'(\rmF), \ 
\calR_{2} = \calR \cap ((\bar\nl_{\tlone}^{\tlS})' \oplus \bar\nl_{\tlzero}^{\tlone} \oplus \ml_{\tlzero})(\rmF), \ 
\calR_{3} = \calR \cap \nl_{\tlzero}^{\tlS}(\rmF).
\]
On majorera alors 
l'intégrale sur $a_{1} \in A_{1}^{\infty}$ 
du supremum sur $a_{0}^{1} \in A_{0,1}^{\infty}(T)$ 
et $y \in \Gamma$ de
\begin{equation}\label{eq:mainThmConv5}
e^{2(\rho_{\tlS}-\rho_{1})(H_{0}(a_{1}a_{0}^{1}))}
 \sigma_{1}^{2}(H_{1}(a_{1})-T) \ 
\sum_{\mathclap{
\begin{subarray}{c}
\eta \in \calR_{1},  \\
\zeta \in \calR_{2}, \, \xi \in \calR_{3}
\end{subarray}
}} \ 
|\Phi_{S}(y,\Ad((a_{0}^{1}a_{1})^{-1})(\zeta + \xi),\Ad((a_{0}^{1}a_{1})^{-1})\eta)|.
\end{equation}

\newcommand{\lenf}{\le_{n,f}}

Soit $n_{0} \in \N$ tel que l'on a
\[
\sum_{\xi \in \calR \smin \{0\}}\|\xi\|^{-n} < \infty, \quad \forall \, n \ge n_{0}.
\]
Soient $n, n_{1} \in \N$ tels que $n > n_{1} \ge n_{0}$. On suppose $n$ pair et on le laissera varier. L'entier $n_{1}$ va être fixé plus tard.
On utilisera la notation suivante pour $A,B \in \R$
\[
A \lenf B
\]
pour signifier qu'il existe une constante positive $c$ qui dépend, éventuellement, seulement de $n$ et de $f$
telle que $A \le c B$.

En faisant des intégrations par partie on a pour tous $X \in \ml_{\tlS}(\A)$ et
$Y \in \nl_{\tlS}^{\tltwo}(\A) \smin \{0\}$ que $\Phi_S(y,X, Y) $ s'exprime comme une somme finie 
d'expressions de type
\[
c_{2}(y)
\| Y \|^{-n}
\Phi_S^{(n)}(y,X,Y)
\]
où 
$c_{2}(y)$ dépend continument de $y$ 
et
\begin{displaymath}
\Phi_{S}^{(n)}(y,X,Y) = \int_{\nl_{\tlS}(\A)}
\partial_{n}f(y^{-1}(X + U)y)\psi 
(\langle Y,U\rangle)dU
\end{displaymath}
pour un opérateur différentiel $\partial_{n}$ convenable sur les 
fonctions de classe $C^{\infty}$ sur $\tlul(\A_{\infty})$. On a donc
\begin{displaymath}
|c_{2}(y)\Phi_{S}^{(n)}(y,X,Y) | \le 
\int\limits_{\nl_{\tlS}(\A)}
\sup_{y \in \Gamma}
| c_{2}(y)
\partial_{n}f(y^{-1}(X + U_{S})y)|dU_{S} =: 
\Psi_{S,n}(X).
\end{displaymath}
La fonction $\Psi_{S,n}$ est à décroissance rapide sur $\ml_{\tlS}(\A)$, donc
pour tous $\zeta \in \calR_{2}$ 
et $\xi \in \calR_{3} \smin \{0\}$ on a
\begin{gather*}
\Psi_{S,n}(\Ad((a_{0}^{1}a_{1})^{-1})(\zeta + \xi)) \lenf \|\Ad((a_{0}^{1}a_{1})^{-1})\zeta\|^{-n}
\|\Ad((a_{0}^{1}a_{1})^{-1})\xi\|^{-n_{1}}, \\
\Psi_{S,n}(\Ad((a_{0}^{1}a_{1})^{-1})\zeta) \lenf \|\Ad((a_{0}^{1}a_{1})^{-1})\zeta\|^{-n}.
\end{gather*}
On obtient alors
\begin{equation}\label{eq:mainThmCvg6}
\begin{split}
\sum_{\mathclap{
\begin{subarray}{c}
\eta \in \calR_{1},  \\
\zeta \in \calR_{2}, \, \xi \in \calR_{3}
\end{subarray}
}} \ 
|\Phi_{S}(y,\Ad((a_{0}^{1}a_{1})^{-1})(\zeta + \xi),\Ad((a_{0}^{1}a_{1})^{-1})\eta)| & \lenf \\
\sum_{\mathclap{
\begin{subarray}{c}
\eta \in \calR_{1},  \\
\zeta \in \calR_{2}, \, \xi \in \calR_{3}
\end{subarray}
}} \
\|\Ad((a_{0}^{1}a_{1})^{-1})\eta\|^{-n}
\Psi_{S,n}(\Ad((a_{0}^{1}a_{1})^{-1})(\zeta + \xi))  & \lenf \\
\dsl \
\sum_{\mathclap{\eta \in \calR_{1}}}\|\Ad((a_{0}^{1}a_{1})^{-1})\eta\|^{-n}
\rb \! \!
\dsl \
\sum_{\mathclap{\zeta \in \calR_{2}}}\|\Ad((a_{0}^{1}a_{1})^{-1})\zeta\|^{-n}
\rb \! \! &
\dsl 
1 + \sum_{\mathclap{\xi \in \calR_{3} \smin \{0\}}}\|\Ad((a_{0}^{1}a_{1})^{-1})\xi\|^{-n_{1}}
\rb.
\end{split}
\end{equation}

On introduit un peu de notations.
Soit $\Sigma_{\tlU}$ l'ensemble de racines 
pour l'action du tore 
$A_{0}$
sur l'algèbre de Lie $\tlul$. Pour $\mu \in \Sigma_{\tlU}$ 
et $X \in \tlul(\A)$
soit $X_{\mu}$ la projection de $X$ sur 
le sous-espace de $\tlul$ de valeur propre $\mu$. 
Notons dans ce cas $\Sigma(X) = \{\mu \in \Sigma_{\tlU} | X_{\mu} \neq 0\}$. 

Fixons $\Sigma \sbs \Sigma_{\tlU}$.
En utilisant l'équivalence de normes sur $\tlul(\A_{\infty})$, on peut 
fixer une constante $c_{\Sigma} >0 $ telle que pour 
tout $X \in \tlul(\A)$ tel que $\Sigma(X) = \Sigma$ on a
\begin{equation}\label{eq:AGineq2}
\|X\| \ge c_{\Sigma} \prod_{\mu \in \Sigma = \Sigma(X)}
\|X_{\mu}\|^{\frac{1}{|\Sigma|}}.
\end{equation}
On fixe maintenant $n_{1} := n_{0}|\Sigma_{\tlU}|$. 

En utilisant le lemme \ref{lem:DeltaSontLesMemes}, on fixe des constantes
$k_{\Sigma, \alpha} \in \R$ pour tout $\al \in \Delta_{0}$,
telles que 
  $\sum_{\mu \in \Sigma} \frac{1}{|\Sigma|} \mu = \sum_{\al \in \Delta_{0}} k_{\Sigma, \al} \al$. 
Pour tout $\calS \sbs \calR$ soit $\Sigma(\calS) = \{\Sigma(\xi)| 
\xi \in \calS \}$. On a alors

\blem\label{lem:SommeUtile} Pour tous $\calS \sbs \calR \smin \{0\}$, $a \in A_{0}^{\infty}$ et $n \ge n_{1}$ on a
\[
\sum_{\xi \in \calS}\|\Ad(a^{-1})\xi\|^{-n} \lenf
\sum_{\Sigma \in \Sigma(\calS)}e^{n\sum_{\al \in \Delta_{0}}k_{\Sigma, \al}\al(H_{0}(a))}.
\]
\bdem
En utilisant l'inégalité (\ref{eq:AGineq2}), on a 
\begin{multline*}
\sum_{\xi \in \calS} \|\Ad(a^{-1})\xi\|^{-n} \! \! = 
\sum_{\mathclap{\Sigma \in \Sigma(\calS)}} \quad
\sum_{\mathclap{\begin{subarray}{c}
\xi \in \calS \\
\Sigma(\xi) = \Sigma
\end{subarray}}}
\|\Ad(a^{-1})\xi\|^{-n} \le  
\sum_{\mathclap{\Sigma\in \Sigma(\calS)}}
c_{\Sigma}
\displaystyle \left(
\! \! \prod_{\mu \in\Sigma} \! \quad
(
\sum_{\mathclap{\begin{subarray}{c}
\xi \in \calR \smin \{0\} \\
\Sigma(\xi) = \{\mu\}
\end{subarray}}}
 \|\Ad (a^{-1})\xi\|^{-nn_{\Sigma,\mu}}) \! \right)   \\
=  \sum_{\Sigma \in \Sigma(\calS)}
 c_{\Sigma}
 \displaystyle \left(
 \prod_{\mu \in \Sigma}e^{nn_{\Sigma,\mu}\mu(H_{0}(a_{0}))}
(\sum_{\begin{subarray}{c}
\xi \in \calR \smin \{0\} \\
\Sigma(\xi) = \{\mu\}
\end{subarray}}
 \|\xi\|^{-nn_{\Sigma,\mu}})\right)  
 \lenf
\sum_{\mathclap{\Sigma \in \Sigma(\calS)}}e^{n\sum_{\al \in \Delta_{0}}k_{\Sigma, \al}\al(H_{0}(a))}.
\end{multline*}
\edem
\elem

  On pose $k_{1} = n_{1} \max \{ k_{\Sigma, \al} | \Sigma \sbs \Sigma_{\tlU}, \al \in \Delta_{0}\}$ 
  et $k_{2} = \min \{ k_{\Sigma, \al} | \Sigma \sbs \Sigma_{\tlU}, \al \in \Delta_{0}, k_{\Sigma, \al} > 0 \}$. 
Puisque on suppose $P_{1} \sbn P_{2}$, on a en particulier $P_{0} \neq U$ et donc $k_{1}, k_{2} >0$.

\bcor\label{cor:somExpMaj1}
 Pour tout $a \in A_{0}^{\infty}$ tel que $\tau_{0}^{2}(H_{0}(a) - T_{1}) = 1$ on a
\[
\dsl 
\sum_{\eta \in \calR_{1}}\|\Ad(a^{-1})\eta\|^{-n}
\rb
\dsl 
\sum_{\zeta \in \calR_{2}}\|\Ad(a^{-1})\zeta\|^{-n} \rb \lenf 
e^{-nk_{2}\sum_{\al \in \Delta_{0}^{2} \smin \Delta_{0}^{1}} \al(H_{0}(a))}.
\]
\bdem
Soit $a \in A_{0}^{\infty}$ comme dans l'énoncé. 
En utilisant le lemme \ref{lem:SommeUtile} pour $\calS  = \calR_{1}$ 
et pour $\calS  = \calR_{2}$ 
on voit qu'il suffit de vérifier pour tout $\Sigma_{1} \in \Sigma(\calR_{1})$ 
et tout $\Sigma_{2} \in \Sigma(\calR_{2})$ qu'on a
\begin{gather*}
e^{n\sum_{\al \in \Delta_{0}}k_{\Sigma_{1}, \al}\al(H_{0}(a))} \lenf
e^{-nk_{2}\sum_{\al \in \Delta_{0}^{2} \smin \Delta_{0}^{S}} \al(H_{0}(a))}, \\
e^{n\sum_{\al \in \Delta_{0}}k_{\Sigma_{2}, \al}\al(H_{0}(a))} \lenf
e^{-nk_{2}\sum_{\al \in \Delta_{0}^{S} \smin \Delta_{0}^{1}} \al(H_{0}(a))}.
\end{gather*}
En effet,  posons $\Theta_{1} = \Delta_{0}^{2} \smin \Delta_{0}^{S}$ 
et $\Theta_{2} = \Delta_{0}^{S} \smin \Delta_{0}^{1}$. 
Alors, 
pour $i=1,2$,
on a $k_{\Sigma_{i}, \al} = 0$ pour tout $\al \in \Delta_{0} \smin \Delta_{0}^{2}$ et 
$k_{\Sigma_{i}, \al} \le 0$
pour tout $\al \in \Delta_{0}^{2}$. De plus, il résulte de la définition 
de $(\bar \nl_{\tlS}^{\tltwo})'$ (resp. $\ml_{\tlS, \tlone}'$)
que $k_{\Sigma_{1},\al} \le -k_{2}$  (resp. $k_{\Sigma_{2},\al} \le -k_{2}$)
pour tout $\al \in \Theta_{1}$ (resp. $\al \in \Theta_{2}$). 
 Il existe une constante $C_{T_{1}}$ telle que $\al(H_{0}(a)) > C_{T_{1}}$ 
 pour tout $\al \in \Delta_{0}^{2}$. 
 Donc, pour tout $\al_{i} \in \Delta_{0} \smin \Theta_{i}$ 
 et tout $\al_{i}' \in \Theta_{i}$ 
 on a
 \[
 e^{nk_{\Sigma_{i}, \al_{i}}\al_{i}(H_{0}(a))} \lenf 1, 
\quad  e^{n k_{\Sigma_{i}, \al_{i}'}\al_{i}'(H_{0}(a))} \lenf
e^{-nk_{2} \al_{i}'(H_{0}(a))},
 \quad i = 1,2.
 \]
 \edem
\ecor

\bcor\label{cor:somExpMaj2}
 Pour tout $a \in A_{0}^{\infty}$ tel que $\tau_{0}^{2}(H_{0}(a) - T_{1}) = 1$ on a
\[
1 + \sum_{\eta \in \calR_{3} \smin \{0\}}\|\Ad(a)^{-1}\xi\|^{-n_{1}}  \lenf 
e^{k_{1}\sum_{\al \in \Delta_{0}^{S}} \al(H_{0}(a))}.
\]
\bdem
Soit $a \in A_{0}^{\infty}$ comme dans l'énoncé. 
Il est clair qu'on a $1 \lenf e^{k_{1}\sum_{\al \in \Delta_{0}^{S}} \al(H_{0}(a))}$.
En utilisant le lemme \ref{lem:SommeUtile} pour $\calS = \calR_{3} \smin \{0\}$ et $n = n_{1}$ 
on voit qu'il suffit de vérifier pour tout $\Sigma \in \Sigma(\calR_{3} \smin \{0\})$ 
qu'on a
\begin{equation}\label{eq:ekuneku}
e^{n_{1}\sum_{\al \in \Delta_{0}}k_{\Sigma, \al}\al(H_{0}(a))} \lenf
e^{k_{1}\sum_{\al \in \Delta_{0}^{S}} \al(H_{0}(a))}.
\end{equation}
Or, on a $k_{\Sigma, \al} = 0$ pour tout $\al \in \Delta_{0} \smin \Delta_{0}^{S}$ et
$n_{1}k_{\Sigma, \al} \le k_{1}$ pour tout $\al \in \Delta_{0}^{S}$ par définition de $k_{1}$. D'autre part, 
il existe une constante $C_{T_{1}}$ telle que $\al(H_{0}(a)) > C_{T_{1}}$ 
 pour tout $\al \in \Delta_{0}^{2}$, 
ce qui démontre (\ref{eq:ekuneku}) et le corollaire. 
\edem
\ecor

En regardant (\ref{eq:mainThmConv5}) et (\ref{eq:mainThmCvg6}), 
et en utilisant les corollaires \ref{cor:somExpMaj1} et \ref{cor:somExpMaj2} on se ramène 
a majorer le supremum sur $a_{0}^{1} \in A_{0,1}^{\infty}(T)$ de
\begin{multline*}
\int\limits_{A_{1}^{\infty}} 
 \sigma_{1}^{2}(H_{1}(a_{1})-T) 
 e^{(\la -nk_{2}\sum_{\al \in \Delta_{0}^{2} \smin \Delta_{0}^{1}} )\al(H_{0}(a_{0}^{1}a_{1}))}
 da_{1} = 
 e^{(\la -nk_{2}\sum_{\al \in \Delta_{0}^{2} \smin \Delta_{0}^{1}} \al)(H_{0}(a_{0}^{1}))} \\
 \cdot 
 \int\limits_{\all_{1}^{2}} 
  e^{-nk_{2}\sum_{\al \in \Delta_{1}^{2}} \al(H_{1}^{2})}
   \int\limits_{\all_{2}} 
 \sigma_{1}^{2}(H_{1}^{2} + H_{2}-T) 
 e^{\la(H_{1}^{2} + H_{2})}
 dH_{2}dH_{1}^{2} 
\end{multline*}
où $\la = 2(\rho_{\tlS}-\rho_{1}) + k_{1}\sum_{\al \in \Delta_{0}^{S}} \al \in \all_{S}^{*}$.

Il est clair qu'il existe une constante $c_{1} >0$ qui ne dépend pas de $n$ telle que:
\[
\sup_{\mathclap{
a_{0}^{1} \in A_{0,1}^{\infty}(T)
}} \ 
 e^{(\la -nk_{2}\sum_{\al \in \Delta_{0}^{2} \smin \Delta_{0}^{1}} \al)(H_{0}(a_{0}^{1}))} \lenf 
 e^{c_{1} \|T\|}.
\]
 
 On a aussi
\blem 
Il existe une constante $c_{2} >0$ telle que pour tout $H_{1}^{2} \in \all_{1}^{2}$ on a
\[
\int_{\all_{2}}\sigma_{1}^{2}(H_{2} + H_{1}^{2} - T)e^{\la(H_{2})}dH_{2} 
\le c_{2}\tau_{1}^{2}(H_{1}^{2} - T)e^{c_{2}\|T\|}e^{c_{2}\sum_{\al \in \Delta_{1}^{2}}\al(H_{1}^{2})}.
\]
\bdem
En vertu du lemme 6.1 de \cite{arthur3} ainsi que du corollaire 6.2 de loc. cit. 
il existe une constante $c' >0$ telle que $\sigma_{1}^{2}(H_{2} + H_{1}^{2}-T)=1$ 
implique
\[
\tau_{1}^{2}(H_{1}^{2} - T) = 1 \quad \text{et} \quad \|H_{2}\| \le c'(1 + \|H_{1}^{2}\| +\|T\|).
\]
En utilisant cela, il existe une constante $c''>0$ telle que $\la(H_{2}) \le c''(\|H_{1}^{2}\| + \|T\|)$. 
On a donc que l'intégrale dans l'énoncé est bornée par
\[
\tau_{1}^{2}(H_{1}^{2} - T)e^{c''(\|H_{1}^{2}\|  +\|T\|)}
\int_{
H_{2} \in \all_{2} | \,
\|H_{2}\| \le c'(1 + \|H_{1}^{2}\| +\|T\|)
}dH_{2} \le c''' \tau_{1}^{2}(H_{1}^{2} - T)e^{c'''(\|H_{1}^{2}\|  +\|T\|)}
\]
pour une constante $c''' > 0$.
On obtient le résultat en utilisant l'équivalence de normes sur $\all_{1}^{2}$, 
en remarquant que $\Delta_{1}^{2}$ est une base de 
$(\all_{1}^{2})^{*}$ et que si $\tau_{1}^{2}(H_{1}^{2} - T) = 1$ alors 
$|\al(H_{1}^{2})| = \al(H_{1}^{2})$ pour tout $\al \in \Delta_{1}^{2}$.
\edem
\elem
 
 Grâce au lemme ci-dessus il nous reste a majorer
 \[
 \int\limits_{\all_{1}^{2}} 
 \tau_{1}^{2}(H_{1}^{2} - T)
  e^{-(nk_{2}- c_{2})\sum_{\al \in \Delta_{1}^{2}} \al(H_{1}^{2})}
 dH_{1}^{2}
 \]
 où $c_{2}$ ne dépend pas de $n$. En prenant $n$ suffisamment grand ceci 
est majoré par $ e^{-(nk_{2}- c_{2})\sum_{\al \in \Delta_{0}^{2} \smin \Delta_{0}^{1}} \al(T)}$.
On rappelle qu'on a un $\varepsilon_{0} >0$ tel que 
$\al(T) > \varepsilon_{0} \|T\|$ pour tout $\al \in \Delta_{0}$. 
 Il est clair donc que pour tout $N >0$ il existe un $n \in \N$ qui dépend de $\varepsilon_{0}$
 tel que
 \[
 e^{(c_{1}+c_{2})(\|T\|)}e^{-(nk_{2}- c_{2})\sum_{\al \in \Delta_{0}^{2} \smin \Delta_{0}^{1}} \al(T)}
 \le e^{\|T\| (c_{1} +c_{2} -(nk_{2}- c_{2})d_{1}^{2}\varepsilon_{0})} \le 
 e^{-N \|T\|}.
 \]
 Ce qui conclut la preuve des théorèmes \ref{thm:MainConv} et \ref{thm:MainConv2}. 
 \edem
 \etheo

\newcommand{\spsl}{\psl{\upla} \sps}
\newcommand{\sbsl}{\sbs_{\upla}}

\section{Les propriétés qualitatives}\label{qualitatives}

Pour une fonction $f \in \mathcal{S}(\tlul(\A))$ 
et $T \in T_{+} + \mathfrak{a}_{0}^{+}$
on note
\begin{displaymath}
k^{T}(x) = 
k_{f}^{T}(x)  = \sum_{\mathfrak{o} \in \mathcal{O}}
k_{f,\mathfrak{o}}^{T}(x), \quad 
x \in U(\rmF)\backslash U(\A).
\end{displaymath}
Grâce au théorème \ref{thm:MainConv} les 
distributions suivantes:
\begin{align*}
J_{\mathfrak{o}}^{T}(f) &= 
\int_{U(\rmF)\backslash U(\A)}k_{f,\mathfrak{o}}^{T}(x)dx, &
\ \mathfrak{o} \in \mathcal{O}, \ T \in T_{+} + \mathfrak{a}_{0}^{+},
\\ 
J^{T}(f) &= \int_{U(\rmF)\backslash U(\A)}k_{f}^{T}(x)dx&
\end{align*}
sont bien définies.

Dans le paragraphe \ref{par:asymptChapt} 
on démontrera que la fonction
$T \mapsto J_{\mathfrak{o}}^{T}(f)$ 
est un polynôme-exponentielle dont le terme purement polynomial, 
noté $J_{\ol}(f)$, ne dépend pas de $T$. 
Pour bien énoncer ce résultat on introduit d'abord dans le paragraphe 
\ref{par:JTforLevis} les distributions $J_{\ol}^{M_{Q},T}$ pour tout sous-groupe 
parabolique standard $Q$ et 
dans le paragraphe \ref{par:fonsPolExp} on étudie les fonctions de type 
polynôme-exponentielle. 
La suite de ce chapitre est consacrée 
aux propriétés des distributions $J_{\ol}$.

\subsection{Une généralisation du théorème \ref{thm:MainConv}}\label{par:JTforLevis}

Soit $G = \prod_{i=1}^{k}\Res_{\rmE/\rmF}\Gl_{n_{i}}$ où $k \in \N$ et 
$n_{i} \in \N^{*}$ pour $i = 1, \ldots ,k$. Soient aussi $U' = U(V', \Phi')$ et 
$\tlU' = U(W', \tlPhi')$ un couple de groupes unitaires munis de l'inclusion 
$U' \hrar \tlU'$ comme dans le paragraphe 
\ref{generalities}. 
On va généraliser le théorème \ref{thm:MainConv} au cas de l'inclusion 
$G \times U' \hrar G \times \tlU'$. 

Notons $\gl = \Lie(G)$, $\ul' = \Lie(U')$ 
et $\tlul' = \Lie(\tlU')$. 
Pour $X \in (\gl\times \tlul')(\rmF)$ soient $X_{1} \in \gl(\rmF)$ 
et $X_{2} \in \tlul'(\rmF)$ tels que $X = X_{1}  + X_{2}$. 
Soit $\calO^{G \times U'}$ la relation d'équivalence sur $(\gl\times \tlul')(\rmF)$ 
définie de la façon suivante. On a $X = X_{1} + X_{2} \sim 
Y = Y_{1} + Y_{2}$ si et seulement si 
les polynômes caractéristiques de $X_{1}$ et $Y_1$ 
co\"incident et si $X_{2}$ et $Y_{2}$ sont dans la même classe 
pour 
la relation d'équivalence dans $\tlul'(\rmF)$ 
décrite dans le paragraphe 
\ref{lesInvs} par rapport à l'inclusion $U' \hrar \tlU'$.

Fixons $P_{0}$ un sous-$\rmF$-groupe parabolique minimal 
de $G \times U'$ et fixons aussi $M_{0}$ une partie de Levi de $P_{0}$. 
Soit $P$ un sous-groupe parabolique standard de $G \times U'$. 
On a $P = P_{G} \times P_{U'}$ où $P_{G}$ (resp. $P_{U'}$) est un sous-groupe parabolique 
standard de $G$ (resp. $U'$) par rapport au sous-groupe parabolique minimal $P_{0} \cap G$ 
(resp. $P_{0} \cap U'$). On pose alors $\tlP = P_{G} \times \tlP_{U'}$ 
où $\tlP_{U'}$ c'est le sous-groupe parabolique de $\tlU'$
associé à $P_{U'}$ par la procédure décrite dans 
 le paragraphe \ref{generalities}. On note $\ml_{\tlP}$ l'algèbre de Lie 
 de la partie de Levi de $\tlP$ contenant $M_{0}$.

Pour une fonction 
$f \in \mathcal{S}((\gl\times \tlul')(\A))$,  
un sous-groupe parabolique standard $P$ de $G \times U'$ 
et une classe $\ol \in \calO^{G \times U'}$ on pose
\[
k_{f,P, \ol}(x) = 
\sum_{\xi \in \ml_{\tlP}(\rmF) \cap \ol}
\int_{\Lie (N_{\tlP}) (\A)}
f(x^{-1}(\xi + U_{P})x)dU_{P}, \quad 
x \in  (G \times U')(\A).
\]
Pour un $T \in (\mathfrak{a}_{P_{0}}^{G \times U'})^{+}$ on pose donc
\[
k^{T}_{f,\ol}(x) = 
\sum_{P \sps P_{0}} (-1)^{d_{P}^{G \times U'}}
\sum_{\delta \in P(\rmF) \bsl (G \times U')(\rmF)}
\htau_{P}^{G \times U'}(H_{P}(\delta x)-T)
k_{f,P,\ol}(\delta x).
\]

\btheo Soit $f \in \mathcal{S}((\gl\times \tlul')(\A))$, 
alors pour tout $T \in T_{+} + (\mathfrak{a}_{P_{0}}^{G \times U'})^{+}$ on a:
\[
\sum_{\ol \in \calO^{G \times U'}}
\int_{(G \times U')(\rmF) \backslash (G \times U')(\A)^{1}}
|k_{f,\ol}^{T}(x)|dx < \infty.
\]
\bdem 
Le résultat découle des théorèmes 3.1 ci-dessus et 3.1 de \cite{chaud}.
\edem 
\etheo

Notons alors:
\begin{displaymath}
J_{\ol}^{G \times U',T}(f) =\int_{(G \times U')(\rmF) \backslash (G \times U')(\A)^{1}}
k_{f,\ol}^{T}(x)dx. 
\end{displaymath}

Revenons au groupe $U$.
Soit $Q$ le sous-groupe parabolique standard de $U$
défini par le $l+1$-uplet 
$(i_{0},i_{1},\ldots, i_{l})$, de façon que $Q$ soit 
le stabilisateur du drapeaux isotrope (\ref{drapParabStand}). 
Les choix qu'on a faits dans le paragraphe \ref{generalities} 
nous permettent d'écrire $V_{i_{l}} = V_{Q}$ et
$V = V_{Q} \oplus Z_{Q} \oplus V_{Q}^{\sharp}$. 
La restriction de 
$\Phi$ \`a $Z_{Q}$ est une forme hermitienne 
non-dégénérée. On a donc
\begin{equation*}
M_{Q} \cong \prod_{j = 0}^{l-1} \Res_{\rmE/\rmF}\Gl_{i_{j+1}-i_{j}} 
\times U(Z_{Q},\Phi_{| Z_{Q}}). 
\end{equation*} 
En particulier $M_{Q}$ est de type considéré au début de ce paragraphe. 

Soit $\ol \in \mathcal{O}$, il existe 
$\ol_{Q,1}, \ldots, \ol_{Q,m} \in \calO^{M_{Q}}$,
où $0 \le m < \infty$,
tels que
\begin{equation}\label{eq:mltlQcapol}
\ml_{\tlQ}(\rmF) \cap \mathfrak{o} = 
\coprod_{i=1}^{m} \ol_{Q,i}.
\end{equation}
On définit alors les distributions $J_{\ol}^{M_{Q},T}$ 
et $J^{M_{Q},T}$ sur 
$\calS(\ml_{\tlQ}(\A))$ par:
\begin{equation}\label{eq:indDistDef}
J^{M_{Q},T}_{\mathfrak{o}}(f) = 
\sum_{i=1}^{m} 
J^{M_{Q},T}_{\ol_{Q,i}}(f), \quad
J^{M_{Q},T}(f) = \sum_{\ol \in \calO}J^{M_{Q},T}_{\mathfrak{o}}(f)
\end{equation}
où pour $\ol_{Q} \in \calO^{M_{Q}}$, $J^{M_{Q},T}_{\ol_{Q}}$ 
c'est la distribution décrite ci-dessus par rapport au sous-groupe parabolique 
minimal $P_{0} \cap M_{Q}$, le sous-groupe de Levi $M_{0}$ 
et le sous-groupe compact maximal $K \cap M_{Q}(\A)$ de $M_{Q}(\A)$ admissible par rapport à $M_{0}$.

Pour $f \in \calS(\tlul(\A))$ on 
pose
\begin{equation}\label{eq:fQdef}
f_{Q}(X) = \int_{K}\int_{\nl_{\tlQ}(\A)}f(k\inv(X+U_{Q})k)dU_{Q}dk, \quad 
X \in \ml_{\tlQ}(\A);
\end{equation}
alors $f_{Q} \in \calS(\ml_{\tlQ}(\A))$. 
Notons que l'application $Q \sps P \mapsto M_{Q} \cap P$
définit une bijection entre les sous-groupes paraboliques de 
$U$ contenus dans $Q$ et les sous-groupes paraboliques 
de $M_{Q}$ contenant $P_{0} \cap M_{Q}$. 
On s'aperçoit alors qu'on a
\begin{equation}\label{eq:JMQisThis}
\begin{split}
J^{M_{Q},T}_{\mathfrak{o}}(f_{Q}) & = 
\int_{M_{Q}(\rmF) \backslash M_{Q}(\A)^{1}}
\sum_{i=1}^{m}k_{f_{Q},\ol_{Q,i}}^{T}(m)dm  \\ &= 
\int_{M_{Q}(\rmF) \backslash M_{Q}(\A)^{1}}
\sum_{P \sbs Q}(-1)^{d_{P}^{Q}} 
\sum_{\eta \in (P \cap M_{Q})(\rmF) \bsl M_{Q}(\rmF)}
\htau_{P}^{Q}(H_{P}(\eta m)-T) \\
& 
\dsl
\sum_{\xi \in \ml_{\tlP}(\rmF) \cap \ol}
\int_{\nl_{\tlP}^{\tlQ}(\A)}
f_{Q}(\Ad((\eta m)^{-1})(\xi + U_{P}^{Q}))dU_{P}^{Q}
\rb dm.
\end{split}
\end{equation}

\subsection{Polynômes-exponentielles}\label{par:fonsPolExp}

Soit $\calV$ un $\R$-espace vectoriel de dimension finie. 
Par un polynôme-exponentielle sur $\calV$ on entend une fonction 
sur $\calV$ de la forme
\[
f(v) = \sum_{\la \in \calV^{*}}e^{\la(v)}P_{\la}(v), \quad v \in \calV
\]
où
$P_{\la}$ est un polynôme sur $\calV$ 
à coefficients complexes, nul 
pour presque tout $\la \in \calV^{*}$. 
On appelle $\la \in \calV^{*}$ tels que $P_{\la} \neq 0$ les exposants 
de $f$ et le polynôme correspondant à $\la = 0$ le terme 
purement polynomial de $f$.
On a alors le résultat d'unicité suivant: 
si $f$ est comme ci-dessus 
et $g = \sum_{\la \in \calV^{*}}e^{\la}Q_{\la}$ est un 
polynôme-exponentielle sur $\calV$ tel que 
$g(v) = f(v)$ pour tout $v \in \calV$ alors 
pour tout $\la \in \calV^{*}$ on a $P_{\la} = Q_{\la}$.

On rappelle qu'à
la fin du paragraphe \ref{par:haarMesSubs} 
on a défini un $\upla_{Q} \in \all_{Q}^{*}$ pour 
tout sous-groupe parabolique standard $Q$.
On aura besoin du lemme suivant.
\blem\label{lem:positiveLinear}
 Soit $Q$ un sous-groupe parabolique standard 
 de $U$. Alors, pour tout  $\varpi^{\vee} \in  \hDelta_{Q}^{\vee}$ 
on a $\upla_{Q}(\varpi^{\vee}) > 0$.

\bdem
Soit $\calS \sbs \all_{Q}^{*}$ l'ensemble des poids 
pour l'action du tore $A_{Q}$ 
sur $V_{Q}$. Alors $\calS$ est une base de $\all_{Q}^{*}$ 
et on note $\calS^{\vee} \sbs \all_{Q}$ sa base duale. 
Comme on a déjà remarqué à la fin du paragraphe 
\ref{par:haarMesSubs}, on a
$\upla_{Q} = \sum_{e^{*} \in \calS}n_{e^{*}}e^{*}$ 
avec $n_{e^{*}} \in \N^{*}$.
D'autre part, tout $\varpi^{\vee} \in \hDelta_{Q}^{\vee}$ 
est une combinaison linéaire à coefficients positifs des 
vecteurs dans $\calS^{\vee}$, d'où le résultat. 
\edem
\elem

Donc, pour tout  $Q$ et tout $R \sps Q$ différent de $U$ 
on a que $\upla_{R}$, vu comme un élément de 
$\all_{Q}^{*}$ grâce à la décomposition (\ref{eq:decomp}), 
est non-nul.

Soient $Q \sbs R$ deux sous-groupes paraboliques standards. 
Notons $v_{Q}^{R}$ 
 (resp. $v_{R}'$)
 le volume dans $\all_{Q}^{R}$ (resp. $\all_{R}$) du 
 parallélotope engendré par $(\hDelta_{Q}^{R})^{\vee}$
 (resp. $\Delta_{R}^{\vee}$).
 Suivant le paragraphe 2 de  \cite{arthur2}, posons
\begin{equation}\label{eq:thetaHatDef}
\hat \theta_{Q}^{R}(\mu) = 
(v_{Q}^{R})^{-1}
\prod_{\varpi \in (\hDelta_{Q}^{R})^{\vee}} \mu(\varpi^{\vee}), \quad
\theta_{R}(\mu) = 
(v_{R}')^{-1}
\prod_{\al \in \Delta_{R}} \mu(\al^{\vee}), \quad \mu \in \all_{0}^{*} \otimes_{\R} \C.
\end{equation}
Quand $R = U$ on écrira $\hat \theta_{Q} = \hat \theta_{Q}^{U}$.

Toujours supposant $R \sps Q$, pour 
$X \in \all_{Q}$ on note $X_{R}$ sa projection à $\all_{R}$ 
selon la décomposition (\ref{eq:decomp}). 
Suivant loc. cit., posons
\begin{equation}\label{eq:GammaQDef}
\Gamma_{Q}'(H,X) = 
\sum_{R \sps Q}(-1)^{d_{Q}^{R}}
\htau_{R}(H_{R}-X_{R})\tau_{Q}^{R}(H), \quad H,X \in \all_{Q}.
\end{equation}

On démontrera le lemme suivant:
\blem\label{lem:pQExplicit}
 Soit $Q$ un sous-groupe parabolique standard de $U$.
Alors, pour tout $R \sps Q$ 
il existe un polynôme $p_{Q,R}$ de degré au plus $d_{Q}$ sur $\all_{R}$ 
tel que la fonction
\[
p_{Q}(X) := \int_{\all_{Q}}e^{\upla_{Q}(H)}\Gamma_{Q}'(H,X)dH, \quad 
X \in \all_{Q}
\]
égale
\[
\sum_{R \sps Q}e^{\upla_{R}(X_{R})}
p_{Q,R}(X_{R})
\]
où $p_{Q,U}(X_{U}) = (-1)^{d_{Q}}\hat \theta_{Q}(\upla_{Q})^{-1}$.
En particulier, la fonction $p_{Q}$ est un polynôme-exponentielle
dont le terme purement polynomial est constant 
et égale $(-1)^{d_{Q}}\hat \theta_{Q}(\upla_{Q})^{-1}$.

\brems On ne prétend pas que les polynômes $p_{Q,R}$ 
sont uniquement déterminés pour tout $R \sps Q$. En effet, il 
arrive que $\upla_{R} = \upla_{R'}$ pour $R \neq R'$. Cependant, $U$ 
est le seul sous-groupe parabolique $R \sps Q$
tel que $\upla_{R} = 0$ d'où l'unicité du
terme $p_{Q,U}$.
\erems

\bdem
Étudions l'intégrale
\begin{equation}\label{entireint}
\int_{\mathfrak{a}_{Q}}
e^{\mu(H)}
\Gamma'_{Q}(H,X)dH, \quad \mu \in \all_{Q}^{*} \otimes_{\R} \C.
\end{equation}
 Il résulte du lemme 2.1 de \cite{arthur2} que 
pour un $X$ fixé, la fonction $H \mapsto \Gamma_{Q}'(H,X)$ 
est à support compact dans $\all_{Q}$. 
L'intégrale ci-dessus est donc bien définie et aussi la fonction 
$p_{Q}$ est bien définie sur $\all_{Q}$.

 D'après le lemme 
2.2  dans \cite{arthur2} on obtient que 
l'intégrale (\ref{entireint}) est une fonction enti\`ere 
de la variable $\mu$ et sa valeur est donnée  
par
\begin{equation}\label{eq:intcalculated}
\sum_{R \supseteq Q}(-1)^{d_{Q}^{R}}
e^{\mu(X_{R})}
\hat \theta_{Q}^{R}(\mu)\inv\theta_{R}(\mu)\inv
\end{equation}
pour tout $\mu\ \in \mathfrak{a}_{Q}^{*}\otimes_{\R}\C$ 
pour lesquelles cela a un sens.

On n'a pas le droit d'utiliser cette formule pour $\upla_{Q}$ 
car certains $\hat \theta_{Q}^{R}(\upla_{Q})$ et 
$\theta_{R}(\upla_{Q})$ vont s'annuler. 
On proc\`ede alors comme suit. 
On fixe un 
$\varepsilon \in \mathfrak{a}_{Q}^{*} \otimes_{\R}\C$ tel que $\hat \theta_{Q}^{R}(\varepsilon) \neq 0$, 
$\theta_{R}(\varepsilon) \neq 0$  pour 
tout $R \supseteq Q$. Alors il en est de m\^eme pour 
$\upla_{Q} + t\varepsilon$ o\`u $t \in \R \smallsetminus \{0\}$ 
suffisamment petit. 
Pour calculer $p_{Q}(X)$ on met $\mu = \upla_{Q} + t\varepsilon$
dans (\ref{eq:intcalculated}) et 
on considère la fonction
\begin{equation*}
\R \ni t \mapsto 
\sum_{R \supseteq Q}(-1)^{d_{Q}^{R}}
e^{(\upla_{Q} + t\varepsilon)(X_{R})}
\hat \theta_{Q}^{R}(\upla_{Q} + t\varepsilon)\inv\theta_{R}
(\upla_{Q} + t\varepsilon)\inv.
\end{equation*}
C'est une fonction analytique et sa valeur en $t = 0$ 
égale $p_{Q}(X)$. Pour calculer $p_{Q}(X)$ il suffit alors de
développer cette fonction en $t=0$ et de prendre le terme constant. 

Fixons $R \sps Q$ et regardons 
\[
\R \ni t \mapsto 
e^{t\varepsilon(X_{R})}
\hat \theta_{Q}^{R}(\upla_{Q} + t\varepsilon)\inv\theta_{R}
(\upla_{Q} + t\varepsilon)\inv.
\]
Notons $p_{Q,R, \varepsilon}(X_{R})$ son terme constant 
dans son développement 
en $t = 0$. Il est clair que $p_{Q,R,\varepsilon}$ 
est un polynôme en la variable 
$X_{R} \in \all_{R}$ de degré au plus $d_{Q}$.

On obtient donc
\begin{equation*}
p_{Q}(X) = 
\sum_{R \supseteq Q}(-1)^{d_{Q}^{R}}e^{\upla_{R}(X_{R})}
p_{Q,R,\varepsilon}(X_{R}),
\end{equation*}
où on utilise le fait que la restriction de $\upla_{Q}$ à $\all_{R}$ 
c'est $\upla_{R}$ pour $R \sps Q$.

Regardons le terme $p_{Q,U, \varepsilon}(X_{U}) $. 
On a 
\[
\lim_{t \rightarrow 0}
e^{t\varepsilon(X_{U})}
\hat \theta_{Q}^{U}(\upla_{Q} + t\varepsilon)\inv\theta_{U}
(\upla_{Q} + t\varepsilon)\inv = 
\hat \theta_{Q}(\upla_{Q})^{-1}
\]
 en vertu du 
lemme \ref{lem:positiveLinear}, d'où le résultat. 
\edem
\elem

\subsection{Le comportement asymptotique en $T$}\label{par:asymptChapt}

On démontre la proposition suivante.

\brop\label{prop:mainQualitProp}
Soient $f \in \calS(\tlul(\A))$, $T' \in T_{+} + \all_{0}^{+}$, 
$\ol \in \calO$ et $T \in T' + \all_{0}^{+}$. Alors
\[
J_{\ol}^{T}(f) = \sum_{Q \sps P_{0}}
p_{Q}(T_{Q} - T'_{Q})e^{\upla_{Q}(T_{Q}')}J_{\ol}^{M_{Q}, T'}(f_{Q}).
\]
où pour un sous-groupe parabolique $Q \sps P_{0}$, 
la fonction $p_{Q}$ est définie dans le lemme \ref{lem:pQExplicit}, 
$\upla_{Q} \in \all_{Q}^{*}$ est défini 
à la fin du paragraphe \ref{par:haarMesSubs}, 
la distribution $J_{\ol}^{M_{Q}, T'}$ 
est définie dans le paragraphe \ref{par:JTforLevis} et 
$f_{Q} \in \calS(\ml_{\tlQ}(\A))$ est définie par (\ref{eq:fQdef}) 
dans le même paragraphe.

\bdem 
Il est démontré dans 
le paragraphe 2
de \cite{arthur2}, que les fonctions $\Gamma_{Q}'$, 
définies par (\ref{eq:GammaQDef}), 
vérifient la relation suivante: 
pour tout sous-groupe parabolique $P$ de $U$, on a
\begin{equation}\label{GammaPrimeRecurrence}
\htau_{P}(H - X) = 
\sum_{Q \supseteq P}(-1)^{d_{Q}}
\htau_{P}^{Q}(H)
\Gamma_{Q}'(H,X), \quad H,X \in \all_{P}.
\end{equation}
Fixons un $T'\in T_{+}+ \mathfrak{a}_{0}^{+} $ 
et soit $T \in T'+ \mathfrak{a}_{0}^{+} $. En utilisant 
l'égalité ci-dessus avec $H = H_{P}(\delta x) -T'$ 
et $X = T-T'$ on a
\begin{multline}\label{eq:useGammaJol}
J_{\mathfrak{o}}^{T}(f) = 
\int_{U(\rmF)\backslash U(\A)}
\sum_{P \sps P_{0}}(-1)^{d_{P}}
\sum_{\delta \in P(\rmF)\backslash U(\rmF)}
\sum_{Q \supseteq P}(-1)^{d_{Q}}
\Psi_{P,Q,\mathfrak{o}}^{T,T'}(\delta x)dx = \\
\sum_{Q \sps P_{0}}\int_{Q(\rmF)\backslash U(\A)}
\sum_{P\subseteq Q}(-1)^{d_{P}^{Q}}
\sum_{\eta \in (P \cap M_{Q})(\rmF)\backslash M_{Q}(\rmF)}
\Psi_{P,Q,\mathfrak{o}}^{T,T'}(\eta x)dx
\end{multline}
o\`u 
\begin{displaymath}
\Psi_{P,Q,\mathfrak{o}}^{T,T'}(x) = k_{P,\mathfrak{o}}(x)\htau_{P}^{Q}(H_{P}(x)-T')
\Gamma'_{Q}(H_{Q}(x)-T',T-T').
\end{displaymath}
Posons
$x = namk$ o\`u $n \in N_{Q}(\rmF)\backslash N_{Q}(\A)$, 
$m \in M_{Q}(\rmF)\backslash M_{Q}(\A)^{1}$, 
$a \in A_{Q}^{\infty}$ et $k \in K$. Donc 
$dx = e^{-2\rho_{Q}(H_{Q}(a))}dndadmdk$. 
On peut donc changer l'intégrale sur 
$Q(\rmF)\backslash U(\A)$ en l'intégrale 
\begin{displaymath}
\int_{A_{Q}^{\infty}}
\int_{M_{Q}(\rmF)\backslash M_{Q}(\A)^{1}}
\int_{K}
\int_{N_{Q}(\rmF)\backslash N_{Q}(\A)}e^{-2\rho_{Q}(H_{Q}(a))}.
\end{displaymath}
On a
\begin{displaymath}
\int_{N_{Q}(\rmF)\backslash N_{Q}(\A)}
\Psi_{P,Q,\mathfrak{o}}^{T,T'}(\eta namk)dn = 
\int_{N_{Q}(\rmF)\backslash N_{Q}(\A)}
\Psi_{P,Q,\mathfrak{o}}^{T,T'}(na\eta mk)dn
\end{displaymath}
car $\eta \in M_{Q}(\rmF)$ 
normalise $N_{Q}(\A)$ sans changer sa mesure et 
il commute avec $A_{Q}^{\infty}$. 
Les facteurs 
de $\Psi_{P,Q,\mathfrak{o}}^{T,T'}(na\eta mk)$ 
($\eta \in (P(\rmF) \cap M_{Q}(\rmF))\backslash M_{Q}(\rmF)$) deviennent
\begin{gather*}
\Gamma'_{Q}(H_{Q}(na\eta mk)-T',T-T') = 
\Gamma'_{Q}(H_{Q}(a)-T',T-T'), \\
\htau_{P}^{Q}(H_{P}(na\eta mk)-T') =  
\htau_{P}^{Q}(H_{P}( \eta m)-T').
\end{gather*}
Quant \`a $k_{P,\mathfrak{o}}$,
en utilisant sa définition (\ref{eq:kPolDef}) 
et 
en faisant les changements de variables 
$(a\inv n\inv (U + \xi)na - \xi) \mapsto U$ et 
$m^{-1}\eta^{-1}U_{Q} \eta m \mapsto U_{Q}$ on obtient
\begin{multline}\label{induction}
\! \! \! \! e^{-2\rho_{\tlQ}(H_{Q}(a))} \! \int\limits_{K} \!
\int\limits_{\mathrlap{[N_{Q}]}}
k_{P,\mathfrak{o}}(na\eta mk)dn dk \! = \!
\int\limits_{K} \! \int\limits_{\mathrlap{\nl_{\tlP}(\A)}} \quad \sum_{\mathrlap{\xi \in \ml_{\tlP}(\rmF) \cap \mathfrak{o}}}
f(k\inv m\inv\eta\inv(\xi + U_{P})\eta mk)dU_{P}dk = \\
\int_{\nl_{\tlP}^{\tlQ}(\A)}\sum_{\xi \in \ml_{\tlP}(\rmF) \cap \mathfrak{o}}
\int\limits_{K}\int\limits_{\mathrlap{\nl_{\tlP}(\A)}} 
e^{2\rho_{\tlQ}(H_{Q}(m))}
f(k\inv (m\inv\eta\inv(\xi + U_{P} ^{Q})\eta m+ U_{Q})k)
dU_{Q}dkdU_{P}^{Q} = \\
\int_{\nl_{\tlP}^{\tlQ}(\A)}\sum_{\xi \in \ml_{\tlP}(\rmF) \cap \mathfrak{o}}
f_{Q}(m\inv\eta\inv(\xi + U_{P}^{Q})\eta m)dU_{P}^{Q}
\end{multline}
o\`u $f_{Q} \in \calS(\ml_{\tlQ}(\A))$ 
est définie par (\ref{eq:fQdef}) dans le 
paragraphe \ref{par:JTforLevis} et on a utilisé le 
fait que $\rho_{\tlQ}(H_{Q}(m)) = 0$ car $m \in M_{Q}(\A)^{1} \sbs M_{\tlQ}(\A)^{1}$.

En regardant la relation (\ref{eq:JMQisThis}), 
on s'aperçoit qu'avec la notation de l'équation 
(\ref{eq:useGammaJol}) on a
\begin{equation*}
\int_{Q(\rmF)\backslash U(\A)}
\sum_{P\subseteq Q}(-1)^{d_{P}^{Q}}
\sum_{\eta \in P(\rmF) \cap M_{Q}(\rmF)\backslash M_{Q}(\rmF)}
\Psi_{P,Q,\mathfrak{o}}^{T,T'}(\eta x)dx = 
J_{\ol}^{M_{Q}, T'}(f_{Q})p_{Q}(T, T')
\end{equation*}
où $p_{Q}(T, T')$ égale
\begin{equation*}
\int\limits_{\mathclap{A_{Q}^{\infty}}} \!
e^{2\rho_{\tlQ}(H_{Q}(a)) - 2\rho_{Q}(H_{Q}(a))}
\Gamma'_{Q}(H_{Q}(a)-T',T-T')da \! = \! \!
\int\limits_{\mathclap{\mathfrak{a}_{Q}}}
e^{\upla_{Q}(H)}
\Gamma'_{Q}(H-T_{Q}',T_{Q}-T_{Q}')dH.
\end{equation*}
Par un changement de variable on obtient 
$p_{Q}(T,T') = e^{\upla_{Q}(T_{Q}')}p_{Q}(T_{Q} - T_{Q}')$ 
où $p_{Q}$ c'est la fonction étudiée dans le lemme 
\ref{lem:pQExplicit}, d'où le résultat. 
\edem
\erop

En utilisant la proposition \ref{prop:mainQualitProp} démontrée 
ci-dessus et le lemme \ref{lem:pQExplicit} qui 
décrit les fonctions $p_{Q}$ explicitement on obtient 
le comportement en $T$ des distributions $J_{\ol}^{T}$ et $J^{T}$.

\begin{theo}\label{mainQualitThm}
Soit $f \in \calS(\tlul(\A))$.
Les fonctions 
$T \mapsto J^{T}_{\mathfrak{o}}(f)$ et
$T \mapsto  J^{T}(f)$ 
 où
$\mathfrak{o} \in \mathcal{O}$ et
$T$ 
parcourt $T_{+} + \mathfrak{a}_{0}^{+}$ 
sont des polynômes-exponentielles
dont les parties purement polynomiales sont constantes 
et données respectivement par  
\begin{gather*}
J_{\ol}(f) := \sum_{Q}(-1)^{d_{Q}}\hat \theta_{Q}(\upla_{Q})^{-1}
e^{\upla_{Q}(T_{Q}')}J_{\ol}^{M_{Q},T'}(f_{Q}), \\
J(f) := \sum_{Q}(-1)^{d_{Q}}\hat \theta_{Q}(\upla_{Q})^{-1}
e^{\upla_{Q}(T_{Q}')}J^{M_{Q},T'}(f_{Q}).
\end{gather*}
pour tout $T' \in T_{+} + \all_{0}^{+}$. 
En particulier, les distributions $J_{\ol}$ et $J$ 
ne dépendent pas de $T'$.
\end{theo}

\brem\label{rem:JMQisPolExp}
 Soit $Q$ un sous-groupe parabolique standard 
de $U$. 
Il résulte de la proposition \ref{prop:mainQualitProp} 
ci-dessus et du théorème 4.2 de \cite{chaud}
que les 
distributions $J_{\ol}^{M_{Q},T}$ et $J^{M_{Q},T}$, 
définies dans le paragraphe \ref{par:JTforLevis}, 
sont des polynômes-exponentielles avec des exposants 
$\upla_{R} - \upla_{Q}$ où $R$ parcourt 
les sous-groupes paraboliques contenus dans $Q$. 
Cependant, 
si $Q \neq U$ le terme purement polynomial n'est pas constant. 
\erem

\subsection{Invariance}\label{par:compConjug}

Dans ce paragraphe on démontrera 
l'invariance par conjugaison des distributions
$J_{\ol}$ et $J$. 

Soient $f \in \calS(\tlul(\A))$ et $y \in U(\A)$. 
Notons $f^{y} \in \calS(\tlul(\A))$ la fonction 
définie par $f^{y}(X) = f(\Ad(y)X)$. On voit donc 
qu'on a
\begin{equation*}
J_{\ol}^{T}(f^{y}) = 
\int_{U(\rmF)\backslash U(\A)}
\sum_{P}(-1)^{d_{P}}
\sum_{\delta \in P(\rmF)\backslash U(\rmF)}
k_{P,\ol}(\delta x)\htau_{P}(H_{P}(\delta xy)-T_{P})dx.
\end{equation*}
Pour $x \in U(\A)$ soit $k_{P}(x)$ un élément 
de $K$ tel que $xk_{P}(x)^{-1} \in P(\A)$. 
Alors, en utilisant l'égalité (\ref{GammaPrimeRecurrence}) 
on a
\begin{equation*}
\htau_{P}(H_{P}(\delta x y)-T_{P}) = 
\sum_{Q \supseteq P}(-1)^{d_{Q}}
\htau_{P}^{Q}(H_{P}(\delta x)-T_{P})
\Gamma_{Q}'(H_{P}(\delta x)-T_{Q},-H_{P}(k_{P}(\delta x)y))
\end{equation*}
d'où
\begin{equation*}
J_{\ol}^{T}(f^{y}) = 
\sum_{Q} \int_{Q(\rmF)\backslash U(\A)}
\sum_{P \subseteq Q}(-1)^{d_{P} - d_{Q}}
\sum_{\eta \in P(\rmF) \cap M_{Q}(\rmF) \backslash M_{Q}(\rmF)}
\Psi_{P,Q,\ol}^{T,y}(\eta x)dx
\end{equation*}
où
\begin{equation*}
\Psi_{P,Q,\ol}^{T,y}(x) = 
k_{P,\ol}(x)\htau_{P}^{Q}(H_{P}(x)-T_{P})
\Gamma_{Q}'(H_{P}(x)-T_{Q}, -H_{P}(k_{P}(x)y)).
\end{equation*}

Soit
$x = namk$ o\`u $n \in N_{Q}(\rmF)\backslash N_{Q}(\A)$, 
$m \in M_{Q}(\rmF)\backslash M_{Q}(\A)^{1}$, 
$a \in A_{Q}^{\infty}$ et $k \in K$. Donc 
$dx = e^{-2\rho_{Q}(H_{Q}(a))}dndadmdk$. 
On a pour $\eta \in M_{Q}(\rmF)$
\begin{equation*}
\Gamma_{Q}'(H_{P}(\eta namk)-T_{Q},-H_{P}(k_{P}(\eta namk)y))
 = 
\Gamma_{Q}'(H_{Q}(a)-T_{Q},-H_{Q}(ky)).
\end{equation*}
Ensuite, en faisant les mêmes opérations 
comme dans (\ref{induction}) au début de la preuve
de la proposition \ref{prop:mainQualitProp} on s'aperçoit  
qu'on a pour $P \subseteq Q$,
$m \in M_{Q}(\rmF)\backslash M_{Q}(\A)^{1}$ 
et $\eta \in M_{Q}(\rmF)$
\begin{multline*}
\! \! \! \! \! \! \!
\int\limits_{A_{Q}^{\infty}} \!
\int\limits_{K} \!
\int\limits_{\mathrlap{[N_{Q}]}}
e^{-2\rho_{Q}(H_{Q}(a))}
k_{\tlP,\ol}(\eta nmak)
\Gamma_{Q}'(H_{P}(\eta namk)-T_{Q},-H_{P}(k_{P}(\eta namk)y))
 dndkda \! = \\
 e^{\upla_{Q}(T_{Q})} \! \! \! \! \!
\int\limits_{\nl_{\tlP}^{\tlQ}(\A)}  \quad
\sum_{\mathclap{\xi \in \ml_{\tlP}(\rmF) \cap \ol}} \quad \ \ 
\int\limits_{K} \!
\int\limits_{\all_{Q}} \!
\int\limits_{\mathrlap{\nl_{\tlQ}(\A)}} \! \!
e^{\upla_Q(H)}
f(k^{-1}(m^{-1}\eta^{-1}(\xi + U_{P}^{Q})\eta m + U_{Q})k)
\Gamma_{Q}'(H,-H_{Q}(ky))
\\
dU_{Q}dHdkdU_{P}^{Q} = 
e^{\upla_Q(T_{Q})}
\int\limits_{\mathclap{\nl_{\tlP}^{\tlQ}(\A)}} \
\sum_{\xi \in \ml_{\tlP}(\rmF) \cap \ol}
f_{Q,y}(m^{-1}\eta^{-1}(\xi + U_{P}^{Q})m\eta)dU_{P}^{Q}
\end{multline*}
où l'on pose
\begin{equation*}
f_{Q,y}(X) = 
\int_{K}
\int_{\nl_{\tlQ}(\A)}
f(k^{-1}(X + U_{Q})k)u'_{Q}(k,y)
dU_{Q}dk, \quad X \in \ml_{\tlQ}(\A)
\end{equation*}
où 
\begin{equation*}
u'_{Q}(k,y) = \int_{\all_{Q}}
e^{\upla_Q(H)}\Gamma_{Q}'(H,-H_{Q}(ky))dH.
\end{equation*}
La fonction $k \mapsto u'_{Q}(k,y)$ étant continue 
on a bien $f_{Q,y} \in \calS(\ml_{\tlQ}(\A))$. 
On obtient le théorème suivant.

\begin{theo}\label{invarianceTheo} Soient $y \in U(\A)$ et 
$f \in \calS(\tlul(\A))$. Les distributions $J_{\ol}^{T}$ vérifient
\begin{equation*}
J_{\ol}^{T}(f^{y}) - J_{\ol}^{T}(f) = 
\sum_{P_{0} \subseteq Q \subsetneq U} 
e^{\upla_{Q}(T_{Q})}J_{\ol}^{M_{Q},T}(f_{Q,y})
\end{equation*}
où les distributions $J_{\ol}^{M_{Q},T}$ 
sur $\calS(\ml_{\tlQ}(\A))$ 
sont définies par (\ref{eq:indDistDef}).
 En particulier, on a
\[
J_{\ol}(f^{y}) = J_{\ol}(f), \quad J(f^{y}) = J(f).
\]

\bdem 
La formule pour la différence $J_{\ol}^{T}(f^{y}) - J_{\ol}^{T}(f)$ est claire 
après les calculs qu'on a faits. Cette formule-ci démontre aussi l'invariance, car
si $Q \subsetneq U$, d'après la remarque \ref{rem:JMQisPolExp}, le terme 
$J_{\ol}^{M_{Q},T}(f_{Q,y})$ est un polynôme-exponentielle d'exposants $\upla_{R} - \upla_{Q}$ 
où $R \sbs Q$. Il en découle que $e^{\upla_{Q}(T_{Q})}J_{\ol}^{M_{Q},T}(f_{Q,y})$ n'a pas de terme 
constant dans ce cas et par conséquent les termes constants de $J_{\ol}^{T}(f^{y})$ 
et $J_{\ol}^{T}(f)$ coïncident. 
\edem
\end{theo}

\subsection{Indépendance des choix}\label{par:noChoixMade}

Dans ce paragraphe on démontrera que la distribution $J_{\ol}$ 
ne dépend d'aucun choix, sauf le choix d'une mesure de Haar 
sur $U(\A)$ et les choix des mesures sur les 
sous-espaces $\calV$ de $\nl_{\tlzero}$, notre choix étant que 
$\calV(\rmF)\bsl \calV(\A)$ soit de volume $1$.

Démontrons d'abord que $J_{\ol}$ ne dépend pas du choix 
du sous-groupe parabolique minimal contenant $M_{0}$. 
Soient $P_{0}' \in \calP(M_{0})$ et $s \in \Omega$ 
tel que $P_{0}' = w_{s}^{-1}P_{0}w_{s}$. 
Notons $J_{P_{0}',\ol}^{T}$ et $J_{P_{0}',\ol}$ les distributions 
construites par rapport à $P_{0}'$.
C'est une conséquence simple de la relation 
(\ref{eq:htauSemist}) qu'on a pour $T \in T_{+} + \all_{0}^{+}$
\[
J_{\ol}^{T} = J_{P_{0}',\ol}^{s^{-1}T + H_{P_{0}}(w_{s}^{-1})}.
\]
En vertu du théorème \ref{mainQualitThm} on voit donc que 
les termes constants de $T \mapsto J_{\ol}^{T}$ et $T \mapsto J_{P_{0}',\ol}^{s^{-1}T + H_{P_{0}}(w_{s}^{-1})}$
coïncident, d'où
$J_{\ol} = J_{P_{0}',\ol}$.

On va démontrer maintenant que $J_{\ol}$ ne dépend pas du choix 
du sous-groupe compact maximal admissible par rapport à $M_{0}$.
Soit $K^{*}$ un autre tel sous-groupe. Pour tout sous-groupe parabolique 
standard $P$ notons 
$H_{P}^{*}$ le prolongement à $U(\A)$ par rapport à $K^{*}$ 
de la fonction $H_{P}$ définie 
sur $P(\A)$ par (\ref{harishChandra}). Comme avant, pour 
tout $x \in U(\A)$ on note $k_{P}(x)$ 
un élément de $K$ tel que $xk_{P}(x)^{-1} \in P(\A)$. 
On a alors $H_{P}^{*}(x) = H_{P}^{*}(xk_{P}(x)^{-1}) + H_{P}^{*}(k_{P}(x))$
et ni $H_{P}^{*}(k_{P}(x))$ ni $H_{P}^{*}(xk_{P}(x)^{-1})$ 
ne dépendent du choix de $k_{P}(x)$. De surcroît
$H_{P}^{*}(xk_{P}(x)^{-1}) = H_{P}(x)$.
Alors, en utilisant l'égalité (\ref{GammaPrimeRecurrence}) 
on a
\begin{equation}\label{gamPrRecApp}
\htau_{P}(H_{P}^{*}(x)-T) = 
\sum_{Q \supseteq P}(-1)^{d_{Q}}
\htau_{P}^{Q}(H_{P}(x)-T)
\Gamma_{Q}'(H_{Q}(x)-T,-H_{Q}^{*}(k_{P}(x))).
\end{equation}

Pour un sous-groupe parabolique standard $Q$ on pose
\begin{equation*}
f_{Q}^{*}(X) = 
\int_{K}
\int_{\nl_{\tlQ}(\A)}
f(k^{-1}(X + U_{Q})k)u^{*}_{Q}(k)
dU_{Q}dk, \quad X \in \ml_{\tlQ}(\A)
\end{equation*}
où 
\begin{equation*}
u^{*}_{Q}(k) = \int_{\all_{Q}}
e^{\upla_Q(H)}\Gamma_{Q}'(H,-H_{Q}^{*}(k))dH.
\end{equation*}
Donc $f_{Q}^{*} \in \calS(\ml_{\tlQ}(\A))$.

On note  $J_{K^{*},\ol}^{T}$ et $J_{K^{*},\ol}$, 
les distributions 
construites par rapport à $K^{*}$. 
En partant alors de l'égalité (\ref{gamPrRecApp}) et en effectuant 
les mêmes opérations que dans le paragraphe \ref{par:compConjug} on trouve
\begin{equation*}
J_{K^{*},\ol}^{T}(f) - J_{\ol}^{T}(f) = 
\sum_{P_{0} \subseteq Q \subsetneq U} 
e^{\upla_{Q}(T)}J_{\ol}^{M_{Q},T}(f_{Q}^{*}).
\end{equation*}
De nouveau, en vertu du théorème \ref{mainQualitThm} 
et de la remarque \ref{rem:JMQisPolExp}, par égalité des termes constants, on a
$J_{\ol} = J_{K^{*},\ol}$.

Par définition de $J_{\ol}$, il est clair qu'elle est indépendante des choix 
des mesures de Haar qu'on a fait dans le paragraphe \ref{par:haarMesSubs}, sauf pour la mesure de Haar sur $U(\A)$ et
les mesures sur les
sous-espaces de $\nl_{\tlzero}$. 

Il nous reste à démontrer l'indépendance du choix de sous-groupe de Levi minimal. 
Soit $M_{0}'$ un sous-groupe de Levi défini sur $\rmF$ minimal de $U$. 
Il existe alors un $y \in U(\rmF)$ tel que $M_{0}' = y^{-1}M_{0}y$. 
On note alors $J_{M_{0}', \ol}$ le terme constant 
de la distribution $J_{M_{0}',\ol}^{T}$ définie par rapport au sous-groupe 
parabolique minimal $y^{-1}P_{0}y$ et le compact maximal $y^{-1}Ky$. 
On trouve alors $J_{\ol}(f) = J_{M_{0}',\ol}(f^{y})$ et le résultat 
découle du théorème \ref{invarianceTheo}.

\subsection{Orbites semi-simples réguli\`eres}\label{par:regOrbsChap}

Soit $\mathcal{O}_{reg} \subseteq \mathcal{O}$ l'ensemble des 
orbites semi-simples réguli\`eres, \cad 
des orbites composées d'éléments 
semi-simples réguliers . Le but de ce paragraphe est 
de démontrer que, sous la condition 
que $\mathfrak{o} \in \mathcal{O}_{reg}$, la distribution 
$J_{\mathfrak{o}}(f)$ s'exprime comme une 
intégrale orbitale de $f$. 
On a d'abord:

\blem\label{regItIs}
 Soient $X \in \tlul(\rmF)$ un élément régulier semi-simple et 
$P$ un sous-groupe parabolique de $U$ différent de $U$. 
Alors $X \notin \ml_{\tlP}(\rmF)$.

\bdem 
Soit $P$ comme dans l'énoncé. 
Le tore $A_{P}$ est non-trivial et centralise  
$\ml_{\tlP}$ donc $X$ ne peut pas être dans $\ml_{\tlP}(\rmF)$
en vertu du point 
\textit{2)} de la proposition \ref{prop:relelt}.
\edem
\elem

Soient $f \in \calS(\tlul(\A))$ et $\ol \in \calO_{reg}$. 
D'après le lemme \ref{regItIs} ci-dessus,
on a $k_{f,P,\ol} \equiv 0$ 
pour tout sous-groupe parabolique standard différent de $U$.
En utilisant alors
le point \textit{2)} de la proposition \ref{prop:relelt}
 ainsi que la proposition \ref{prop:orbRegConj}, 
 on trouve
\[
k_{U,\mathfrak{o}}^{T}(x) = k_{U,\mathfrak{o}}(x)
= \sum_{\xi \in \tlul(\rmF) \cap \mathfrak{o}} 
f(x\inv \xi x) = 
\sum_{\delta \in U(\rmF)} 
f(x\inv\delta^{-1} X_{1}\delta x),
\]
où $X_{1} \in \ol$ quelconque. 
On obtient donc la proposition suivante:
\brop 
Pour tous  $f \in \calS(\tlul(\A))$, $T \in T_{+} + \all_{0}^{+}$, 
$\ol \in \calO_{reg}$
et $X_{1} \in\mathfrak{o}$ on a
\begin{displaymath}
J_{\mathfrak{o}}^{T}(f) = J_{\mathfrak{o}}(f) = 
\int_{U(\A)}f(x^{-1}X_{1}x)dx
\end{displaymath}
où l'intégrale est absolument convergente.
\erop 
\section{Formule des traces infinitésimale}\label{FourierTransChap}

Il résulte de l'analyse faite dans le paragraphe \ref{lesInvs} qu'on 
a la décomposition de $\tlul$ en sous-$\rmF$-espaces stables sous l'action de $U$ 
suivante:
\begin{equation}\label{eq:decBete}
\tlul= \tl_{1} \oplus \tl_{2} \oplus \tl_{3}
\end{equation}
où $\tl_{1} = \ul$, $\tl_{2} \cong \Res_{\rmE/\rmF}(V)$ et $\tl_{3} = \Lie(U(D_{0}, \tlPhi|_{D_{0}}))$. 
Soit $\tl \sbs \tlul$ un $\rmF$-sous-espace défini 
comme une somme directe de certains d'entre les $\tl_{i}$, $i = 1,2,3$. 
Il y a donc huit possibilités pour $\tl$. 
Puisque chaque $\tl_{i}$ est $U(\rmF)$-stable et la restriction 
de la forme bilinéaire $\bilif$, 
définie par (\ref{eq:kiliform}), 
à $\tl_{i}$ est non-dégénérée, il en est de même pour $\tl$. 
Pour $X \in \tlul(\A)$ soit $X_{\tl}$ la projection 
de $X$ à $\tl(\A)$ selon la décomposition 
(\ref{eq:decBete}) ci-dessus.

Fixons $\psi$ 
un caractère non-trivial de $\rmF \bsl \A$. 
Pour $\tl$ comme ci-dessus, notons $\calF_{\tl}$ 
l'opérateur sur $\mathcal{S}(\tlul(\A))$ suivant
\[
\calF_{\tl}(f)(X) = \int_{\tl(\A)}f(X - X_{\tl} + Y_{\tl})
\psi (\langle X_{\tl}, Y_{\tl} \rangle )dY_{\tl}, \quad 
f \in \mathcal{S}(\tlul(\A)), \ X \in \tlul(\A)
\]
où $dY_{\tl}$ c'est la mesure de Haar sur $\tl(\A)$ pour 
laquelle le volume de $\tl(\rmF) \bsl \tl(\A)$ vaut $1$. 

\begin{theo}\label{thm:RTFJRI} Pour tout
$f \in \mathcal{S}(\tlul(\A))$ on a
\begin{displaymath}
\sum_{\ol \in \calO} 
J_{\ol}(f) = 
\sum_{\ol \in \calO} 
J_{\ol}(\calF_{\tl}(f)).
\end{displaymath}

\bdem 
Soit $T \in T_{+} + \mathfrak{a}_{0}^{+}$. 
En utilisant l'identité (\ref{chgmtordre}) on a
que $J^{T}(f) - J^{T}(\calF_{\tl}(f))$ vaut
\begin{multline*}
\int\limits_{[U]}
F^{U}(x,T) (k_{U,U}(x,f) - k_{U,U}(x, \calF_{\tl}(f))) + \\
 \sum_{P_{1} \subsetneq P_{2}} 
\sum_{\delta \in P_{1}(\rmF) \backslash U(\rmF)}
\chi_{P_{1},P_{2}}^{T}(\delta x)
(k_{P_{1},P_{2}}(\delta x,f) - k_{P_{1},P_{2}}(\delta x,\calF_{\tl}(f)))dx
\end{multline*}
où pour $P_{1} \sbs P_{2}$ et $\upphi \in \mathcal{S}(\tlul(\A))$ 
on pose
\[
k_{P_{1},{P_2}}(x, \upphi) = 
\sum_{P_{1} \subseteq P \subseteq P_{2}}
(-1)^{d_{P}}k_{P,\upphi}(x), \ 
x \in P_{1}(\rmF)\backslash U(\A).
\]
On a alors pour tout $x \in U(\A)$
\begin{displaymath}
k_{U,U}(x,f) =
\sum_{\xi \in \tlul(\rmF)}f(x\inv \xi x) = 
\sum_{\xi \in \tlul(\rmF)}\calF_{\tl}(f)(x\inv \xi x) = 
k_{U,U}(x, \calF_{\tl}(f)) 
\end{displaymath}
grâce \`a la formule sommatoire de Poisson. 
Fixons $\varepsilon_{0} >0$.
En utilisant le théorème \ref{thm:MainConv2}
pour $f$ et $\calF_{\tl}(f)$ on a pour 
tout $N > 0$
\begin{displaymath}
|J^{T}(f) - J^{T}(\calF_{\tl}(f))|= O(e^{-N\|T\|})
\end{displaymath}
si $T \in T_{+} + \mathfrak{a}_{0}^{+}$ est tel 
que $\forall \al \in \Delta_{0}$, $\al(T) > \varepsilon_{0}\|T\|$. 
D'après la proposition \ref{prop:mainQualitProp}, 
la différence
$J^{T}(f) - J^{T}(\calF_{\tl}(f))$ 
égale une constante plus une somme 
de polynômes-exponentielles en $T$ qui, 
en vertu du lemme \ref{lem:positiveLinear},  
tendent vers $\infty$ quand la norme de 
$T \in \mathfrak{a}_{0}^{+}$ tend vers $\infty$. L'égalité ci-dessus 
implique alors $J^{T}(f) = J^{T}(\calF_{\tl}(f))$ donc en particulier 
$J(f) = J(\calF_{\tl}(f))$ ce qu'il fallait démontrer. 
 \edem
\end{theo}
\section{Orbites relativement semi-simples régulières}\label{sec:orbRrssS}

\newcommand{\IoB}{I_{0}}
\newcommand{\IoBpm}{I_{0}}
\newcommand{\IoA}{I'}
\newcommand{\IoApm}{I'}
\newcommand{\sbse}{\subseteq_{\epsilon}}
\newcommand{\sbspm}{\subseteq_{\pm}}
\newcommand{\acalI}{|\mathcal{I}|}
\newcommand{\acalK}{|\mathcal{K}|}
\newcommand{\acalA}{|\mathcal{A}|}
\newcommand{\acalIB}{|\mathcal{C}_{B}|}
\newcommand{\calIB}{\mathcal{C}_{B}}
\newcommand{\acalJ}{|\mathcal{J}|}
\newcommand{\acalJj}{|\mathcal{J}_{1}|}
\newcommand{\acalJd}{|\mathcal{J}_{2}|}
\newcommand{\acalJt}{|\mathcal{J}_{3}|}
\newcommand{\acalJc}{|\mathcal{J}_{4}|}
\newcommand{\acalJp}{|\mathcal{J}_{5}|}
\newcommand{\acalIj}{|\mathcal{I}_{1}|}
\newcommand{\acalId}{|\mathcal{I}_{2}|}
\newcommand{\acalIt}{|\mathcal{I}_{3}|}
\newcommand{\acalIc}{|\mathcal{I}_{4}|}
\newcommand{\acalIp}{|\mathcal{I}_{5}|}
\newcommand{\acalZ}{|\mathcal{Z}|}
\renewcommand{\upla}{\uprho}

Soit
\(\mathcal{O}_{rs} \subseteq \mathcal{O}\) l'ensemble 
des classes contenant un élément \(
\begin{pmatrix}
B & u \\
u^{\sharp} & d \\
\end{pmatrix}\) tel que le polynôme caractéristique de 
$B$ est séparable (i.e. $B$ est semi-simple régulier dans $\ul(\rmF)$). 
On appelle de telles 
classes \textit{relativement semi-simples régulières}. 
Si \(\mathfrak{o} \in \mathcal{O}_{rs}\), alors tout élément 
\(
\begin{pmatrix}
B_{0} & u_{0} \\
u_{0}^{\sharp} & d \\
\end{pmatrix} \in \mathfrak{o}\) a la propriété 
que $B_{0}$ soit semi-simple régulier.

Soit $\ol \in \mathcal{O}_{rs}$.
Le but de cette section est de donner une expression explicite 
pour $J_{\ol}$ ce qu'on achève par le théorème
\ref{thm:theThmOrbs}. Voici le plan de la section:
après avoir introduit quelques notations dans le paragraphe 
suivant \ref{par:noteEnsem}, on décrit 
la décomposition de $\ol$ en $U(\rmF)$-orbites 
dans le paragraphe \ref{par:orbitesDansClasse}. 
On introduit encore 
un peu plus de notation dans 
\ref{par:defsOrbs}. Dans le paragraphe \ref{par:LeResultRss}
on définit une expression $j_{\ol}(x)$ pour laquelle on a
\begin{equation}\label{eq:JolIsjol}
J_{\ol}(f) = \int_{U(\rmF)\bsl U(\A)} j_{\ol}(x)dx.
\end{equation}
En supposant cela, on
donne la preuve du théorème \ref{thm:theThmOrbs} omettant 
les preuves des énoncés techniques. 
Dans la section \ref{par:noyuTronqNouv} on introduit 
un nouveau noyau tronqué $j_{f,\ol}^{T}(x)$ 
tel que 
$\int_{[U]} j_{f,\ol}^{T}(x)dx = J_{\ol}^{T}(f)$. 
Ce résultat nous permet de démontrer (\ref{eq:JolIsjol}) 
 dans le paragraphe suivant 
\ref{par:IntReprDeJ}. 
Dans le paragraphe \ref{par:resDeConv} on démontre les résultats de convergence nécessaires et on finit la preuve dans le 
dernier paragraphe \ref{par:holoRes} 
où on étudie certaines fonctions zêtas.

\subsection{Notations}\label{par:noteEnsem}

On utilisera les lettres $I, J$, avec de possibles indices, 
pour noter des 
sous-ensembles finis de $\N^{*}$. 
Soit $I \sbs \N^{*}$ fini. On pose 
$-I = \bigcup_{i \in I} \{-i\}$.
On dit que $\calI$ est un $\epsilon$-sous-ensemble 
de $I$, si $\calI \sbs I \cup -I$ et si pour tout $i \in \calI$ 
on a $-i \nin \calI$.
Dans ce cas on 
écrit $\calI \sbse I$. La notation est 
un peu abusive car $\calI$ n'est pas forcément un 
sous-ensemble de $I$. 
On définit aussi 
$\acalI = \{|i| | i \in \calI \} \sbs \N^{*}$ 
et
$\calI^{\sharp} \sbse I$ comme $I^{\sharp} = -\calI$.
On réserve les lettres $\calI$, $\calJ$ et $\calK$ et seulement ces 
trois lettres avec de possibles indices, 
pour des $\epsilon$-sous-ensembles. 

On utilisera aussi la notation abrégée suivante
soit $I' \sbs \N^{*}$ et $\calI, \calJ \sbse I'$, 
on écrira $\calJ \cup \calI \sbse I'$ pour signifier 
que la réunion ensembliste $\calI \cup \calJ$ est aussi 
un $\epsilon$-sous-ensemble (ce qui n'est pas toujours vrai). 
On utilise le symbole $\sqcup$ pour noter la réunion
disjointe, donc $I \sqcup J = I'$ implique 
$I \cap J = \varnothing$.
On écrira aussi, pour $\calI \sbse I'$ fixés
\[
\sum_{\acalJ = I'} := 
\sum_{\begin{subarray}{c}
\calJ \sbse I' \\
\acalJ = I'
\end{subarray}}
, \ 
\sum_{\calK \sqcup \calJ \sbs \calI} := 
\sum_{\begin{subarray}{c}
\calK, \calJ \sbs \calI \\
\calK \cap \calJ = \varnothing
\end{subarray}}
, \ 
\sum_{\acalK \sqcup \acalJ = I'} := 
\sum_{\begin{subarray}{c}
\calK, \calJ \sbse I' \\
\acalK \sqcup \acalJ = I'
\end{subarray}}.
\]

\subsection{Orbites dans une classe relativement semi-simple régulière}\label{par:orbitesDansClasse}

Dans ce paragraphe on décrit la décomposition en orbites 
d'une classe relativement semi-simple régulière.

D'abord, on rappelle la description des $U(\rmF)$-orbites semi-simples régulières
dans $\ul(\rmF)$. Donnons-nous
\begin{itemize}
\item Un ensemble fini  $I = \{1,\ldots, m\} \sbs \N$. 
\item Pour tout $i \in I$ une extension 
finie $\rmF_{i}$ de $\rmF$ et on note 
$\rmE_{i}$ la $\rmE$-algèbre étale
$\rmF_{i} \otimes_{\rmF} \rmE$. 

\item Pour $i \in I$, 
on note $\sigma_{i}$ l'automorphisme de $\rmE_{i}$
défini comme  $\Id_{\rmF_{i}} \otimes \sigma$. On se donne 
alors des éléments $b_{i} \in \rmE_{i}$ et $c_{i} \in \rmF_{i}^{*}$ tels 
que $\sigma_{i}(b_{i}) = -b_{i}$.
\item On note 
$\rmE_{I} := \bigoplus_{i \in I}\rmE_{i}$ et 
$b_{I} = (b_{i})_{i \in I} \in \rmE_{I} $. 
Alors on suppose que $b_{I}$ engendre la $\rmE$-algèbre étale 
$\rmE_{I}$. 

\item $\dim_{\rmE}V = \sum_{i \in I}[\rmE_{i}:\rmE]$. 

\item Pour tout $i \in I$ on définit une forme hermitienne 
non-dégénérée 
sur le $\rmE$-espace $\rmE_{i}$ 
relativement à l'extension $\rmE/\rmF$ par
$\Phi_{i}(x_{i},y_{i}) = \Tr_{\rmE_{i}/\rmE}
(\sigma_{i}(x_{i})y_{i} c_{i})$, $x_{i},y_{i} \in \rmE_{i}$. Cela induit la forme 
$\Phi_{I} = \sum_{i \in I} \Phi_{i}$ sur $\rmE_{I}$.

\item Pour $w \in \rmE_{I}$ et $i \in I$ on note $w_{i}$ 
la projection orthogonale de $w$ à $\rmE_{i}$. 
On définit $B_{I} \in \Lie (U(\rmE_{I}, \Phi_{I}))$ par 
$B_{I}(w) = \sum_{i \in I}b_{i}w_{i} = b_{I}w$. 
On voit alors que $B_{I}$ est semi-simple régulier
et on a $\rmE_{I} \cong \rmE[X]/P_{B_{I}}(X)$ où 
$P_{B_{I}} \in \rmE[X]$ c'est le polynôme caractéristique de $B_{I}$. 
\end{itemize}

Supposons qu'on a un isomorphisme
\begin{equation}\label{iotaIso}
\iota : \rmE_{I} \irar V
\end{equation}
de $\rmE$-espaces hermitiens. Il induit naturellement 
des isomorphismes, qu'on note 
aussi $\iota$, 
de groupes et d'algèbres de Lie correspondants.
Dans ce cas, la $U(\rmF)$-classe de conjugaison 
de $\iota(B_{I}) \in \ul(\rmF)$ ne dépend pas de l'isomorphisme choisi. 
Toute 
$U(\rmF)$-orbite semi-simple régulière dans $\ul(\rmF)$ est obtenue par une 
telle construction. On note alors simplement $B = \iota(B_{I})$. 

On considère $\iota$ fixé, il nous donne l'action de $U(\rmF)$ et 
$\ul(\rmF)$
sur $\rmE_{I}$. Explicitement, pour $g \in U(\rmF)$ 
et $v \in \rmE_{I}$ on écrira simplement $gv$ 
pour $\iota^{-1}(g)v$, de même pour les éléments de $\ul(\rmF)$. En particulier, 
on écrit $Bv$ pour désigner $B_{I}v$.

Soit $I_{1} \sbs I$ l'ensemble des 
$i \in I$ tels que $\rmE_{i}$ n'est pas 
un corps. Pour $i \in I_{1}$ on a alors que 
$\rmF_{i}$ est une extension de $\rmE$. On se donne une 
 inclusion $\iota_{i}: \rmE \hookrightarrow \rmF_{i}$. 
 Soit $\rmF_{-i}$ la $\rmE$-algèbre 
 dont le groupe multiplicatif et additif égale celui de $\rmF_{i}$ 
 et où l'action de $\rmE$ est donné par l'inclusion 
 $\rmE \ni e \mapsto \iota_{i}(\sigma(e)) \in \rmF_{-i}$.
On a alors une $\rmE$-algèbre étale $\rmF_{-i} \oplus \rmF_{i}$ 
munie de l'action de 
$\Gal(\rmE/\rmF)$ donnée par
$\sigma((a,b)) = (b,a)$ et 
l'on a une forme hermitienne 
\[
\Phi_{\rmF_{-i} \oplus \rmF_{i}}(x_{-i} + x_{i}, x_{-i}' + x_{i}') = 
\Tr_{\rmF_{i}/\rmE}( c_{i} x_{-i}x_{i}') + 
\Tr_{\rmF_{-i}/\rmE}(c_{i} x_{-i}'x_{i})
\]
où $x_{-i},x_{-i}'\in \rmF_{-i}$ et 
$x_{i},x_{i}'\in \rmF_{i}$.
Dans ce cas, on a l'isomorphisme de
$\rmE$-algèbres 
$\rmE_{i} = \rmF_{i} \otimes_{\rmF} \rmE \cong \rmF_{-i} \oplus 
\rmF_{i}$ 
donné par $x \otimes e \mapsto (x\sigma(e), xe)$ 
préservant toutes les structures mentionnées. 
On fixe cet isomorphisme. 
Les $\rmE$-sous-espaces $\rmF_{-i}$ et 
$\rmF_{i}$ de $\rmE_{i}$ sont des sous-espaces isotropes maximaux 
et ce sont
les seuls sous-espaces non-triviaux 
stables par $b_{i}$ et donc par $B$.

Quitte à conjuguer $B$ on peut supposer 
que $B \in \ml_{P_{I_{1}}}(\rmF)$ pour un sous-groupe parabolique 
standard, noté $P_{I_{1}}$, et qu'aucun $M_{P_{I_{1}}}(\rmF)$-conjugué 
de $B$ n'appartienne à $\ml_{R}(\rmF)$ 
pour un sous-groupe parabolique $R \sbn P_{I_{1}}$. 
Alors, si $P_{I_{1}}$ est défini comme 
le stabilisateur du drapeau isotrope 
$\{0\}  = V_{i_{0}} \subseteq  \cdots \subseteq V_{i_{l}}$ 
on a $\# I_{1}= l$ et donc, quitte à réindexer, 
on peut supposer que 
$I_{1} = \{1,2, \ldots, l\}$. 
Finalement, quitte à conjuguer $\iota$ 
par un élément de groupe de Weyl de $U$, 
on peut supposer que 
$\iota$
induit des isomorphismes 
$V_{i_{j}}/V_{i_{j-1}} \cong \rmF_{-j}$ 
et $V^{\sharp}_{i_{j}}/V^{\sharp}_{i_{j-1}} \cong \rmF_{j}$ 
pour $j \in I_{1}$ et 
$Z_{i_{l}} \cong \prod_{j \in I \smin I_{1}} \rmE_{j}$. 

Pour tout $\calI \sbse I_{1}$ on note 
$\rmF_{\calI} = \prod_{i \in \calI} \rmF_{i}$ 
et $1_{\calI}$ l'unité de $\rmF_{\calI}^{*}$. Alors, quand 
$\calI$ parcourt les $\epsilon$-sous-ensembles de $I_{1}$, 
les $\rmF_{\calI}$ parcourent tous les sous-espaces isotropes 
de $\rmE_{I}$ stables par $B$. 
On note aussi pour tout $I' \sbs I$, $\rmE_{I'} = 
\prod_{i \in I'} \rmE_{i}$.

Fixons une fois pour toutes une classe $\mathfrak{o} \in \mathcal{O}_{rs}$
telle qu'il existe un  $X \in \ol$ de type  $\matx{B}{u}{u^{\sharp}}{d}$. 
Remarquons que $d$ ne dépend que de $\ol$, on va le noter $d_{\ol}$.
Tout élément de $\ol$ est donc $U(\rmF)$-conjugué à un élément 
du type 
$\matx{B}{u'}{(u')^{\sharp}}{d_{\ol}} \in \ol$. 
On se propose de décrire les classes de 
$U(\rmF)$-conjugaison dans $\ol$.

Notons $T_{I}$ le centralisateur de $B$ dans $U$. C'est un sous-tore maximal de $U$.
Pour 
tout $I' \sbs I$ soit $T_{I'}$ le plus grand sous-tore 
de $T_{I}$ qui agit trivialement sur $\rmE_{I \smin I'}$. 
On a donc que $T_{I'}(\rmF)$ s'identifie 
aux éléments $ u \in \rmE_{I'}^{*}$ 
tels que $u_{i}\sigma_{i}(u_{i})$ c'est l'élément neutre de $\rmE_{i}^{*}$ pour tout $i \in \rmI'$. Si $I' = \{i\}$ on écrit simplement 
$T_{i}$. En particulier, si $i \in I_{1}$ 
on a $T_{i}(\rmF) = \{(t_{-i}, t_{i}) \in \rmF_{-i}^{*} \times \rmF_{i}^{*} | 
t_{-i} = t_{i}^{-1} \} \cong \rmF_{i}^{*}$.

On introduit l'ensemble 
$V_{\ol} = \{u \in \rmE_{I}|
-\Phi_{I}(u,B^{i-1}u) = A_{i}(X) \ \forall \ 1 \le i \le n\}$  où
$X \in \mathfrak{o}$ quelconque et les invariants $A_{i}(X)$ 
ont été définis dans le paragraphe \ref{lesInvs}.
Comme $T_{I}(\rmF)$ agit sur $\rmE_{I}$, commute à $B$ et laisse $\Phi$ invariant, il agit aussi sur $V_{\ol}$.
On voit que l'ensemble des orbites dans $V_{\ol}$ sous l'action 
de $T_{I}(\rmF)$
est en bijection avec l'ensemble des classes de $U(\rmF)$-conjugaison dans 
$\mathfrak{o}$, la bijection étant induite 
par l'application
$V_{\ol} \ni u \mapsto \matx{B}{u}{u^{\sharp}}{d_{\ol}}$.

\blem\label{lem:chaudNagginLem} 
Il existe un $\al_{I} \in \rmF_{I}$ tel que 
pour tout $u \in \rmE_{I}$
on a
\[
u \in V_{\ol} \Longleftrightarrow 
\sigma_{i}(u_{i})u_{i} = \al_{i} \ \forall \ i \in I.
\]

\bdem 
Pour tous $u, u' \in V_{\ol}$ et $k \in \N$ on a
\[
\Phi(u,B^{k}u) = 
\sum_{i \in I}\Tr_{\rmE_{i}/\rmE}(\sigma_{i}(u_{i})u_{i}b_{i}^{k}c_{i}) = 
\sum_{i \in I}\Tr_{\rmE_{i}/\rmE}(\sigma_{i}(u_{i}')u_{i}'b_{i}^{k}c_{i}) =\Phi(u',B^{k}u')
\]
d'où
\[
\Tr_{\rmE_{I}/ \rmE}
\dsl b^{k}_{I} \dsl\sum_{i \in I}
c_{i}(\sigma_{i}(u_{i})u_{i} - \sigma_{i}(u_{i}')u_{i}') \rb \rb = 0, \quad 
\forall \ k \in \N.
\]
D'après la proposition (18.3) dans \cite{bInvolutions}, 
la forme
$\Tr_{\rmE_{I}/ \rmE}$ est non-dégénérée, et puisque
les puissances de $b_{I}$ 
engendrent $\rmE_{I}$ sur $\rmE$ on obtient
\begin{equation}
\sigma_{i}(u_{i})u_{i} = \sigma_{i}(u_{i}')u_{i}' \ \forall \ i \in I, \ 
\forall \ u,u' \in V_{\ol}.
\end{equation}
On pose donc $\al_{i} = \sigma_{i}(u_{i})u_{i}$ où $u \in V_{\ol}$ 
quelconque. Il reste à démontrer 
que si $u \in \rmE_{I}$ est tel que 
$\sigma_{i}(u_{i})u_{i} = \al_{i}$ pour tout $i \in I$ 
alors $u_{i} \in V_{\ol}$. Pour cela il suffit de faire le même calcul dans le sens inverse.
\edem
\elem

Soit $\al_{I} \in \rmE_{I}$ comme dans le lemme précédant. 
Notons $I_{2} \sbs I$ l'ensemble de $i \in I$ tels 
que $\al_{i} = 0$ et posons $I_{0} = I_{1} \cap I_{2}$.

\brop\label{prop:orbitsRRSS}
Il existe une unique $T_{I}(\rmF)$-orbite 
dans $V_{\ol}$ composée des $u \in V_{\ol}$ 
tels que $u_{i} = 0$ pour tout $i \in I_{2}$. On choisit 
$\xi_{\varnothing}$ 
un représentant de cette orbite. Alors,  
les $T_{I}(\rmF)$-orbites dans $V_{\ol}$ 
sont en bijection avec les $\epsilon$-sous-ensembles de $I_{0}$, 
le représentant de l'orbite correspondant à $\calI \sbse I_{0}$ 
étant 
$\xi_{\calI} := \xi_{\varnothing} + 1_{\calI}$. Les orbites 
de dimension maximale correspondent aux $\calI \sbse I_{0}$ 
tels que $\acalI = I_{0}$.

\bdem 
On voit que $u$ et $u'$
sont $T_{I}(\rmF)$-conjugués 
si et seulement $u_{i}$ et $u_{i}'$ sont $T_{i}(\rmF)$-conjugués 
pour tout $i \in I$. 

Soient $u, u' \in V_{\ol}$ et $i \in I \smin I_{0}$. 
On prétend que $u_{i}$ et $u_{i}'$ sont $T_{i}(\rmF)$-conjugués. 
En effet, si $i \in I_{2} \smin I_{0}$ on a $u_{i} = u_{i}' = 0$ 
d'après le lemme \ref{lem:chaudNagginLem}.
Sinon, en vertu de ce lemme-là, 
on a $u_{i}, u'_{i} \in \rmE_{i}^{*}$ 
et $t_{i} := u_{i}/u_{i}'$ appartient à 
$T_{i}(\rmF)$ et vérifie $t_{i}u_{i} = u_{i}'$.

Soit $i \in I_{0}$ et $u \in V_{\ol}$. On écrira $u_{i}^{-i}$ (resp. $u_{i}^{i}$) 
la projection de $u_{i} \in \rmE_{i}$ à $\rmF_{-i}$ (resp. $\rmF_{i}$).
On a alors $u_{i}\sigma(u_{i}) = (u_{i}^{-i}u_{i}^{i}, u_{i}^{-i}u_{i}^{i})= 0$, 
on voit alors qu'au moins l'un de $u_{i}^{-i}$, $u_{i}^{i}$
vaut zéro. 
Pour finir, il suffit de montrer que si 
$\calI \sbse I_{0}$ et 
$u, u' \in V_{\ol}$ sont tels 
que pour $i \in I_{0} \cup -I_{0}$ 
on a 
$u_{|i|}^{i} \in \rmF_{i}^{*}$ si et seulement si 
$ i \in \calI$ (de même pour $u'$) alors 
$u$ et $u'$ sont $T_{I}(\rmF)$ conjugués et ne le sont pas sinon. 
Cela est clair, car pour $i \in I_{0}$, 
$T_{i}(\rmF)$ agit simplement transitivement sur $\{0\} \times \rmF_{i}^{*} \sbs \rmE_{i}$ et sur
$\rmF_{-i}^{*} \times \{0\} \sbs \rmE_{i}$ 
et $0 \in \rmE_{i}$ est $T_{i}(\rmF)$-conjugué 
à lui même seulement. 
\edem
\erop

\subsection{Quelques définitions associées aux orbites}\label{par:defsOrbs}

D'après la proposition \ref{prop:orbitsRRSS} ci-dessus les  
$X_{\calI} := 
\matx{B}{\xi_{\varnothing} +1_{\calI}}{(\xi_{\varnothing} +1_{\calI})^{\sharp}}{d_{\ol}}$,
où $\calI \sbse I_{0}$, sont des représentants des orbites 
pour l'action de $U(\rmF)$ à $\ol$. On les considère fixés désormais.

Pour un $\rmF$-sous-groupe $H$ de $U$, $X \in \tlul(\rmF)$ 
et une $\rmF$-algèbre $R$, notons $H(R,X)$ le groupe des $R$-points 
du stabilisateur de $X$ dans $H$.
Alors, pour tout $\calI \sbse I_{0}$ et toute $\rmF$-algèbre $R$ on a:
\begin{equation}\label{eq:stabXcalI}
U(R,X_{\calI}) = T_{I_{2} \smin \acalI}(R).
\end{equation}

Pour tout $I' \sbs I_{1}$ soit $M_{I'}$ le sous-groupe de 
Levi de $U$ défini comme le stabilisateur des espaces $\rmF_{-i}$ 
et $\rmF_{i}$ pour tout $i \in I'$. Donc, en utilisant la notation du paragraphe \ref{par:orbitesDansClasse}, 
on a en particulier $M_{I_{1}} = M_{P_{I_{1}}}$. On pose 
$\all_{I'} = \all_{P_{I'}}$ et $\all_{I', \C}^{*} = \all_{I'}^{*} \otimes_{\R} \C$ 
où $P_{I'} \in \calP(M_{I'})$ quelconque. Donc, on a $\all_{I'} \sbs \all_{I_{1}}$.
Soient $I'' \sbs I'$,
$Q \in \calP(M_{I''})$ et $P \in \calP(M_{I'})$ tels que $Q \sps P$, alors 
$\all_{Q} = \all_{I''}$, $\all_{P} = \all_{I'}$ et $\all_{P}^{Q} = \all_{I' \smin I''}$.
Grâce à la décomposition (\ref{eq:decomp}), on obtient alors 
$\all_{I'} = \all_{I''} \oplus \all_{I' \smin I''}$. 
Notons finalement $\la_{I'}$ la projection d'un 
$\la \in \all_{I_{1},\C}^{*}$  à 
$\all_{I',\C}^{*}$.

Soit $A_{\IoB}$ le sous-tore de $T_{\IoB}$, déployé sur $\rmF$ et maximal pour cette propriété. 
Pour $i \in I_{0} \cup -I_{0}$,
soit $\upla_{i} \in \all_{\IoB}^{*}$ le caractère par lequel 
$A_{\IoB}$ 
agit sur $\rmF_{i}$ (l'inclusion $A_{\IoB} \hrar M_{\IoB}$ induit l'isomorphisme 
$\Hom_{\rmF}(A_{\IoB}, \Gm) \otimes_{\Z} \R \cong \all_{\IoB}^{*}$).
On a donc $\upla_{i} = - \upla_{-i}$. 
Soit $\calI \sbse \IoBpm$.
On a que $\{\upla_{i}\}_{i \in \calI}$ est une base de 
$\all_{\acalI}^{*}$. 
On pose 
\[
\upla_{\calI} = \sum_{i \in \calI}\upla_{i}.
\]
Il est facile de voir que si  $Q \in \calF(M_{\IoB})$ est tel que 
$V_{Q} = \rmF_{\calI}$ alors 
$\underline{\rho}_{Q} \in \all_{\IoB}^{*}$, défini à la fin du paragraphe 
\ref{par:haarMesSubs}, égale $\upla_{\calI}$.

Soient $\{e_{i}^{\vee}\}_{i \in I_{0} \cup -I_{0}} \sbs \all_{\IoB}$
les vecteurs tels que
\[
\upla_{j}(e_{i}^{\vee}) = 
\begin{cases} 
1 \text{ si } j=i,\\
-1 \text{ si } j= -i,\\
0 \text{ sinon, }
\end{cases} \quad i,j \in I_{0} \cup -I_{0}.
\]
Il est clair que pour tous $\calI \sbse \IoBpm$ 
et $I' \sbs \IoB$ la projection de
$\sum_{i \in \calI} a_{i} e_{i}^{\vee}$ à $\all_{I'}$ 
égale 
$\sum_{i \in \calI, \ |i| \in I'} a_{i} e_{i}^{\vee}$. 

 On introduit aussi
le cône ouvert $\all_{\calI}^{*}$ défini comme
\[
\all_{\calI}^{*} = 
\displaystyle \left \{
\sum_{i \in \calI} a_{i}\upla_{i}| a_{i} > 0
\right \} 
\sbs \all_{\acalI}^{*} \sbs \all_{\IoB}^{*}.
\]
Donc, pour qu'un $\la \in \all_{\acalI,\C}^{*}$ vérifie
$\Rel(\la) \in \all_{\calI}^{*}$ il faut et suffit que
$\Rel (\la(e_{i}^{\vee})) >0$ pour tout $ i \in \calI$.

Pour $\calI \sbse \IoBpm$ 
soit 
$\ind_{\calI}$ la fonction caractéristique de 
$H \in \all_{I_{0}}$ tels que 
\[
\upla_{i}(H) \le 0, \  \forall \ i \in \calI \cap I_{0} \quad 
\text{et} \quad \upla_{i}(H) < 0, \  \forall \ i \in \calI \cap -I_{0}.
\] 
On a alors pour tous $I' \sbs I_{0}$ et $H \in \all_{I_{0}}$
\begin{equation}\label{eq:decOf1}
1 = \sum_{\acalI =  I'}\ind_{\calI}(H).
\end{equation}
Pour $\calJ_{1} \sqcup \calJ_{3} \sbse I_{0}$ 
et $\calJ_{2} \sbs \calJ_{3}$ on utilisera la notation suivante
\[
\calJ_{13} := \calJ_{1} \cup \calJ_{3}, \ 
\calJ_{3 \smin 2} := \calJ_{3} \smin \calJ_{2}.
\]

\subsection{Le résultat principal}\label{par:LeResultRss}

Dans ce paragraphe on énonce et on démontre le théorème \ref{thm:theThmOrbs}. 
Cependant, certains résultats seront seulement énoncés avec les renvois  
vers leurs démonstrations dans les paragraphes suivants. 

Il est clair que si l'orbite d'un $X_{\calI}$ où $\calI \sbse I_{0}$ 
intersecte non-trivialement 
$\ml_{\tlP}(\rmF)$ où $P$ est un sous-groupe parabolique standard
de $U$
alors celui-ci est conjugué 
à un élément de $\calF(M_{I_{1}})$. Pour tout
$Q \in \calF(M_{I_{1}})$ soit $\calI_{Q} \sbse I_{1}$ l'unique 
$\epsilon$-ensemble 
tel que $V_{Q} = \rmF_{\calI_{Q}}$. 
\blem\label{orbsAsSets}
 Soient $Q \in \calF(M_{I_{1}})$ et $\calI \sbse \IoB$. Alors $X_{\calI} \in \ml_{\tlQ}(\rmF)$ 
si et seulement si $\calI_{Q} \sbse \IoB$ et $\acalI \cap |\calI_{Q}| = \varnothing$.
\elem
\bdem
En vertu du lemme \ref{shortLemmeOrb} on a, 
en utilisant la notation du paragraphe \ref{generalities}, 
que $X_{\calI} \in \ml_{\tlQ}(\rmF)$ si et seulement si 
$\xi_{\varnothing} + 1_{\calI} \in Z_{Q}$. 
Or, $\xi_{\varnothing} + 1_{\calI} \in Z_{Q}$ si et seulement si 
$\xi_{\varnothing} \in Z_{Q}$ et $\acalI \cap |\calI_{Q}| = \varnothing$. 
En outre, la première condition est équivalente à dire que
pour tout $i \in \calI_{Q}$ la $i$-composante de $\xi_{\varnothing}$ 
vaut zéro, donc $\calI_{Q} \sbse I_{0}$, d'où le résultat. 
\edem

On rappelle que, avec la notation de la fin du paragraphe
\ref{par:prelimstrace}, si $P_{1}$, $P$ sont des sous-groupes paraboliques standards 
alors pour $s \in \Omega(\all_{1}, P)$ on note 
$s^{-1}P$ le sous-groupe parabolique semi-standard $Q \in \calF(M_{1})$ 
égale $w_{s}^{-1}Pw_{s}$.
On démontre maintenant:
\blem\label{lem:orbitCapP} 
 Soient $P$ un sous-groupe parabolique standard de $U$ et 
$\calI \sbs_{\epsilon} I_{0}$. 
Alors, l'intersection de la $U(\rmF)$-orbite de $X_{\calI}$ 
avec $\ml_{\tlP}(\rmF)$ égale
\[
\coprod_{
\begin{subarray}{c}
s \in \Omega(\all_{I_{0}}, P) \\
\acalI \cap |\calI_{s^{-1}P}| = \varnothing
\end{subarray}
}
\coprod_{\eta \in M_{P}(\rmF, \Ad(w_{s})X_{\calI}) \bsl M_{P}(\rmF)}
\{\Ad(\eta^{-1}w_{s})X_{\calI}\}.
\]
\bdem
Supposons que 
$\Ad(\gamma^{-1})X_{\calI} \in \ml_{\tlP}(\rmF)$ 
pour un $\gamma \in U(\rmF)$. En particulier donc, 
grâce au lemme \ref{shortLemmeOrb}, 
on a $\Ad(\gamma^{-1})B \in \ml_{P}(\rmF)$. 
En raisonnant comme 
au début du paragraphe 5.2 dans \cite{chaud} on voit 
qu'il existe un unique $s \in \Omega(\mathfrak{a}_{I_{1}},P)$ 
(voir le paragraphe \ref{par:prelimstrace}) 
et un unique $\delta \in M_{P}(\rmF, \Ad(w_{s})X_{\calI}) \bsl M_{P}(\rmF)$
tels que 
$\gamma^{-1} = \delta^{-1} w_{s}\). 
Soit $Q = s^{-1}P$, alors  
$X_{\calI} \in \ml_{\tlQ}(\rmF)$. 
En vertu du lemme \ref{orbsAsSets} ci-dessus on a $Q \in \calF(M_{\IoB}, P)$.
 Donc, en utilisant la bijection 
(\ref{eq:semiStConjBij}) et l'unicité de $s$ on obtient 
$s \in \Omega(\mathfrak{a}_{I_{0}},P)$ d'où le résultat.
\edem
\elem

Pour une fonction $\upphi$ sur $\tlul(\A)$ et 
$x \in U(\A)$ on définit
\[
\upphi_{x}(X) := \upphi(\Ad(x^{-1})X).
\]

Fixons $f \in \calS(\tlul(\A))$.
Soit \(P\) un sous-groupe parabolique semi-standard de \(U\).
D'après le lemme \ref{shortLemmeOrb} \textit{b)} 
on a un isomorphisme $N_{P}$-équivariant
\begin{equation}\label{lisomorphism}
\nl_{\tlP} \cong \nl_{P} \oplus \Res_{\rmE/\rmF}(V_{P}).
\end{equation}
On pose 
\begin{equation}\label{eq:fhatP}
f^{P}(X) = \int_{V_{P}(\A)}f(X + Y_{P})dY_{P}, \quad 
X \in \tlul(\A).
\end{equation}
Supposons en plus que $P$ est standard et posons
\begin{equation}\label{eq:IP}
I_{P,\ol}(x) = 
\sum_{\xi \in \ml_{\tlP}(\rmF) \cap \ol}
\sum_{\eta \in N_{P}(\rmF)}f_{\eta x}^{P}(\xi), \quad 
x \in P(\rmF)\bsl U(\A),
\end{equation}
où $f_{x}^{P} = (f_{x})^{P}$. 
En vertu du lemme \ref{lem:orbitCapP} on a alors
\[
I_{P,\ol}(x) = \sum_{s \in \Omega(\all_{I_{0}}, P)}
\sum_{
\begin{subarray}{c}
\calI \sbs_{\epsilon}\IoBpm \\
\acalI \cap |\calI_{s^{-1}P}| = \varnothing
\end{subarray}}
\sum_{\eta \in P(\rmF, \Ad w_{s} X_{\calI}) \bsl P(\rmF)}
f_{\eta x}^{P}(\Ad (w_{s})X_{\calI}).
\]
En utilisant
\[
U(\rmF,X)\backslash  U(\rmF) \xrightarrow \sim 
U(\rmF,\Ad w_{s}X) \backslash  U(\rmF) \text{ grâce \`a }  
\eta_{0} \mapsto w_{s}\eta_{0},
\]
on a
\begin{equation}\label{IPrewritten1} 
\sum_{\delta \in P(\rmF) \bsl U(\rmF)}
I_{P,\ol}(\delta x)
= 
\sum_{s \in \Omega(\all_{I_{0}}, P)}
\sum_{
\begin{subarray}{c}
\calI \sbs_{\epsilon}\IoBpm \\
\acalI \cap |\calI_{s^{-1}P}| = \varnothing
\end{subarray}} 
\sum_{\delta \in U(\rmF,X_{\calI}) \bsl U(\rmF)}
f_{\delta x}^{s^{-1}P}(X_{\calI}).
\end{equation}
Fixons un caractère additif continu non-trivial $\psi$ sur $\rmF \bsl \A $.
Pour une fonction $\upphi \in \calS(\tlul(\A))$ et 
$\calJ \sbse I_{0}$ on définit la transformée 
de Fourier de $\upphi$ par rapport à $\calJ$
\[
\hat \upphi^{\calJ}(X) = 
\int_{\A_{\calJ}}
\upphi \dsl 
\matx{B'}{u'(1-1_{\calJ^{\sharp}}) + u_{\calJ}}{(u'(1-1_{\calJ^{\sharp}}) + u_{\calJ})^{\sharp}}{d'}
\rb \psi ( \langle u_{\calJ}, u'1_{\calJ^{\sharp}} \rangle )du_{\calJ}
\]
où $X = \matx{B'}{u'}{(u')^{\sharp}}{d'} \in \tlul(\A)$, $1$ 
c'est l'unité dans $\rmE_{I}^{*}$, $du_{\calJ}$ c'est la mesure 
de Haar sur $\A_{\calJ} := \rmF_{\calJ} \otimes_{\rmF} \A$
pour laquelle $\rmF_{\calJ} \bsl \A_{\calJ}$ est de volume $1$ et 
$\bilif$ c'est l'accouplement défini par (\ref{eq:kiliform}). 
En fait, avec nos identifications, 
pour tout $u, u' \in V(\A) = \rmE_{I} \otimes_{\rmF} \A$ 
on a $\langle u, u' \rangle = - \Phi_{I}(u,u') - \Phi_{I}(u',u)$.

Soient alors $s \in \Omega(\all_{I_{0}}, P)$ et $\calI \sbse I_{0}$ tel que $\acalI \cap |\calI_{s^{-1}P}| = \varnothing$. 
On a $f_{x}^{s^{-1}P}(X_{\calI}) = \widehat{(f_{x})}^{\calI_{s^{-1}P}}(X_{\calI})$. 
On utilisera la notation abrégée  $\hf_{x}^{\calI_{s^{-1}P}} := \widehat{(f_{x})}^{\calI_{s^{-1}P}}$ 
(à ne pas confondre avec $(\widehat{f}^{\calI_{s^{-1}P}})_{x}$ qu'on n'utilisera pas).
En utilisant (\ref{eq:stabXcalI}), on voit donc qu'on peut réécrire (\ref{IPrewritten1}) comme
\begin{equation}\label{eq:Ipolrewritten}
\sum_{s \in \Omega(\all_{I_{0}}, P)}
\sum_{
\begin{subarray}{c}
\calI \sbs_{\epsilon} \IoBpm \\
\acalI \cap |\calI_{s^{-1}P}| = \varnothing
\end{subarray}} 
\sum_{\delta \in T_{I_{2} \smin \acalI}(\rmF) \bsl U(\rmF)}
\hf_{\delta x}^{\calI_{s^{-1}P}}(X_{\calI}).
\end{equation}
Les sommes \eqref{IPrewritten1} 
et \eqref{eq:Ipolrewritten} 
sont convergentes grâce au lemme suivant. 
 \blem[cf. \ref{lem:forFucksSake}] 
 Soient $\calI, \calJ \sbse \IoB$ tels que $\acalI \cap \acalJ = \varnothing$. 
 Alors, pour tout $x \in U(\A)$ on a:
 \[
\sum_{\delta \in U(X_{\calI}, \rmF) \bsl U(\rmF)} |\hf^{\calJ}_{\delta x}(X_{\calI})| < \infty.
 \]
 \elem 

On pose:
\[
j_{\ol}(x) =
\sum_{P \sps P_{0}}(-1)^{d_{P}}
\sum_{\delta \in P(\rmF) \bsl U(\rmF)}
I_{P,\ol}(\delta x), \quad x \in U(\rmF) \bsl U(\A).
\]

\brop[cf. \ref{prop:newjol}]\label{prop:newjol0} On a
\[
\int_{U(\rmF) \bsl U(\A)}|j_{\ol}(x)|dx < \infty
\quad \text{et} \quad \int_{U(\rmF) \bsl U(\A)}j_{\ol}(x)dx = J_{\ol}(f).
\]
\erop
\nident
En utilisant le résultat ci-dessus,
on a $\int_{[U]}j_{\ol}(x)dx = J_{\ol}(f)$ où, grâce à
la formule (\ref{eq:Ipolrewritten}), on a
\[
j_{\ol}(x) = \sum_{P}(-1)^{d_{P}}
\sum_{s \in \Omega(\all_{I_{0}}, P)}
\sum_{
\begin{subarray}{c}
\calI \sbs_{\epsilon} \IoBpm \\
\acalI \cap |\calI_{s^{-1}P}| = \varnothing
\end{subarray}} 
\sum_{\delta \in T_{I_{2} \smin \acalI}(\rmF) \bsl U(\rmF)}
\hf_{\delta x}^{\calI_{s^{-1}P}}(X_{\calI}).
\]
En inversant l'ordre de sommation on a aussi
\[
j_{\ol}(x) =
\sum_{\calI \sqcup \calJ \sbse I_{0}}
\ \ \ \
\mu_{\calJ}
\sum_{\mathclap{\delta \in T_{I_{2} \smin \acalI}(\rmF) \bsl U(\rmF)}}
\ \ 
\hf^{\calJ}_{\delta x}(X_{\calI})
\]
où 
\[
\mu_{\calJ} = \sum_{P \sps P_{0}}
(-1)^{d_{P}}\sum_{\begin{subarray}{c}
s \in \Omega(\all_{I_{0}}, P)\\
\calI_{s^{-1}P} = \calJ
\end{subarray}
}1.
\]
\blem\label{lem:signCombi} Soit $\calJ \sbse I_{0}$. Alors $\mu_{\calJ} = (-1)^{\# \calJ}$.
\bdem
Posons $m = \# \calJ$ et
\[
a_{k}^{m} = \#\{(\calJ_{0}, \calJ_{1}, \ldots, \calJ_{k})| 
\varnothing = \calJ_{0} \sbn \calJ_{1} \sbn \cdots \sbn 
\calJ_{k} = \calJ\}, \quad k \in \N.
\]
Il est clair que $a_{k}^{m}$ ne dépend que de $m = \# \calJ$. 
On voit, en invoquant la bijection \eqref{eq:semiStConjBij} du paragraphe \ref{par:prelimstrace}, 
qu'on a $\mu_{\calJ} = \sum_{k=0}^{m}(-1)^{k}a_{k}^{m}$. 
Il est facile de voir qu'on a la relation de récurrence suivante:
\[
a_{k}^{m} = \sum_{1 \le j \le m} {{m}\choose{j}}a^{m-j}_{k-1}, \quad 
k \in \N^{*}.
\]
En utilisant cette relation on vérifie $\mu_{\calJ} = (-1)^{m}$ par récurrence.
\edem
\elem

On vient d'obtenir alors
\begin{equation}\label{eq:jIknewAllAlong}
j_{\ol}(x) = 
\sum_{\calI \sqcup \calJ \sbse I_{0}}
(-1)^{\# \calJ}
\sum_{\mathclap{\delta \in T_{I_{2} \smin \acalI}(\rmF) \bsl U(\rmF)}}
\ \ 
\hf^{\calJ}_{\delta x}(X_{\calI}).
\end{equation}

\blem[cf. \ref{lem:UpsPros}]\label{lem:UpsPros0}
 Soient $\calJ, \calJ_{1}, \calJ_{2} \sbse I_{0}$ 
tels que $\calJ_{1} \sqcup \calJ_{2} \sbs \calJ$ 
et $\acalJ = I_{0}$.
L'intégrale suivante
\begin{equation*}
\brLa_{\calJ_{1},\calJ_{2}}^{\calJ}(f)(\la) = 
\int\limits_{\mathclap{
T_{I_{2} \smin |\calJ_{12}|}(\rmF)\bsl U(\A)}}
\ind_{(\calJ \smin \calJ_{2}) \cup \calJ_{2}^{\sharp}}(H_{0}(x))
e^{\la(H_{0}(x))}
\hf_{x}^{\calJ \smin \calJ_{1}}(X_{\calJ_{1} \cup \calJ_{2}^{\sharp}})dx, \quad 
\la \in \all_{I_{0}, \C}^{*}
\end{equation*}
converge absolument et uniformément sur tous les compacts 
d'un ouvert de $\all_{I_{0},\C}^{*}$ contenant $0$
et
admet un prolongement méromorphe
à $\all_{I_{0},\C}^{*}$, noté aussi
$\brLa_{\calJ_{1},\calJ_{2}}^{\calJ}(f)$.
\elem
\noindent Ici $H_{0}$ c'est juste $H_{P_{0}}$. Comme on l'a déjà remarqué dans le paragraphe \ref{par:defsOrbs}, 
on a que $\all_{I_{0}}$ est contenu dans $\all_{I_{1}} = \all_{P_{I_{1}}}$ qui est contenu 
dans $\all_{0}$ car $P_{I_{1}}$ est un sous-groupe parabolique standard. 

On aura besoin du lemme suivant:
\blem\label{lem:PoissonInlcExcl} Soient $\uphi \in \calS(\tlul(\A))$ et $\calI \sbse I_{0}$. Alors, 
pour tout $x \in U(\A)$ on a:
\begin{gather}
\sum_{\delta \in T_{\acalI}(\rmF)}\uphi_{\delta x}(X_{\calI}) = 
\sum_{\calI_{2} \sbs \calI_{1} \sbs \calI}(-1)^{\#(\calI \smin \calI_{1})}
\sum_{\delta \in T_{|\calI_{2}|}(\rmF)}
\huphi_{\delta x}^{\calI_{1}}(X_{\calI_{2}^{\sharp}}), \label{eq:PoissonInlcExcl1}
\\
\sum_{\delta \in T_{\acalI}(\rmF)}\huphi^{\calI}_{\delta x}(X_{\calI^{\sharp}}) = 
\sum_{\calI_{2} \sbs \calI_{1} \sbs \calI}(-1)^{\#(\calI \smin \calI_{1})}
\sum_{\delta \in T_{|\calI_{2}|}(\rmF)}
\huphi_{\delta x}^{\calI \smin \calI_{1}}(X_{\calI_{2}}) \label{eq:PoissonInlcExcl2},
\end{gather}
où, comme toujours, $\huphi_{x}^{\calI'} = \widehat{(\uphi_{x})}^{\calI'}$ pour 
$\calI' \sbse I_{0}$.
\bdem
Démontrons d'abord l'égalité (\ref{eq:PoissonInlcExcl1}). En utilisant 
le principe d'inclusion-exclusion on a
\[
\sum_{\delta \in T_{\acalI}(\rmF)}\uphi_{\delta x}(X_{\calI}) = 
\sum_{\calI_{1} \sbs \calI}(-1)^{\#(\calI \smin \calI_{1})}
\sum_{\calI_{2} \sbs \calI_{1}}\sum_{\delta \in T_{|\calI_{2}|}(\rmF)}
\uphi_{\delta x}(X_{\calI_{2}}).
\]
D'autre part, pour tout $\calI_{1} \sbs \calI$, la formule sommatoire de Poisson nous donne
\[
\sum_{\calI_{2} \sbs \calI_{1}}\sum_{\delta \in T_{|\calI_{2}|}(\rmF)}
\uphi_{\delta x}(X_{\calI_{2}}) = 
\sum_{\calI_{2} \sbs \calI_{1}}\sum_{\delta \in T_{|\calI_{2}|}(\rmF)}
\huphi_{\delta x}^{\calI_{1}}(X_{\calI_{2}^{\sharp}})
\]
d'où (\ref{eq:PoissonInlcExcl1}). L'égalité (\ref{eq:PoissonInlcExcl2}) se démontre de la même façon. 
\edem
\elem

\blem\label{lem:distIsUpsilon} On a 
\[
J_{\ol}(f) = \sum_{\acalJ = I_{0}}
\sum_{\calJ_{1} \sqcup \calJ_{2} \sbs \calJ}
(-1)^{\# (\calJ \smin \calJ_{12})}
\brLa_{\calJ_{1},\calJ_{2}}^{\calJ}(f)(0).
\]

\bdem 
En utilisant l'égalité (\ref{eq:PoissonInlcExcl2}) 
du lemme \ref{lem:PoissonInlcExcl} on a pour tout $x \in U(\A)$
\begin{multline}\label{eq:disUpsSom}
\sum_{\acalJ = I_{0}}
\sum_{\calJ_{1} \sqcup \calJ_{2} \sbs \calJ}
(-1)^{\# (\calJ \smin \calJ_{12})}
\sum_{\eta \in T_{|\calJ_{12}|}(\rmF)}
\hf_{\eta x}^{\calJ \smin \calJ_{1}}
(X_{\calJ_{1} \cup \calJ_{2}^{\sharp}}) 
\ind_{(\calJ \smin \calJ_{2}) \cup \calJ_{2}^{\sharp}}(H_{0}(\eta x))
= \\
\sum_{\acalJ = I_{0}}
\sum_{\calJ_{1} \sqcup \calJ_{2} \sbs \calJ}
\sum_{\calJ_{4} \sbs \calJ_{3} \sbs \calJ_{2}}
(-1)^{\#(\calJ \smin \calJ_{13})}
\sum_{\eta \in T_{|\calJ_{14}|}(\rmF)}
\hf_{\eta x}^{\calJ \smin \calJ_{13}}
(X_{\calJ_{14}})
\ind_{(\calJ \smin \calJ_{2}) \cup \calJ_{2}^{\sharp}}(H_{0}(\eta x)).
\end{multline}
On change l'ordre de sommation en mettant 
$\calJ_{0} := \calJ \smin \calJ_{13}$ et $\calI_{0} = \calJ_{14}$.
On obtient alors
\[
\sum_{\calI_{0} \sqcup \calJ_{0} \sbse I_{0}}
(-1)^{\# \calJ_{0}}
\sum_{\eta \in T_{|\calI_{0}|}(\rmF)}
\hf_{\eta x}^{\calJ_{0}}
(X_{\calI_{0}}) \dsl \sum_{\calK} \ind_{\calK}(H_{0}(\eta x)) \rb
\]
où la somme porte sur les $\calK \sbse I_{0}$ 
tels que  $\hf_{\eta x}^{\calJ_{0}}
(X_{\calI_{0}})\ind_{\calK}(H_{0}(\eta x))$ 
apparaît dans la somme (\ref{eq:disUpsSom}) ci-dessus. 
On prétend que tout $\acalK = I_{0}$ est de cette forme exactement une fois pour $\calI_{0}$ et $\calJ_{0}$ fixés. Fixons-les et 
soit $\calK \sbse I_{0}$ tel que $\acalK = I_{0}$. 
On vérifie alors que pour 
\begin{gather*}
\calJ_{1} := \calK \cap \calI_{0}, \ 
\calJ_{2} := (\calK \smin (\calJ_{0} \cup \calI_{0}))^{\sharp}, \ 
\calJ_{3} := (\calK \smin (\calJ_{0} \cup \calJ_{0}^{\sharp} \cup  \calI_{0}))^{\sharp}, \\
\calJ_{4} := \calK^{\sharp} \cap \calI_{0}, \ 
\calJ := (\calK \cap (\calJ_{0} \cup \calI_{0}))  
\cup (\calK \smin (\calJ_{0} \cup \calI_{0}))^{\sharp}
\end{gather*}
l'expression 
$\hf_{\eta x}^{\calJ_{0}}(X_{\calI_{0}})\ind_{\calK}(H_{0}(\eta x))$ 
apparaît dans la somme (\ref{eq:disUpsSom}).
Inversement, il apparaît une seule fois car $\calK$ détermine 
les ensembles $\calJ_{1}$, $\calJ_{2}$, $\calJ_{3}$, $\calJ_{4}$ 
et $\calJ$
uniquement. 

En utilisant l'identité (\ref{eq:decOf1}) on voit qu'on 
a démontré
\begin{multline*}
\sum_{\acalJ = I_{0}}
\sum_{\calJ_{1} \sqcup \calJ_{2} \sbs \calJ}
(-1)^{\# (\calJ \smin \calJ_{12})}
\sum_{\mathclap{\eta \in T_{|\calJ_{12}|}(\rmF)}}
\hf_{\eta x}^{\calJ \smin \calJ_{1}}
(X_{\calJ_{1} \cup \calJ_{2}^{\sharp}}) 
\ind_{(\calJ \smin \calJ_{2}) \cup \calJ_{2}^{\sharp}}(H_{0}(\eta x))
=\\
\sum_{\mathclap{\calI \sqcup \calJ \sbse I_{0}}}
(-1)^{\# \calJ}
\sum_{\mathclap{\eta \in T_{|\calI|}(\rmF)}}
\hf_{\eta x}^{\calJ}
(X_{\calI}).
\end{multline*}
En regardant le côté droit de cette égalité et
en utilisant la formule (\ref{eq:jIknewAllAlong}), 
on voit que, en vertu de la proposition \ref{prop:newjol0}, 
l'intégrale de cette expression sur 
$T_{I_{2}}(\rmF) \bsl U(\A)$ égale $J_{\ol}(f)$.
Or, le lemme 
\ref{lem:UpsPros0} dit que 
l'intégrale 
du côté gauche sur le même quotient 
donne le résultat cherché. 
\edem
\elem

On introduit maintenant les fonctions zêta.

\brop[cf. \ref{prop:zetaDefProps}]\label{prop:zetaDefProps1}
Soit $\calJ \sbse I_{0}$.
Alors, l'intégrale 
\[
\vol(T_{I_{2} \smin I_{0}}(\rmF) \bsl T_{I_{2} \smin I_{0}}(\A))
\int_{U(\A,X_{\calJ}) \bsl U(\A)}f(x^{-1}X_{\calJ}x)
e^{\la(H_{0}(x))}dx, \quad \la \in \all_{\acalJ,\C}^{*}
\]
converge absolument 
et uniformément sur tous les compacts d'un 
ouvert non-vide de $\all_{\acalJ,\C}^{*}$ 
et admet un prolongement méromorphe 
à $\all_{\acalJ,\C}^{*}$, noté $\zeta_{\calJ}(f)$. Les 
fonctions $\zeta_{\calJ}(f)$ et 
$\brLa_{\calJ_{1}, \calJ_{2}}^{\calJ}(f)$ vérifient la relation suivante:
\[ 
\sum_{\acalJ = I_{0}}
\sum_{\calJ_{1} \sqcup \calJ_{2} \sbs \calJ}
(-1)^{\#(\calJ \smin \calJ_{12})}
\brLa_{\calJ_{1}, \calJ_{2}}^{\calJ}(f) = 
\sum_{\acalJ = I_{0}}\zeta_{\calJ}(f).
\]
\erop
\bdem
La première partie est démontré dans la proposition \ref{prop:zetaDefProps}. 
La deuxième assertion c'est le lemme \ref{lem:sumZetasUpsilon} ci-dessous.
\edem

On est prêt à démontrer le résultat principal de cette section.

\begin{theo}\label{thm:theThmOrbs} 
Pour tout $f \in \calS(\tlul(\A))$ la somme 
$\sum_{\acalJ = I_{0}} \zeta_{\calJ}(f)$ est holomorphe en $\la = 0$ et l'on a:
\[
J_{\ol}(f) = \dsl \sum_{\acalJ = I_{0}} \zeta_{\calJ}(f)\rb (0).
\]

\bdem 
Le théorème découle du résultat 
d'holomorphie donné dans le lemme \ref{lem:UpsPros0} 
et de l'égalité démontrée dans le 
lemme \ref{lem:distIsUpsilon} accouplée avec l'égalité de 
la proposition \ref{prop:zetaDefProps1}.
\edem
\end{theo}

\subsection{Deuxième formule pour le noyau tronqué}\label{par:noyuTronqNouv}

\blem\label{keyBJ}
 Soient \(P\) un sous-groupe parabolique 
de \(U\),
\(X\in \ml_{\tlP}(\rmF) \cap \mathfrak{o}\) et
$p : \nl_{\tlP} \rightarrow \nl_{P}$ 
la projection donnée par l'isomorphisme 
(\ref{lisomorphism}). 
L'application suivante:
\begin{equation*}
N_{P}(R) \ni \eta  \mapsto p(\Ad(\eta^{-1})X - X) \in \nl_{P}(R)
\end{equation*}
est une bijection entre $N_{P}(R)$ et $\nl_{P}(R)$ pour toute \(\rmF\)-alg\`ebre \(R\).

\bdem  
Soit \(X = 
\begin{pmatrix}
B & u \\
u^{\sharp} & d \\
\end{pmatrix}\) la décomposition de $X$ 
comme dans (\ref{xisamatrix}). 
On voit alors que l'application décrite dans le lemme 
c'est juste 
$N_{P}(R) \ni \eta \mapsto \Ad(\eta^{-1})B-B 
\in \nl_{P}(R)$. 
Le fait que cette application est une bijection 
pour $B$ régulier semi-simple c'est le contenu 
du lemme 2.3 de \cite{chaud}. 
\edem
\elem

On considère $f \in \calS(\tlul(\A))$ fixée.
Dans le paragraphe \ref{par:LeResultRss}, équation (\ref{eq:fhatP}), 
on a introduit la fonction $f^{P}$ ainsi que $f_{x}(X) := f(\Ad(x^{-1})X)$. 
Par $f_{x}^{P}$ on entend toujours $(f_{x})^{P}$.

Le corollaire suivant découle directement du lemme \ref{keyBJ}.
\bcor\label{cor:lemcor1}
 Soient \(P\) un sous-groupe parabolique standard de \(U\),
 \(\xi \in  \ml_{\tlP}(\rmF) \cap \mathfrak{o}\) et $x \in U(\A)$, alors
\begin{displaymath}
\sum_{\eta \in N_{P}(\rmF)} f_{\eta x}^{P}(\xi)= 
\sum_{\zeta \in \nl_{P}(\rmF)}
f_{x}^{P}(\xi +\zeta).
\end{displaymath}
\ecor

\bcor\label{cor:lemcor2}
 Soient \(P\) un sous-groupe parabolique standard de \(U\), 
\(\xi \in  \ml_{\tlP}(\rmF) \cap \mathfrak{o}\) et $x \in U(\A)$. Alors
\begin{displaymath}
\int_{N_{P}(\A)} 
f_{nx}^{P}(\xi)dn = \int_{\nl_{\tlP}(\A)}
f_{x}(\xi+U_{P})dU_{P}.
\end{displaymath}

\bdem 
Le corollaire découle du corollaire 2.5 de 
\cite{chaud} et du fait que \(N_{P}(\A)\) normalise 
\(V_{P}(\A)\) sans changer la mesure de Haar. 
\edem
\ecor

Pour $T \in \all_{0}^{+}$ posons
\begin{displaymath}
j_{\mathfrak{o}}^{T}(x) = j_{f,\mathfrak{o}}^{T}(x) = 
\sum_{P \sps P_{0}}(-1)^{d_{P}}
\sum_{\delta \in P(\rmF)\backslash U(\rmF)}
\hat \tau_{P}(H_{P}(\delta x)-T)I_{P,\mathfrak{o}}(\delta x), \ 
x \in U(\rmF)\backslash U(\A)
\end{displaymath}
où $I_{P, \ol}$ est définie dans la section \ref{par:LeResultRss} 
par (\ref{eq:IP}).
La fonction $j_{\mathfrak{o}}^{T}$ 
est une variante de $k_{\ol}^{T}$ définie 
au début de la section \ref{sec:convergence}. 


\begin{theo}\label{thm:convThmNotMain}
 Pour tout \(T \in T_{+} + \mathfrak{a}_{0}^{+}\) on a
\begin{displaymath}
\int_{U(\rmF)\backslash U(\A)}|j_{\mathfrak{o}}^{T}(x)|dx < \infty.
\end{displaymath}

\bdem 
En procédant comme au début 
de la preuve du théorème (\ref{thm:MainConv}) 
et en utilisant la notation du début de ce théorème, on montre 
que l'intégrale
$\int_{[U]}|j_{\mathfrak{o}}^{T}(x)|dx$
est majorée par la somme sur les sous-groupes 
paraboliques standards \(P_{1} \subseteq S \subseteq P_{2}\) 
de \(U\) de
\[
\int_{P_{1}(\rmF) \bsl U(\A)}
\chi^{T}_{1,2}(x)
\sum_{\xi \in (\ml_{\tlS, \tlone})'(\rmF) \cap \mathfrak{o}}
|\sum_{S \subseteq P \subseteq P_{2}}
(-1)^{d_{P}} \sum_{\zeta \in \nl_{\tlS}^{\tlP}(\rmF)}
\sum_{\eta \in N_{P}(\rmF)}f_{\eta x}^{P}(\xi+\zeta)| dx.
\]
En utilisant le corollaire \ref{cor:lemcor1} 
et ensuite la formule sommatoire de 
Poisson on s'aperçoit que la somme entre la valeur absolue dans l'intégrale ci-dessus égale
\begin{equation}\label{newKernelConv}
\sum_{S \subseteq P \subseteq P_{2}}
(-1)^{d_{P}}
\sum_{\zeta_{1} \in \bar \nl_{\tlS}^{\tlP}(\rmF)}
\sum_{\zeta_{2} \in \bar \nl_{P}(\rmF)}
\widetilde{\phi}_{S}(x,\xi,\zeta_{1}+\zeta_{2}) 
\end{equation}
o\`u
\begin{displaymath}
\widetilde{\phi}_{S}(x,X,Y) = 
\int_{\nl_{\tlS}(\A)}f_{x}(X+U_{S})\psi(\langle U_{S},Y\rangle)dU_{S}, 
\ x \in U(\A), X \in \ml_{\tlS}(\A), 
Y \in \bar \nl_{\tlS}(\A).
\end{displaymath}
Soit $P$ un sous-groupe parabolique 
contenu entre \(S\) et 
\(P_{2}\). Pour $P \sbs R \sbs P_{2}$ notons
\begin{displaymath}
\bar \nl_{R,2}' = 
(\bar \nl_{R}^{2} \cap (\bar \nl_{\tlR}^{\tltwo})') 
\oplus \bar \nl_{2} = 
\bar \nl_{R} 
\smallsetminus \dsl \bigcup_{R \subseteq Q \subsetneq P_{2}}
\bar \nl_{R}^{Q} \oplus \bar \nl_{2}\rb.
\end{displaymath}
On a donc
la décomposition $\bar \nl_{P} = \coprod_{P \subseteq  R \subseteq P_{2}}
\bar \nl_{R,2}'$
et en utilisant la décomposition de $\bar \nl_{\tlS}^{\tlP}(\rmF)$ 
donnée par
(\ref{locFermLin}), on voit que 
l'expression (\ref{newKernelConv}) égale
\begin{equation*}
\begin{split}
\sum_{S \subseteq Q \subseteq P\subseteq R \subseteq P_{2}}
(-1)^{d_{P}}
\sum_{\zeta_{1} \in (\bar \nl_{\tlS}^{\tlQ})'(\rmF)}
\sum_{\zeta_{2} \in \bar \nl_{R,2}'(\rmF)}
\widetilde{\phi}_{S}(x,\xi,\zeta_{1}+\zeta_{2}) & = \\
\sum_{S \subseteq Q \subseteq R \subseteq P_{2}}
\sum_{\zeta_{1} \in (\bar \nl_{\tlS}^{\tlQ})'(\rmF)}
\sum_{\zeta_{2} \in \bar \nl_{R,2}'(\rmF)}
\widetilde{\phi}_{S}(x,\xi,\zeta_{1}+\zeta_{2})
(\sum_{Q \subseteq P \subseteq R}(-1)^{d_{P}}) &= \\
\sum_{S \subseteq R \subseteq P_{2}}
(-1)^{d_{R}}
\sum_{\zeta_{1} \in (\bar \nl_{\tlS}^{\tlR})'(\rmF)}
\sum_{\zeta_{2} \in \bar \nl_{R,2}'(\rmF)}
\widetilde{\phi}_{S}(x,\xi,\zeta_{1}+\zeta_{2})&,
\end{split}
\end{equation*}
où l'on a utilisé l'identité (\ref{basicidentity}) dans 
la derni\`ere égalité. 

Pour \(P_{1} \sbs S \subseteq  P_{2}\), posons
\begin{displaymath}
\Psi_{S}(x,X,Y) = 
\sum_{\zeta_{2} \in \bar \nl_{2}(\rmF)}
\widetilde{\phi}_{S}(x,X,Y+\zeta_{2}), \quad 
x \in U(\A), X \in \ml_{\tlS}(\A), 
Y \in \bar \nl_{\tlS}^{\tltwo}(\A).
\end{displaymath}
Alors \(\Psi_{S}(x,X,Y) \in 
\mathcal{S}((\ml_{\tlS} \oplus \bar \nl_{\tlS}^{\tltwo})(\A))\)
pour un \(x\) fixé et on a pour $S \sbs R \sbs P_{2}$:
\begin{displaymath}
\sum_{\xi \in \ml_{\tlS, \tlone}'(\rmF) \cap \mathfrak{o}}
|\sum_{\zeta_{1} \in (\bar \nl_{\tlS}^{\tlR})'(\rmF)}
\sum_{\zeta_{2} \in \bar \nl_{R,2}'(\rmF)}
\widetilde{\phi}_{S}(x,\xi,\zeta_{1}+\zeta_{2})| \le 
\sum_{\xi \in (\ml_{\tlone}^{\tlS})'(\rmF) \cap \mathfrak{o}}
\sum_{\zeta_{1} \in (\bar \nl_{\tlS}^{\tltwo})'(\rmF)}
|\Psi_{S}(x,\xi,\zeta_{1})|
\end{displaymath}
car 
\( (\bar \nl_{\tlS}^{\tlR})'\oplus 
(\bar \nl_{R}^{2} \cap (\bar \nl_{\tlR}^{\tltwo})')  
\subseteq (\bar \nl_{\tlS}^{\tltwo})'\). 

On se ramène alors à borner, pour $P_{1} \sbs S \sbs P_{2}$ fixés:
\begin{displaymath}
\int_{P_{1}(\A)\backslash U(\A)}
\chi^{T}_{P_{1},P_{2}}(x)
\sum_{\xi \in \ml_{\tlS, \tlone}'(\rmF) \cap \mathfrak{o}}
\sum_{\zeta_{1} \in (\bar \nl_{\tlS}^{\tltwo})'(\rmF)}
|\Psi_{S}(x,\xi,\zeta_{1})|dx.
\end{displaymath}
Cette intégrale est identique à (\ref{mainThmConv3}) du théorème 
\ref{thm:MainConv}. Cela conclut la preuve.
\edem
\end{theo}

Au début de la section \ref{qualitatives} nous avons introduit 
les distributions $J_{\ol}^{T}$. Pour des classes $\ol$ 
dans $\calO_{rs}$ cette distribution s'exprime comme suit.
\brop\label{prop:smalljInt}
Pour \(T \in T_{+} + \mathfrak{a}_{0}^{+}\), 
\(f \in \mathcal{S}(\tlul(\A))\) et \(\mathfrak{o} \in \mathcal{O}_{rs}\) on a
\begin{displaymath}
J^{T}_{\mathfrak{o}}(f) = \int_{U(\rmF) \bsl U(\A)}
j_{\mathfrak{o}}^{T}(x)dx.
\end{displaymath}

\bdem 
Dans la preuve on utilisera la notation 
du chapitre 
\ref{sec:convergence} introduite au début de la preuve du 
théorème \ref{thm:MainConv}. Donc, en raisonnant comme au début 
de la preuve de ce théorème-là, on voit que
\[
\int\limits_{\mathclap{U(\rmF) \bsl U(\A)}}j_{f, \mathfrak{o}}^{T}(x)dx = \! \! \! \! \sum_{P_{2} \sps P_{1} \sps P_{0}}
\qquad
\int\limits_{\mathclap{P_{1}(\rmF) \bsl U(\A)}}
\chi_{1,2}^{T}(x)
\sum_{P_{1} \sbs P \sbs P_{2}}
(-1)^{d_{P}}
\sum_{\xi \in \ml_{\tlP}(\rmF) \cap \ol}
\sum_{\eta \in N_{P}(\rmF)}
f_{\eta x}^{P}(\xi)dx.
\]
Fixons $P_{2} \sps P_{1} \sps P_{0}$ et 
décomposons l'intégrale sur 
$P_{1}(\rmF)\bsl U(\A)$ 
en une double intégrale 
sur $x \in M_{1}(\rmF)N_{1}(\A)\bsl U(\A)$ 
et $n_{1} \in N_{1}(\rmF) \bsl N_{1}(\A)$. 
Ensuite on fait passer cette dernière intégrale à l'intérieure 
de la somme sur $P$. On peut le faire 
car la fonction
\[
N_{1}(\A) \ni n_{1} \mapsto 
\sum_{\xi \in \ml_{\tlP}(\rmF) \cap \ol}
\sum_{\eta \in N_{P}(\rmF)}
f_{\eta n_{1} x}^{P}(\xi)
\] 
est $N_{1}(\rmF)$-invariante est continue
donc bornée sur le compact 
$N_{1}(\rmF)\bsl N_{1}(\A)$.
Pour $P$
et $x \in M_{1}(\rmF)N_{1}(\A)\bsl U(\A)$
fixés, on regarde alors
\[
\int_{[N_{1}]}
\sum_{\xi \in \ml_{\tlP}(\rmF) \cap \ol}
\sum_{\eta \in N_{P}(\rmF)}
f_{\eta n_{1} x}^{P}(\xi)dn_{1}.
\]
Comme le volume de $N_{P}(\rmF) \bsl N_{P}(\A)$ vaut $1$ 
et $N_{P} \sbs N_{1}$, cette expression vaut
\begin{equation*}
\begin{split}
\int\limits_{[N_{1}]}
\sum_{\xi \in \ml_{\tlP}(\rmF) \cap \ol}
\int\limits_{[N_{P}]}
\sum_{\eta \in N_{P}(\rmF)}
f_{\eta n n_{1} x}^{P}( \xi)dndn_{1} &= \\
\int\limits_{[N_{1}]}
\sum_{\xi \in \ml_{\tlP}(\rmF) \cap \ol}
\int\limits_{N_{P}(\A)}
f_{n n_{1} x}^{P}( \xi)dndn_{1} & = 
\int\limits_{\mathclap{[N_{1}]}}
k_{\ol,P}(n_{1}x)dn_{1}.
\end{split}
\end{equation*}
La dernière égalité découle du corollaire 
\ref{cor:lemcor2}. On intervertit de nouveau 
la somme qui porte sur $P$ avec l'intégrale 
sur $N_{1}(\rmF) \bsl N_{1}(\A)$ et l'on 
recombine cette dernière avec l'intégrale sur 
$M_{1}(\rmF)N_{1}(\A)\bsl U(\A)$. On retrouve donc
\[
\int_{U(\rmF) \bsl U(\A)}j_{f, \mathfrak{o}}^{T}(x)dx = 
\sum_{P_{2} \sps P_{1} \sps P_{0}}
\int_{P_{1}(\rmF)\bsl U(\A)}
\chi_{1,2}^{T}(x)k_{1,2,\ol}(x)dx.
\]
Chaque intégrale dans la somme ci-dessus converge d'après le théorème 
\ref{thm:MainConv2} ce qui justifie l'intégration et de surcroît, 
la somme elle-même
égale $J_{\ol}^{T}(f)$ 
grâce à l'identité (\ref{chgmtordre}) et 
la définition 
de $J_{\ol}^{T}(f)$ donnée 
au début du chapitre \ref{qualitatives}
 ce qu'il fallait démontrer. 
 \edem
 \erop
 
 Plaçons nous maintenant dans le cadre du paragraphe \ref{par:JTforLevis}. 
 On veut généraliser le théorème \ref{thm:convThmNotMain} et la proposition 
 \ref{prop:smalljInt} au cas $G \times U'$.
Notons $\calO^{G \times U'}_{rs}$ l'ensemble de classes 
contenant un élément $X_{1} + X_{2} \in \gl(\rmF) \times \tlul'(\rmF)$ 
tel que le polynôme caractéristique de $X_{1}$ est séparable 
et tel que $X_{2} \in \tlul'(\rmF)$ appartient à une classe relativement semi-simple régulière
 dans le contexte d'inclusion $U' \hrar \tlU'$. 
 
Pour $f \in \calS((\gl \times \tlul')(\A))$, $\ol \in \calO^{G \times U'}_{rs}$ 
et un sous-groupe parabolique standard $P$ de $G \times U'$ 
soit
 \[
I_{f,P,\ol}(x) = 
\sum_{\xi \in \ml_{\tlP}(\rmF) \cap \ol}
\sum_{\eta \in N_{P}(\rmF)}\int_{V_{P}'(\A)}
f(\Ad((\eta x)^{-1})(\xi + Y_{P}'))dY_{P}'
 \]
 où $x \in P(\rmF) \bsl (G \times U')(\A)^{1}$
 et $V_{P}'$ c'est le plus grand sous-espace isotrope 
 de $V'$ stabilisé par $P$.
 
 Pour $T \in (\all_{P_{0}}^{G \times U'})^{+}$ 
 et $x \in (G \times U')(\rmF) \bsl (G \times U')(\A)^{1}$ on pose aussi
 \[
j_{f,\ol}^{T}(x) = 
\sum_{P \sps P_{0}}(-1)^{d_{P}^{G \times U'}}
\sum_{\mathclap{\delta \in P(\rmF) \bsl (G \times U')(\rmF)}}
\
\htau_{P}^{G \times U'}(H_{P}(\delta x) - T_{P})
I_{P,\ol}(\delta x).
 \] 
 La preuve du théorème \ref{thm:convThmNotMain} s'étend sans problème 
 dans ce cas donnant
 \[
\int_{(G \times U')(\rmF) \bsl (G \times U')(\A)^{1}}
|j_{\ol}^{T}(x)|dx < \infty 
\]
pour $T$ suffisamment régulier. De même la preuve 
de la proposition \ref{prop:smalljInt} 
s'étend aussi bien et l'on obtient, 
avec la notation du paragraphe \ref{par:JTforLevis}, 
pour $\ol\in \calO^{G \times U'}_{rs}$
\[
\int_{(G \times U')(\rmF) \bsl (G \times U')(\A)^{1}}
j_{f,\ol}^{T}(x)dx = 
J_{\ol}^{G \times U',T}(f).
\]

Soient maintenant $Q$ un sous-groupe parabolique standard de $U$, 
$\ol \in \calO_{rs}$ et 
$\ol_{Q,1}, \ldots, \ol_{Q,m} \in \calO^{M_{Q}}$ 
comme dans l'équation (\ref{eq:mltlQcapol}). Alors $\ol_{Q,i} \in \calO^{M_{Q}}_{rs}$
pour $i = 1, \ldots, m$. 
Si $P \sbs Q$, 
alors $V_{P}^{Q} := V_{P} \cap Z_{Q}$ c'est le plus grand sous-espace isotrope 
de $Z_{Q}$ stabilisé par $P \cap M_{Q}$.
On a dans ce cas l'analogue de 
l'égalité (\ref{eq:JMQisThis}) suivant pour tout $f \in \calS(\tlul(\A))$:
\begin{multline}\label{eq:JMQisThisRs}
J^{M_{Q},T}_{\mathfrak{o}}(f_{Q}) \! = \! \!
\int_{M_{Q}(\rmF) \backslash M_{Q}(\A)^{1}}
\sum_{i=1}^{m}j_{f_{Q},\ol_{Q,i}}^{T}(m)dm \! = \!
\int_{M_{Q}(\rmF) \backslash M_{Q}(\A)^{1}}
\sum_{P \sbs Q}(-1)^{d_{P}^{Q}} \\
\sum_{\mathrlap{\eta \in (P \cap M_{Q})(\rmF) \bsl M_{Q}(\rmF)}} 
\htau_{P}^{Q}(H_{P}(\eta m)-T) \!
\dsl
\sum_{\xi \in \ml_{\tlP}(\rmF) \cap \ol} \!
\sum_{\delta \in N_{P}^{Q}(\rmF)} \!
\int_{V_{P}^{Q}(\A)}
f_{Q}(\Ad((\delta \eta m)^{-1})(\xi + Y_{P}^{Q}))dY_{P}^{Q} \!
\rb dm.
\end{multline}

\subsection{Expression intégrale de $J_{\ol}$}\label{par:IntReprDeJ}

Dans ce paragraphe, on démontre 
la proposition \ref{prop:newjol0}.

\brop\label{prop:newjol} On a
\[
\int_{U(\rmF) \bsl U(\A)}|j_{\ol}(x)|dx < \infty
\quad \text{et} \quad \int_{U(\rmF) \bsl U(\A)}j_{\ol}(x)dx = J_{\ol}(f).
\]

\bdem 
Pour tout sous-groupe parabolique standard $Q$ de 
$U$ soit $\brtau_{Q}$ la fonction caractéristique de 
$H \in \all_{Q}$ tels que $\al(H) \le 0$ 
pour tout $\al \in \Delta_{Q}$. 
Il résulte du lemme combinatoire de Langlands 
(Proposition 1.7.2 de \cite{labWal}) que pour tout 
sous-groupe parabolique standard $P$ de $U$ on a
\[
\sum_{Q \sps P}\htau_{P}^{Q}\brtau_{Q} = 1.
\]
En utilisant cette identité, on a 
pour tout $T \in T_{+} + \all_{0}^{+}$
\begin{equation*}
\begin{split}
& j_{\ol}(x) =
\sum_{P}(-1)^{d_{P}}
\sum_{\delta \in P(\rmF) \bsl U(\rmF)}
I_{P,\ol}(\delta x) = \\
& \sum_{P}(-1)^{d_{P}} \ \ 
\sum_{\mathclap{\delta \in P(\rmF) \bsl U(\rmF)}} \ \
I_{P,\ol}(\delta x)
(\sum_{Q \sps P}\htau_{P}^{Q}(H_{P}(\delta x) - T_{P})
\brtau_{Q}(H_{Q}(\delta x) - T_{Q})) = \\
& \sum_{Q}(-1)^{d_{Q}} \ \ 
\sum_{\mathclap{\delta \in Q(\rmF) \bsl U(\rmF)}} \ \ 
\brtau_{Q}(H_{Q}(\delta x) - T_{Q})
\sum_{P \sbs Q} \ 
(-1)^{d_{Q}^{P}}
\sum_{\mathclap{\eta \in (M_{Q} \cap P)(\rmF) \bsl M_{Q}(\rmF)}}
I_{P,\ol}(\eta \delta x)\htau_{P}^{Q}(H_{P}(\eta \delta x) - T_{P}),
\end{split}
\end{equation*}
où les sommes sont absolument convergentes.
Il suffit de montrer que pour tout $Q$ l'intégrale
\begin{equation}\label{eq:suffitPourQ}
\int_{Q(\rmF) \bsl U(\A)}
\brtau_{Q}(H_{Q}(x) - T_{Q})
\sum_{P \sbs Q}
(-1)^{d_{Q}^{P}} \ \ 
\sum_{\mathclap{\eta \in (M_{Q} \cap P)(\rmF) \bsl M_{Q}(\rmF)}} \
I_{P,\ol}(\eta x)\htau_{P}^{Q}(H_{P}(\eta x) - T_{P})dx
\end{equation}
converge absolument. 
L'analyse va être analogue à celle de la preuve du théorème 
\ref{mainQualitThm}.

Posons
$x = namk$ o\`u $n \in N_{Q}(\rmF)\backslash N_{Q}(\A)$, 
$m \in M_{Q}(\rmF)\backslash M_{Q}(\A)^{1}$, 
$a \in A_{Q}^{\infty}$ et $k \in K$. Donc 
$dx = e^{-2\rho_{Q}(H_{Q}(a))}dndadmdk$. 

Fixons $P \sbs Q$.
Pour $n$, $m$, $k$ et $a$ comme ci-dessus et $\eta \in M_{Q}(\rmF)$ on a
$\brtau_{Q}(H_{Q}(\eta namk) - T_{Q}) = \brtau_{Q}(H_{Q}(a) - T_{Q})$
et, en faisant les changements de variable 
$a^{-1}na \mapsto n$ et $a^{-1}Y_{P}a \mapsto Y_{P}$, 
\begin{equation}\label{eq:passAuLevi}
\begin{split}
\int_{K}\int_{[N_{Q}]}
I_{P,\ol}(\eta n amk)\htau_{P}^{Q}(H_{P}(\eta namk) - T_{P})dndk  &= \\ 
\int_{K}\int_{[N_{Q}]}
I_{P,\ol}(\eta n amk)\htau_{P}^{Q}(H_{P}(\eta m) - T_{P})dndk & = \\
\htau_{P}^{Q}(H_{P}(\eta m) - T_{P}) e^{2\rho_{\tlQ}(H_{Q}(a))}
\int\limits_{K} \sum_{\xi \in \ml_{\tlP}(\rmF) \cap \ol} 
\sum_{\delta \in N_{P}^{Q}(\rmF)} &
\int\limits_{N_{Q}(\A)}
\int\limits_{V_{P}(\A)}
f_{n \delta \eta mk}(\xi + Y_{P})dY_{P}dndk.
\end{split}
\end{equation}

Fixons $\xi \in \ml_{\tlP}(\rmF) \cap \ol$.
Pour une fonction $\upphi \in \calS(\tlul(\A))$ 
regardons
\begin{equation}\label{eq:intAuxRss}
\int_{N_{Q}(\A)}
\int_{V_{P}(\A)}
\upphi(n^{-1} (\xi + Y_{P})
 n)dY_{P}dn.
\end{equation}
On a 
$V_{P} = V_{Q} \oplus V_{P}^{Q}$, où $V_{P}^{Q} = V_{P} \cap Z_{Q}$.
Pour $Y_{P} \in V_{P}$ soient $Y_{Q} \in V_{Q}$ 
et $Y_{P}^{Q} \in V_{P}^{Q}$ tels que $Y_{P} = Y_{Q} + Y_{P}^{Q}$. 
Donc, puisque $Y_{P}^{Q} \in Z_{Q}$, si 
$n \in N_{Q}$ on a $Y_{P}^{Q} - nY_{P}^{Q}n^{-1} \in V_{Q}$.
L'intégrale (\ref{eq:intAuxRss}) ci-dessus égale donc
\begin{multline*}
\int_{V_{P}^{Q}(\A)}
\int_{N_{Q}(\A)}
\int_{V_{Q}(\A)}
\upphi(n^{-1} (\xi + Y_{Q} + Y_{P}^{Q} - nY_{P}^{Q}n^{-1})n
+
Y_{P}^{Q})dY_{Q}dndY_{P}^{Q} = \\
\int_{V_{P}^{Q}(\A)}
\int_{N_{Q}(\A)}
\int_{V_{Q}(\A)}
\upphi(n^{-1} (\xi + Y_{Q})n
+
Y_{P}^{Q})dY_{Q}dndY_{P}^{Q}.
\end{multline*}
En utilisant le corollaire \ref{cor:lemcor2} et en faisant un changement 
de variable ceci égale
\[
\int_{V_{P}^{Q}(\A)}
\int_{\nl_{\tlQ}(\A)}
\upphi(\xi + U_{Q}+Y_{P}^{Q})dU_{Q} dY_{P}^{Q}.
\]
On voit alors que (\ref{eq:passAuLevi}) devient:
\[
e^{2\rho_{\tlQ}(H_{Q}(a))} \htau_{P}^{Q}(H_{P}(\eta m) - T_{P}) \! \!
\sum_{\xi \in \ml_{\tlP}(\rmF) \cap \ol}
\sum_{\delta \in N_{P}^{Q}(\rmF)}
\int\limits_{V_{P}^{Q}(\A)}
f_{Q}(\Ad((\delta \eta m)^{-1})(\xi + Y_{P}^{Q}))dY_{P}^{Q}
\]
où $f_{Q} \in \calS(\ml_{\tlQ}(\A))$ est définie 
par (\ref{eq:fQdef}) dans le 
paragraphe \ref{par:JTforLevis}.

En utilisant alors l'égalité (\ref{eq:JMQisThisRs}), 
on s'aperçoit que l'intégrale (\ref{eq:suffitPourQ})
égale
$
J_{\ol}^{M_{Q},T}(f_{Q})
$
fois
\begin{equation}\label{eq:notQuitePQ}
\int_{\all_{Q}}e^{\underline{\rho}_{Q}(H)}\brtau_{Q}(H-T_{Q})dH.
\end{equation}
La convergence de l'intégrale dans le théorème va alors découler
de la convergence de l'intégrale (\ref{eq:notQuitePQ}) ci-dessus. 
Or, cette intégrale converge en vertu du lemme 
\ref{lem:positiveLinear} et donne 
précisément 
$\hat \theta_{Q}(\underline{\rho}_{Q})^{-1} e^{\underline{\rho}_{Q}(T_{Q})}$
où $\hat \theta_{Q} = \hat \theta_{Q}^{U} $ est définie 
par (\ref{eq:thetaHatDef}) dans le paragraphe \ref{par:fonsPolExp}.

On vient d'obtenir la convergence ainsi que 
pour tout $T \in \all_{0}^{+}$ 
suffisamment régulier
\[
\int_{[U]}j_{\ol}(x)dx = 
\sum_{Q}(-1)^{d_{Q}}
\hat \theta_{Q}(\underline{\rho}_{Q})^{-1}e^{\underline{\rho}_{Q}(T_{Q})}J_{\ol}^{M_{Q},T}(f_{Q}).
\]
D'après le théorème \ref{mainQualitThm} la somme ci-dessus égale 
$J_{\ol}(f)$, ce qu'il fallait démontrer. 
\edem
\erop

\subsection{Résultats de convergence}\label{par:resDeConv}

Pour toute place $v$ de $\rmF$ on note $\rmF_{v}$ son
complété à $v$. Soit $S_{\infty}$ l'ensemble des places archimédiennes
de $\rmF$.
Par définition d'un sous-groupe compact 
maximal admissible de $U(\A)$ par rapport à $M_{0}$
(voir paragraphe 1 de \cite{arthur2}), on a 
$K = \prod_{v}K_{v}$ où $K_{v}$ est un 
sous-groupe compact 
maximal admissible de $U(\rmF_{v})$ par rapport à $M_{0}$. 
Pour toute place $v$ de $\rmF$, 
on choisit une mesure de Haar $dx_{v}$ 
sur $U(\rmF_{v})$ de façon que la mesure de Haar sur $U(\A)$ 
soit leur produit. De même, 
pour tout $i \in I$ on choisit une mesure de Haar sur $T_{i}(\A)$ 
et des mesures de Haar sur $T_{i}(\rmF_{v})$ pour toute place 
$v$ de $\rmF$ de façon que la mesure de Haar sur $T_{i}(\A)$ 
soit leur produit. Pour tout $I' \sbs I$ et toute place $v$ de $\rmF$
on met la mesure produit sur $T_{I'}(\rmF_{v}) = \prod_{i \in I'}T_{i}(\rmF_{v})$ 
et sur $T_{I'}(\A) = \prod_{i \in I'}T_{i}(\A)$. On fixe
aussi sur $K_{v}$ la mesure de Haar de masse totale $1$.

On suit \cite{labWal}, paragraphe 3.2. 
On fixe une $\rmF$-base de $V$ et à partir d'elle, 
pour toute place $v$ de $\rmF$,
on définit des normes notées $| \cdot |_{v}$ sur 
$V(\rmF_{v}) := 
V \otimes_{\rmF} \rmF_{v} \cong \rmE_{I} \otimes_{\rmF} \rmF_{v}$ 
et $\End(V(\rmF_{v}))$. Grâce à l'inclusion 
$\ul(\rmF_{v}) \hrar \End_{\rmF_{v}}(V(\rmF_{v}))$ 
on obtient une norme $|\cdot |_{v}$ sur $\ul(\rmF_{v})$ par restriction. 
On fixe une norme, notée aussi  $| \cdot|_{v}$, sur $U(\rmF_{v})$ 
grâce à l'inclusion $U(\rmF_{v}) \hrar \End_{\rmF_{v}}(V(\rmF_{v})) \oplus \End_{\rmF_{v}}(V(\rmF_{v}))$ 
donnée par $U(\rmF_{v}) \ni g \mapsto (g, \ps{t}g^{-1})$. 
On a donc 
$|x_{v}|_{v} = |x_{v}^{-1}|_{v}$ 
pour tout $x_{v} \in U(\rmF_{v})$. On peut et on va supposer en plus 
que $|k_{v}|_{v} = 1$ pour tout $k_{v} \in K_{v}$.
On en 
déduit des hauteurs sur $V(\A) := V \otimes_{\rmF} \A$, $U(\A)$ et $\ul(\A)$ 
par $\| \cdot \| := \prod_{v}|\cdot|_{v}$.
Il existe alors des constantes $c_{0},c_{1}, c_{2} >0$ 
telles que pour tous $x_{1}, x_{2} \in U(\A)$ et 
tout $v \in V(\A)$ on a
\[
\|x_{1}^{-1}\| = \|x_{1}\|, \quad 
\|x_{1}\| \ge c_{0}, \quad 
\|x_{1}x_{2}\| \le c_{1}\|x_{1}\|\|x_{2}\|, \quad 
\|x_{1}v\| \le c_{2}\|x_{1}\|\|v\|.
\]

On note aussi qu'il existe des constantes $c_{3}, r_{3} > 0$
telles que
\begin{equation}\label{eq:expInequality}
|e^{\la(H_{0}(x))}| \le c_{3} e^{\|\Rel(\la)\|}\|x\|^{r_{3}}, \quad 
\forall \ x \in U(\A), \ \la \in \all_{I_{0},\C}^{*}.
\end{equation}

Notons finalement que pour simplifier l'exposition ci-dessous et éviter l'expression 
$0^{-1}$, on pose $\|0\| = 1$ où $0 \in V(\A)$.

\blem\label{lem:tAndVectors}
 Soit $\calI \sbse I_{0}$. Il existe des constantes $r_{0}, c_{0} > 0$ telles que
pour tout $t \in T_{\acalI}(\A)$ et tout $x \in U(\A)$ on a
\[
\ind_{\calI}(H_{0}(tx))\|t\| \le c_{0} \|t^{-1}1_{\calI}\|^{r_{0}}\|x\|^{r_{0}}.
\]
\bdem
C'est une conséquence de définitions des hauteurs ci-dessus ainsi que de 
la manière par laquelle $T_{\acalI}(\A)$ agit sur $\rmF_{\calI} \otimes_{\rmF} \A$,
qu'il existe des constantes $c', r' >0$ telles que 
\begin{equation}\label{eq:tAndVectors}
\|t\| \le c' \max(\|t^{-1}1_{\calI}\|^{r'}, \|t1_{\calI}\|^{r'}).
\end{equation}
Si l'on suppose $\ind_{\calI}(H_{0}(tx)) = 1$, il résulte de la définition de 
$\ind_{\calI}$ donnée au paragraphe \ref{par:defsOrbs} qu'il existe 
des constantes $c'', r'' >0$ telles que 
$\|t^{-1}1_{\calI}\| \ge c'' \|x\|^{-r''}$ et 
$\|t1_{\calI}\| \le c'' \|x\|^{r''}$. En utilisant ceci 
dans (\ref{eq:tAndVectors}) ci-dessus et en utilisant le fait que 
$\|x\| \ge c_{0}$ on trouve le résultat voulu.
\edem
\elem

\blem\label{lem:txiNothingConv}
Il existe un $r_{0} > 0$ tel que pour tout $r \ge r_{0}$ 
on a
\[
\int_{T_{I_{2} \smin I_{0}}(\rmF) 
\bsl T_{I \smin I_{0}}(\A)}
\|t^{-1}\xi_{\varnothing}\|^{-r}dt < \infty.
\]
\bdem
En utilisant la définition de l'ensemble $I_{2}$ donnée avant la proposition 
\ref{prop:orbitsRRSS} ainsi que la description de $\xi_{\varnothing}$ donnée dans cette proposition-là, 
on regarde les projections de $t^{-1}\xi_{\varnothing}$ à 
$\rmE_{i} \otimes_{\rmF} \A$ pour tout $i \in I \smin I_{0}$ 
et on s'aperçoit qu'il existe 
une constante $c >0 $ telle que
l'intégrale dans le lemme est 
majorée par le produit sur $i \in I \smin I_{2}$ de
\[
\int_{0}^{\infty} \min(t^{cr}, t^{-cr}) dt 
\]
fois le volume de $T_{I_{2} \smin I_{0}}(\rmF) \bsl T_{I_{2} \smin I_{0}}(\A)$. 
En prenant $r$ suffisamment grand on obtient donc la convergence. 
\edem
\elem

\blem\label{lem:t1calIConv}
Soit $\calI \sbse I_{0}$. 
Il existe un $r_{0} >0$ tel que 
pour tout $r \ge r_{0}$ il existe un $r' > r$ et une constante $c >0$ 
telles que
 pour tout $x \in U(\A)$ on a
\[
\int_{T_{\acalI}(\A)}
\ind_{\calI}(H_{0}(tx))
\|x^{-1}t^{-1}1_{\calI}\|^{-r} dt 
\le c \inf_{\mathmakebox[0.5cm]{t \in T_{I}(\A)}}\|tx\|^{r'}.
\]

\bdem 
On a
\[
\|x^{-1}t^{-1}1_{\calI}\| \ge c_{2}^{-1} \|x\|^{-1} \|t^{-1}1_{\calI}\|.
\]
On voit donc qu'il existe des constantes $c >0 $, $c' >0 $ telles que
l'intégrale du lemme est 
majorée par le produit sur $i \in \calI$ de
\[
\|x\|^{r}\int_{c'\|x\|^{-r_{3}}}^{\infty}t^{-cr}dt
\]
où $r_{3}$ est comme dans l'inégalité 
(\ref{eq:expInequality}), ce qui converge pour $r$ suffisamment grand. 
On obtient le résultat voulu en remarquant que l'intégrale du lemme 
est $T_{I}(\A)$-invariante à gauche en tant qu'une fonction de $x \in U(\A)$.
\edem
\elem

\blem\label{lem:ChaudnMajorAd} Soit $v$ une place de $\rmF$. Il existe des constantes
$c_{B}, N_{B} >0$ telles que pour 
tout $x_{v} \in U(\rmF_{v})$
on a
\[
\inf_{\mathmakebox[0.7cm]{t_{v} \in T_{I}(\rmF_{v})}}|t_{v}x_{v}|_{v} 
\le c_{B}(1 + |\Ad((x_{v})^{-1})B|_{v})^{N_{B}}.
\]

\bdem 
Supposons d'abord $T_{I}$ déployé sur $\rmF_{v}$. 
Dans ce cas $T_{I}$ est une partie de Levi d'un $\rmF_{v}$-sous-groupe 
de Borel de $U$ de la forme $T_{I}N_{v}$ où $N_{v}$ c'est sa partie 
unipotente. De plus $B \in \Lie(T_{I})(\rmF_{v})$.
En utilisant la décomposition d'Iwasawa on écrit 
$x_{v} = t_{v}n_{v}k_{v}$ où $t_{v} \in T_{I}(\rmF_{v})$, 
$n_{v} \in N_{v}(\rmF_{v})$ et $k_{v} \in K_{v}$.
Soient $U_{v}, U_{v}' \in \Lie(N_{v})(\rmF_{v})$ 
tels que $n_{v} = \exp(U_{v})$ et $\Ad(n_{v}^{-1})B = B+ U_{v}'$.
D'après le lemme 2.1 de \cite{chaud2}, il existe un polynôme 
$Q$ sur $\Lie(N_{v})$ 
qui ne dépend que de $B$, tel que 
$U_{v} = Q(U_{v}')$. Il est clair qu'il existe des 
constantes $c_{B}'$, $c_{B}''$, $N_{B}'$, qui ne dépendent 
que de $B$, telles que $(1 + |U_{v}|_{v})$ est plus petit que
\[
c_{B}'(1  + |B + U_{v}'|_{v})^{N_{B}'} = 
c_{B}'(1 +|\Ad((t_{v}n_{v})^{-1})B|_{v})^{N_{B}'} \le 
c_{B}''(1 +|\Ad((x_{v}^{-1})B|_{v})^{N_{B}'}.
\]
D'autre part, il existe des constants $c_{0}$, $N_{0}$ telles que:
\[
|n_{v}|_{v} \le c_{0}(1 + |U_{v}|_{v})^{N_{0}}
\] 
d'où
\[
\inf_{\mathmakebox[0.7cm]{t_{v} \in T_{I}(\rmF_{v})}}|t_{v}x_{v}|_{v} \le 
|n_{v} k_{v}|_{v} \le 
\sup_{k_{v} \in K_{v}} |k_{v}|_{v} |n_{v}|_{v} \le c_{B}(1 + |\Ad((x_{v})^{-1})B|_{v})^{N_{B}}.
\]

Dans le cas général, soit $\rmF'_{v}$ une extension finie de 
$\rmF_{v}$ qui déploie $T_{I}$. On prolonge la hauteur 
$| \cdot |_{v}$ à $U(\rmF_{v}')$. Grâce à (4.6) de \cite{arthur6}, 
il existe des constantes $c$, $N$ telles que:
\[
\inf_{\mathmakebox[0.8cm]{t_{v} \in T_{I}(\rmF_{v})}}|t_{v}x_{v}|_{v} \le 
c \inf_{\mathmakebox[0.8cm]{t_{v}' \in T_{I}(\rmF'_{v})}}|t_{v}'x_{v}|_{v}^{N} 
\]
et le résultat suit du cas déployé. 
\edem
\elem
 
 \blem\label{lem:mainLemCvgRegss} 
 Soient $\calJ_{2} \sbs \calJ_{3} \sbse I_{0}$
 et $\calJ_{1} \sbse I_{0}$ tels que $\acalJj \cap \acalJt = \varnothing$.
 Alors, 
 pour tout $f \in \calS(\tlul(\A))$
et 
 tout compact $C \sbs \all_{|\calJ_{1}\cup \calJ_{2}|}^{*}$
 l'intégrale 
\[
\int\limits_{\mathclap{
T_{I_{2} \smin I_{0}}(\rmF)T_{I_{0} 
\smin |\calJ_{12}|}(\A) \bsl U(\A)}}  \ 
\ind_{\calJ_{1} \cup \calJ_{2}^{\sharp}}(H_{0}(x))
e^{(\la + \upla_{\calJ_{3 \smin 2}^{\sharp}})(H_{0}(x))}
\hf_{x}^{\calJ_{3}}(X_{\calJ_{1} \cup \calJ_{2}^{\sharp}})
dx
\]
converge absolument et est majorée indépendamment de $\la \in \all_{|\calJ_{1}\cup \calJ_{2}|,\C}^{*}$ 
tel que $\Rel(\la) \in C$.

 \brems Soit $t \in T_{I_{2} \smin I_{0}}(\rmF)T_{I_{0}\smin |\calJ_{12}|}(\A)$. 
 Pour tout $x \in U(\A)$ on a alors
 \[
 \ind_{\calJ_{1} \cup \calJ_{2}^{\sharp}}(H_{0}(tx)) = \ind_{\calJ_{1} \cup \calJ_{2}^{\sharp}}(H_{0}(x)), \ 
 e^{(\la + \upla_{\calJ_{3 \smin 2}^{\sharp}})(H_{0}(tx))} \! = \!
 e^{(\la + \upla_{\calJ_{3 \smin 2}^{\sharp}})(H_{0}(x))}e^{\upla_{\calJ_{3 \smin 2}^{\sharp}}(H_{0}(t))}
 \]
 et 
 \[
 \hf_{tx}^{\calJ_{3}}(X_{\calJ_{1} \cup \calJ_{2}^{\sharp}}) = 
 e^{-\upla_{\calJ_{3 \smin 2}^{\sharp}}(H_{0}(t))}\hf_{x}^{\calJ_{3}}(X_{\calJ_{1} \cup \calJ_{2}^{\sharp}})
 \]
 d'où $T_{I_{2} \smin I_{0}}(\rmF)T_{I_{0}\smin |\calJ_{12}|}(\A)$-invariance de l'intégrale considérée. 
 \erems 
 
\bdem 
Pour $\calJ \sbse I_{0}$ et $I' \sbs I$ soient $\A_{\calJ} = \rmF_{\calJ} \otimes_{\rmF} \A$ 
et $\A_{I'} = \prod_{i \in I'} \rmE_{i} \otimes_{\rmF} \A$.
Soit $h_{x}$ la fonction sur $\ul(\A) \times \A_{I \smin |\calJ_{3}|} \times \A_{\calJ_{3}^{\sharp}}$ 
définie pour tout $x \in U(\A)$ de façon suivante:
\[
h_{x}(B', u_{I \smin |\calJ_{3}|}, u_{\calJ_{3}^{\sharp}}) = 
\int\limits_{\mathclap{\A_{\calJ_{3}}}}
f_{x} \dsl \! \!
\matx{B'}{u_{I \smin |\calJ_{3}|} + u_{\calJ_{3}}}{(u_{I \smin |\calJ_{3}|} + u_{\calJ_{3}})^{\sharp}}{d_{\ol}}
\! \!\rb  \! \psi ( \langle u_{\calJ_{3}}, u_{\calJ_{3}^{\sharp}} \rangle ) du_{\calJ_{3}}.
\]
Alors $h_{x} \in \calS(\ul(\A) \times \A_{I \smin |\calJ_{3}|} \times \A_{\calJ_{3}^{\sharp}})$.
Soit $P \in \calP(M_{I_{1}})$ un sous-groupe stabilisant 
le drapeau $\rmF_{\calJ_{3}} \sbs \rmF_{\calJ_{13}}$.
Soient 
$m \in M_{I_{1}}(\A)$, $n \in N_{P}(\A)$ et $k \in K$. En faisant le changement de variable 
$(mn)^{-1}u_{\calJ_{3}} \mapsto u_{\calJ_{3}}$ on obtient 
que $|\hf_{mnk}^{\calJ_{3}}(X_{\calJ_{1} \cup \calJ_{2}^{\sharp}})|$ vaut
\[
|e^{\upla_{\calJ_{3}}(H_{P}(m))}
h_{k}(\Ad((mn)^{-1})B, (1 - 1_{\calJ_{3}})((mn)^{-1}(1_{\calJ_{1}} + \xi_{\varnothing})), 1_{\calJ_{3}^{\sharp}}((mn)^{-1}1_{\calJ_{2}^{\sharp}}) ))|
\]
où $1$ c'est l'unité dans $\rmE_{I}$.

On voit donc que pour tout $r > 0$ il existe une fonction $\Upphi \in \calS(\ul(\A))$ telle que
\begin{multline*}
|\hf_{mnk}^{\calJ_{3}}(X_{\calJ_{1} \cup \calJ_{2}^{\sharp}})| \le \\
e^{\upla_{\calJ_{3}}(H_{P}(m))} \Upphi(\Ad((mn)^{-1})B)\|(1 - 1_{\calJ_{3}})((mn)^{-1}(1_{\calJ_{1}} + \xi_{\varnothing}))\|^{-r}
\| 1_{\calJ_{3}^{\sharp}}((mn)^{-1}1_{\calJ_{2}^{\sharp}})\|^{-r}.
\end{multline*}
Il existe des 
constantes $c, c_{0}, c_{1}, c_{2}, c_{3} >0$ qui ne dépendent pas de $r$ telles que
\begin{multline*}
\|(1 - 1_{\calJ_{3}})((mn)^{-1}(1_{\calJ_{1}} + \xi_{\varnothing}))\|^{-r}
\| 1_{\calJ_{3}^{\sharp}}((mn)^{-1}1_{\calJ_{2}^{\sharp}})\|^{-r} \le \\
c
\|n\|^{c_{3}r}\|m^{-1}\xi_{\varnothing}\|^{-c_{0}r}
\|m^{-1}1_{\calJ_{1}}\|^{-c_{1}r}
\|m^{-1}1_{\calJ_{2}^{\sharp}}\|^{-c_{2}r}.
\end{multline*}

On va montrer alors la convergence de l'intégrale suivante:
\begin{multline}\label{eq:convIwasInev}
\int_{T_{I}(\A) \bsl M_{I_{1}}(\A)}
\int_{N_{P}(\A)} \int_{K}\Upphi(\Ad((mn)^{-1})B)
\|n\|^{c_{3}r}
\int\limits_{\mathclap{T_{I_{2} \smin I_{0}}(\rmF) 
\bsl T_{|\calJ_{12}| \cup (I \smin I_{0})}(\A)}}
\ind_{\calJ_{1}}(H_{0}(tmnk))
\ind_{\calJ_{2}^{\sharp}}(H_{0}(tmnk)) \\
\|(tm)^{-1}\xi_{\varnothing}\|^{-c_{0}r}
\|(tm)^{-1}1_{\calJ_{1}}\|^{-c_{1}r}
\|(tm)^{-1}1_{\calJ_{2}^{\sharp}}\|^{-c_{2}r}
|e^{(\la + \upla_{\calJ_{3 \smin 2}^{\sharp}})(H_{0}(tmnk))+ \upla_{\calJ_{3}}(H_{P}(tm))}| \\
dtdk dndm.
\end{multline}

Majorons l'exponentielle qui apparaît sous l'intégrale. 
Notons qu'on a $\upla_{\calJ_{3}} = - \upla_{\calJ_{3 \smin 2}^{\sharp}} + \upla_{\calJ_{2}}$. 
On va majorer alors d'abord 
$|e^{\la(H_{0}(tmn)) + \upla_{\calJ_{2}}(H_{P}(tm))}|$ et 
puis l'expression $e^{\upla_{\calJ_{3 \smin 2}^{\sharp}}
(H_{0}(tmnk) - H_{P}(tm))}$.
Fixons donc $m$, $n$ et $k$ comme ci-dessus. 
Soit $t \in T_{|\calJ_{12}| \cup (I \smin I_{0})}(\A)$ tel que 
$\ind_{\calJ_{1}}(H_{0}(tmnk))
\ind_{\calJ_{2}^{\sharp}}(H_{0}(tmnk)) =1$. 
Notons $t_{1}$ (resp. $t_{2}$) sa projection à $T_{|\calJ_{1}|}(\A)$ 
(resp. $T_{|\calJ_{2}|}(\A)$). En utilisant la propriété (\ref{eq:expInequality}) 
ainsi que le lemme \ref{lem:tAndVectors}
on voit qu'il existe des constantes $c', c'', r_{1}', r_{1} >0$ telles que pour tout 
$\Rel(\la) \in C$ on a
\begin{equation*}
\begin{split}
|e^{\la(H_{0}(tmn)) + \upla_{\calJ_{2}}(H_{P}(tm))}| & = 
|e^{\la(H_{0}(t_{1}t_{2}mn)) + \upla_{\calJ_{2}}(H_{P}(t_{1}t_{2}m))}| \\
& \le c'\|n\|^{r_{1}'} \|m\|^{r_{1}'} \|t_{1}\|^{r_{1}'}\|t_{2}\|^{r_{1}'}  \\
\le 
\|t_{1}^{-1}1_{\calJ_{1}}\|^{r_{1}}\|t_{2}^{-1}1_{\calJ_{2}^{\sharp}}\|^{r_{1}}  & = 
c''\|n\|^{r_{1}} \|m\|^{r_{1}}\|t^{-1}1_{\calJ_{1}}\|^{r_{1}}\|t^{-1}1_{\calJ_{2}^{\sharp}}\|^{r_{1}}.  
\end{split}
\end{equation*}
En utilisant
(\ref{eq:expInequality}), vu que
la fonction $M_{I_{1}}(\A) \ni m \mapsto e^{ \quad \mathclap{\upla_{\calJ_{3 \smin 2}^{\sharp}}} \ 
(H_{0}(tmnk) - H_{P}(tm))}$ est $T_{I}(\A)$-invariante à gauche, 
on trouve $c'', r_{2} >0$ telles que
\[
e^{\upla_{\calJ_{3 \smin 2}^{\sharp}}
(H_{0}(tmnk) - H_{P}(tm))} \le c'' \|n\|^{r_{2}}\inf_{\mathmakebox[0.5cm]{t \in T_{I}(\A)}}\|t m\|^{r_{2}}.
\]

Remarquons qu'il existe $c_{4} >0$ tel que pour $\calI = \calJ_{1}, \calJ_{2}^{\sharp}$
on a
\[
\|t^{-1}1_{\calI} \| \le c_{4}\|m\|\|n\|\|(tmnk)^{-1}1_{\calI} \|,  \quad 
\|(tm)^{-1}1_{\calI} \|^{-1} \le c_{4}\|n\| \|(tmnk)^{-1}1_{\calI} \|^{-1}.
\]
En plus $\|(tm)^{-1}\xi_{\varnothing}\| \le c_{2}\|m\| \|t^{-1}\xi_{\varnothing}\|$. 
Donc, en prenant $r$ suffisamment grand pour qu'on puisse appliquer les lemmes \ref{lem:txiNothingConv} 
et \ref{lem:t1calIConv}, on voit que l'intégrale 
sur le tore dans (\ref{eq:convIwasInev}) est majorée par 
$\|n\|^{r_{3}} \|m\|^{r_{3}} \inf_{t \in T_{I}(\A)}\|tm\|^{r_{4}}$ 
pour certaines constantes positives $r_{3}, r_{4}, r_{5}$. 
Puisqu'elle est $T_{I}(\A)$-invariante à gauche, 
en tant qu'une fonction de $m$, on voit qu'elle est simplement majorée 
par 
$\|n\|^{r_{3}}\inf_{t \in T_{I}(\A)}\|tm\|^{r_{3} + r_{4}}$. 

On vient de se ramener à majorer l'intégrale suivante:
\begin{equation}\label{eq:almostVaradarajan}
\int_{T_{I}(\A) \bsl M_{I_{1}}(\A)}
\int_{N_{P}(\A)} \Upphi(\Ad((mn)^{-1})B)
\|n\|^{r'}
\inf_{\mathmakebox[0.5cm]{t \in T_{I}(\A)}}\|tm\|^{r''}dndm
\end{equation}
pour certaines constantes positives $r', r'' >0$.

Quitte à majorer $\Upphi$, on peut supposer que 
$\Upphi = \otimes_{v} \Upphi_{v}$ 
où $\Upphi_{v}$ 
est une fonction dans la classe de Schwartz 
sur $\ul(\rmF_{v})$ si $v \in S_{\infty}$ 
et $\Upphi_{v}$ est une fonction caractéristique 
d'un compact ouvert $\kl_{v}$ dans
$\ul(\rmF_{v})$ si $v$ est une place finie.
On fixe un ensemble fini
 de places $S$ qui contient l'ensemble
 $S_{\infty}$ et qui vérifie 
 en plus pour toute place $v \nin S$
 \begin{enumerate}
 \item $ \forall \ k_{v} \in K_{v}, \ \Ad(k_{v})B \in \kl_{v}$,
 \item $\forall x_{v} \in U(F_{v})$ tel que 
 $\Ad(x_{v}^{-1})B \in \kl_{v}$ on a $x_{v} \in T_{I}(\rmF_{v})K_{v}$.
 \end{enumerate}
C'est possible, comme il est expliqué dans le paragraphe 
5.2 de \cite{chaud2}. 

L'intégrale (\ref{eq:almostVaradarajan}) est alors majorée par le produit sur
 $v \in S$ de
\[
\int_{T_{I}(\rmF_{v}) \bsl M_{I_{1}}(\rmF_{v})}
\int_{N_{P}(\rmF_{v})} \Upphi_{v}(\Ad((m_{v}n_{v})^{-1})B)
|n_{v}|_{v}^{r'}
\inf_{\mathmakebox[0.5cm]{t_{v} \in T_{I}(\rmF_{v})}}|t_{v}m_{v}|_{v}^{r''}dn_{v}dm_{v}
\]
ce qui converge si $v \in S \smin S_{\infty}$ car toutes les intégrales sont 
sur des ensembles compacts. 
Si $v \in S_{\infty}$ cette intégrale égale
\begin{equation*}\label{eq:Varadarajan}
\int_{T_{I}(\rmF_{v}) \bsl M_{I_{1}}(\rmF_{v})}
\int_{\nl_{P}(\rmF_{v})} \Upphi_{v}(\Ad(m_{v}^{-1})B + U_{P,v})
|n_{v}|_{v}^{r'}
\inf_{\mathmakebox[0.5cm]{t_{v} \in T_{I}(\rmF_{v})}}|t_{v}m_{v}|_{v}^{r''}dU_{P,v}dm_{v}
\end{equation*}
où $\Ad(n_{v}^{-1})\Ad(m_{v}^{-1})B = \Ad(m_{v}^{-1})B + U_{P,v}$. 
Comme on l'a remarqué dans la preuve du lemme \ref{lem:ChaudnMajorAd}, 
il résulte du lemme 2.1 de \cite{chaud2} qu'il existe des constantes 
$N_{1}$, $N_{2}$ telles que
\[
|n_{v}|_{v} \le (1 + |\Ad(m_{v}^{-1})B|_{v})^{N_{1}}(1 + |U_{P,v}|_{v})^{N_{2}}.
\]
En utilisant le lemme \ref{lem:ChaudnMajorAd} ainsi que le 
fait que $\Upphi_{v}$ est dans la classe de Schwartz on se ramène a majorer:
\[
\int_{T_{I}(\rmF_{v}) \bsl M_{I_{1}}(\rmF_{v})}
(1 + |\Ad(m_{v}^{-1})B|_{v})^{-N'}dm_{v}
\int_{\nl_{P}(\rmF_{v})} (1 + |U_{P,v}|_{v})^{-N''}dU_{P,v}.
\]
pour $N'$, $N''$ arbitrairement grands. 
Il est clair que pour $N''$ suffisamment grand la deuxième intégrale est convergente. 
La première converge pour $N'$ suffisamment grand en vertu 
du théorème I.3.9 
de \cite{varadarajan}. 
\edem
 \elem
 
 Il nous reste à justifier la convergence de l'expression \eqref{eq:Ipolrewritten}. 
 
 \blem\label{lem:forFucksSake} Soient $\calI, \calJ \sbse \IoB$ tels que $\acalI \cap \acalJ = \varnothing$. 
 Alors, pour tout $x \in U(\A)$ on a:
 \[
\sum_{\delta \in U(X_{\calI}, \rmF) \bsl U(\rmF)} |\hf^{\calJ}_{\delta x}(X_{\calI})| < \infty.
 \]
 \elem 

\bdem 
Soit $P \in \calP(M_{I_{1}})$ un sous-groupe parabolique stabilisant le drapeau 
$\calF_{\calJ} \sbs \calF_{\calJ \cup \calI}$. 
En raisonnant comme dans le lemme \ref{lem:mainLemCvgRegss} et en utilisant la notation 
de ce lemme, on trouve des fonctions $\Phi \in \calS(\ul(\A))$ 
et $\phi \in \calS(\A_{I \smin \acalJ})$ telles que 
pour tout $p \in P(\A)$ et $k \in K$ on a:
\[
|\hf^{\calJ}_{pk}(X_{\calI})| 
\le e^{\upla_{\calJ}(H_{P}(p))}
|\Phi(\Ad((pk)^{-1})B) \phi((1-1_{\calJ})p^{-1} (1_{\calI} + \xi_{\varnothing}))|.
\]

On suppose que $\Phi = \otimes_{v} \Phi_{v}$ (resp. $\phi = \otimes_{v} \phi_{v}$) 
où $\Phi_{v}$ (resp. $\phi_{v}$) est une fonction dans la classe de Schwartz 
dans $\ul(\rmF_{v})$ (resp. $\rmF_{I \smin \acalJ} \otimes_{\rmF} \rmF_{v}$) 
si $v \in \calS_{\infty}$ et $\Phi_{v}$ (resp. $\phi_{v}$) est une fonction caractéristique
d'un compact ouvert $\kl_{v}$ dans
$\ul(\rmF_{v})$ (resp. $\fl_{v}$ dans $\rmF_{I \smin \acalJ} \otimes_{\rmF} \rmF_{v}$) 
si $v$ est une place finie de $\rmF$. 
Soit $S$ l'ensemble des places de $\rmF$ vérifiant les conditions (1) et (2) 
données dans la preuve du lemme \ref{lem:mainLemCvgRegss}.

Soit $\|\cdot \|_{\infty}$, la norme sur $\tlul(\A_{\infty})$ fixée dans le paragraphe \ref{bschwartz}.
On prétend qu'il existe des constantes $n, c > 0$ telles que
pour tout $\delta \in U(\rmF)$ on a:
\begin{equation}\label{eq:wellFoundCle}
\sum_{\eta \in T_{I_{2} \smin \acalI}(\rmF) \bsl T_{I}(\rmF)}
|\hf^{\calJ}_{\delta \eta}(X_{\calI})| \le c \|\Ad(\delta^{-1})B\|_{\infty}^{n}|\Phi(\Ad(\delta^{-1})B)|.
\end{equation}
Or, soit $\delta \in U(\rmF)$. 
Si $\Phi(\Ad(\delta^{-1})B) = 0$ le résultat est évident. 
Sinon, puisque multiplication à gauche de $\delta$ par un élément de $T_{I}(\rmF)$ ne change 
pas l'inégalité \eqref{eq:wellFoundCle}, en utilisant la définition de l'ensemble $S$, 
on peut supposer que $\delta_{v} \in K_{v}$ pour tout $v \nin S$, 
où $\delta = \prod_{v} \delta_{v}$. 
Décomposons $\delta$ selon la décomposition d'Iwasawa 
$\delta = pk$, où $p = \prod_{v} p_{v} \in P(\A)$ et $k = \prod_{v} k_{v} \in K$. 
On a:
\[
\sum_{\mathclap{\eta \in T_{I_{2} \smin \acalI}(\rmF) \bsl T_{I}(\rmF)}} \ \ 
|\hf^{\calJ}_{\delta \eta}(X_{\calI})| \le 
e^{\upla_{\calJ}(H_{P}(p))}
|\Phi(\Ad(\delta^{-1})B)| ( \ \ \sum_{\mathclap{\xi \in \rmF_{I \smin \acalJ}}}  \ \ 
|\phi((1-1_{\calJ})p^{-1}\xi )|).
\]
Soit $\calR \sbs \rmF_{I \smin \acalJ}$ un $\Z$-réseau tel que 
pour tout $\xi \in \rmF_{I \smin \acalJ}$ si 
$\phi(\xi) \neq 0$ alors $\xi \in \calR$ (voir la preuve du théorème \ref{thm:MainConv2}). 
Soit aussi $m > 0$ tel que $\sum_{\xi \in \calR} \|\xi\|^{-m} < \infty$.
Il existe un $m_{0} \in \N^{*}$ tel que pour tout $\xi \in \rmF_{I \smin \acalJ}$
si $\phi((1-1_{\calJ})p^{-1} \xi ) \neq 0$ on a 
$\xi  \in m_{0}^{-1} \calR$. Tout diviseur premier de $m_{0}$, vu comme une place finie de $\Q$, 
divise un élément de $S \smin S_{\infty}$  
car $p_{v}$ est l'identité pour $v \nin S$. 
Donc, puisque $\phi \in \calS(\A_{\calI \smin \acalJ})$, il existe des constantes positives $c_{0}$, $c_{1}$ et $r$ qui ne dépendent que de $\phi$ et $S$ 
telles que:
\[
\sum_{\mathclap{\xi \in \rmF_{I \smin \acalJ}}}  \ \ 
|\phi((1-1_{\calJ})p^{-1}\xi )|
\le c_{0} \sum_{\mathclap{\xi \in m_{0}^{-1} \calR}} \ \ 
\|p_{\infty}^{-1}\xi\|^{-m} \le c_{1} (\prod_{v \in S}|p_{v}|_{v}^{r})\sum_{\xi \in \calR} \|\xi\|^{-m}
\]
où $p_{\infty} = \prod_{v \in S_{\infty}} p_{v} \in P(\A_{\infty})$.
Si on remplace $p$ par $\eta p$ où $\eta \in T_{I}(\rmF)$ 
est tel que $\eta_{v} \fl_{v} \sps \fl_{v}$ pour $v \nin S$, 
rien ne change ci-dessus, sauf peut être l'entier $m_{0}$.
Quitte à élargir $S$, des tels $\eta$ sont denses dans $\prod_{v \in S}T_{I}(\rmF_{v})$ et on trouve, 
en passant par le lemme \ref{lem:ChaudnMajorAd}:
\[
\sum_{\mathclap{\xi \in \rmF_{I \smin \acalJ}}}  \ \ 
|\phi((1-1_{\calJ})p^{-1}\xi )| \le 
c_{2} (\prod_{v \in S} \inf_{t_{v} \in T_{I}(\rmF_{v})}|t_{v}p_{v}|_{v}^{r}) \le 
c_{3}\prod_{v \in S} (1 + |\Ad(\delta_{v}^{-1})B|_{v})^{rN_{B}}
\]
pour certains $c_{2}, c_{3} > 0$ indépendantes de $\delta$. 
Puisque on suppose $\Phi(\Ad(\delta^{-1})B) \neq 0 $, il existe une constante $C > 0$ qui ne dépend que de
$\Phi$ telle que $\prod_{v \in S \smin S_{\infty}} (1 + |\Ad(\delta_{v}^{-1}B)|_{v})^{rN_{B}} < C$. 
Pour la même raison, toujours en passant par le lemme \ref{lem:ChaudnMajorAd}, 
il existe des constantes $C', r' >0$ telles que 
$e^{\upla_{\calJ}(H_{P}(p))} \le C' \|\Ad(\delta^{-1}B)\|_{\infty}^{r'}$.
Il est clair alors qu'on a \eqref{eq:wellFoundCle}. On vient de montrer que la somme 
dans l'énoncé est majorée par la somme
\[
\sum_{\delta \in U(B, \rmF) \bsl U(\rmF)}  \| \Ad(\delta^{-1})B\|_{\infty}^{n}  |\Phi(\Ad(\delta^{-1})B)|
\]
pour un $n > 0$, ce qui converge. 
\edem

\subsection{Résultats d'holomorphie}\label{par:holoRes}

Pour $I' \sbs I_{0}$ 
notons $c_{I'}$ le volume de 
$T_{I'}(\rmF) \bsl T_{I'}(\A)^{1}$ et 
$v_{I'}$ le volume dans $\all_{I'}$ du 
parallélotope déterminé par les vecteurs 
$\{e_{i}^{\vee}\}_{i \in \calI}$ où 
$\calI \sbse I'$ est tel que $\acalI = I'$. 
Cela ne dépend pas du choix de $\calI$.

\blem\label{lem:upeta} 
Soit $\calI \sbse I_{0}$.  Alors, 
pour tout $x \in U(\A)$, l'intégrale
\[
\int_{T_{\acalI}(\rmF) \bsl T_{\acalI}(\A)}
\ind_{\calI}(H_{0}(tx))e^{\la(H_{0}(tx))}dt, \quad 
\la \in \all_{\acalI,\C}^{*}
\]
converge absolument et uniformément sur tous les 
compacts de $\Rel(\la) \in \all_{\calI}^{*}$. Elle admet 
un prolongement méromorphe, noté 
$\upeta_{\calI}$, égale à
\[
\upeta_{\calI}(\la) = 
c_{\acalI}v_{\acalI}\prod_{j \in \calI}\la(e_{j}^{\vee})^{-1}.
\]
En particulier, $\upeta_{\calI}$ ne dépend pas de $x \in U(\A)$.
\bdem 
Calcul direct. 
\edem
\elem

 Soient $\calJ_{2} \sbs \calJ_{3} \sbse I_{0}$
 et $\calJ_{1} \sbse I_{0}$ tels que $\acalJj \cap \acalJt = \varnothing$. 
 Notons $\Upsilon_{\calJ_{1}, \calJ_{2}, \calJ_{3}}(f)$ 
 la fonction holomorphe sur $\all_{|\calJ_{1} \cup \calJ_{2}|,\C}^{*}$ 
 définie par l'intégrale 
 considérée dans le lemme \ref{lem:mainLemCvgRegss} ci-dessus.

\blem\label{lem:LamPros} Soient $\calJ, \calJ_{1}, 
\calJ_{2},\calJ_{3} \sbse I_{0}$ 
tels que 
$\calJ_{1} \sqcup \calJ_{3} \sbs \calJ$ et
$\calJ_{2} \sbs \calJ_{3}$.
L'intégrale suivante
\begin{equation*}
\int\limits_{
T_{(I_{2} \smin I_{0}) \cup |\calJ \smin \calJ_{12}|}(\rmF)
T_{I_{0} \smin \acalJ}(\A)\bsl U(\A)} \! \! \! \! \! \! \! \! \! \! \! \! \! \! \!
\ind_{\calJ_{1} \cup (\calJ \smin \calJ_{1})^{\sharp}}(H_{0}(x))
e^{\la(H_{0}(x))}
\hf_{x}^{\calJ_{3}}(X_{\calJ_{1} \cup \calJ_{2}^{\sharp}})dx, \quad 
\la \in \all_{\acalJ, \C}^{*}
\end{equation*}
converge absolument et uniformément sur tous les compacts de
\begin{equation*}
\Rel(\la) \in 
\all_{|\calJ_{1} \cup \calJ_{2}|}^{*} \times 
(\upla_{\calJ_{3 \smin 2}^{\sharp}} + 
\all_{(\calJ \smin \calJ_{12})^{\sharp}}^{*})
\end{equation*}
et elle admet un prolongement méromorphe
à $\all_{\acalJ,\C}^{*}$, noté
$\La_{\calJ_{1},\calJ_{2},\calJ_{3}}^{\calJ}(f)$, 
qui vérifie
\begin{equation*}
\La_{\calJ_{1},\calJ_{2},\calJ_{3}}^{\calJ}(f)(\la) = 
\upeta_{\calJ_{3 \smin 2}^{\sharp}}
(\la_{|\calJ_{3 \smin 2}|} + 
\upla_{\calJ_{3 \smin 2}})
\upeta_{(\calJ \smin \calJ_{13})^{\sharp}}
(\la_{|\calJ \smin \calJ_{13}|})
\Upsilon_{\calJ_{1},\calJ_{2},\calJ_{3}}(f)(\la_{|\calJ_{12}|}).
\end{equation*}

\bdem 
On intègre d'abord 
sur $T_{|\calJ \smin \calJ_{12}|}(\rmF) 
\bsl T_{|\calJ \smin \calJ_{12}|}(\A)$
ce qui donne
\[
\int_{[T_{|\calJ \smin \calJ_{12}|}]}
\ind_{(\calJ \smin \calJ_{12})^{\sharp}}(H_{0}(tx))
e^{(\la_{|\calJ \smin \calJ_{12}|} + \upla_{\calJ_{3 \smin 2}})(H_{0}(tx))}
dt
\]
ce qui converge, en vertu du lemme \ref{lem:upeta}, pour 
$\Rel(\la_{|\calJ \smin \calJ_{12}|}) 
\in (\upla_{\calJ_{3\smin 2}^{\sharp}} + 
\all_{(\calJ \smin \calJ_{12})^{\sharp}}^{*})$ et admet 
un prolongement méromorphe égale 
à $\upeta_{(\calJ \smin \calJ_{12})^{\sharp}}
(\la_{|\calJ \smin \calJ_{12}|} + \upla_{\calJ_{3 \smin 2}})$.

L'intégrale qui reste à calculer c'est précisément 
l'intégrale considérée dans le lemme \ref{lem:mainLemCvgRegss}, \cad 
$\Upsilon_{\calJ_{1},\calJ_{2},\calJ_{3}}(f)(\la_{|\calJ_{12}|})$, 
d'où la convergence voulue. Finalement, 
on a 
\[
\upeta_{(\calJ \smin \calJ_{12})^{\sharp}}
(\la_{|\calJ \smin \calJ_{12}|} + \upla_{\calJ_{3 \smin 2}}) = 
\upeta_{\calJ_{3 \smin 2}^{\sharp}}
(\la_{|\calJ_{3 \smin 2}|} + 
\upla_{\calJ_{3 \smin 2}})
\upeta_{(\calJ \smin \calJ_{13})^{\sharp}}
(\la_{|\calJ \smin \calJ_{13}|})
\]
d'où le résultat.
\edem
\elem

En vertu du lemme \ref{lem:LamPros} ci-dessus, 
si $\calJ_{3} = \calJ \smin \calJ_{1}$ 
la fonction $\La_{\calJ_{1},\calJ_{2},\calJ \smin \calJ_{1}}^{\calJ}(f)$ 
est holomorphe en $\la = 0$ et on la note simplement $\Lambda_{\calJ_{1}, \calJ_{2}}^{\calJ}(f)$. 
Dans ce cas il y a une autre représentation intégrale de 
$\Lambda_{\calJ_{1}, \calJ_{2}}^{\calJ}(f)$ convergente 
sur un ouvert qui contient zéro.

\blem\label{lem:UpsPros}
 Soient $\calJ, \calJ_{1}, \calJ_{2} \sbse I_{0}$ 
tels que $\calJ_{1} \sqcup \calJ_{2} \sbs \calJ$.
L'intégrale suivante
\begin{equation*}
\int\limits_{
T_{(I_{2} \smin I_{0}) \cup |\calJ \smin \calJ_{12}|}(\rmF)
T_{I_{0} \smin \acalJ}(\A)\bsl U(\A)} \! \! \! \! \! \! \! \! \! \! \! \! \! \! \! \! \! \! \! \! 
\ind_{(\calJ \smin \calJ_{2}) \cup \calJ_{2}^{\sharp}}(H_{0}(x))
e^{\la(H_{0}(x))}
\hf_{x}^{\calJ \smin \calJ_{1}}(X_{\calJ_{1} \cup \calJ_{2}^{\sharp}})dx, \quad 
\la \in \all_{\acalJ, \C}^{*}
\end{equation*}
converge absolument et uniformément sur tous les compacts de
\begin{equation*}
\Rel(\la) \in 
\all_{|\calJ_{1} \cup \calJ_{2}|}^{*} \times 
(\upla_{(\calJ \smin \calJ_{12})^{\sharp}} + 
\all_{\calJ \smin \calJ_{12}}^{*}).
\end{equation*}
En particulier, elle converge pour $\la = 0$.
L'intégrale, en tant qu'une fonction de variable 
$\la$, admet un prolongement méromorphe
à $\all_{\acalJ,\C}^{*}$, noté
$\brLa_{\calJ_{1},\calJ_{2}}^{\calJ}(f)$,
qui vérifie
\begin{gather*}
\brLa_{\calJ_{1},\calJ_{2}}^{\calJ}(f)(\la) = 
\upeta_{\calJ \smin \calJ_{12}}
(\la_{|\calJ \smin \calJ_{12}|} + 
\upla_{\calJ \smin \calJ_{12}})
\Upsilon_{\calJ_{1},\calJ_{2}, \calJ \smin \calJ_{1}}(f)(\la_{|\calJ_{12}|}), \\
\brLa_{\calJ_{1},\calJ_{2}}^{\calJ}(f)(\la) = 
(-1)^{\# (\calJ \smin \calJ_{12})}\La_{\calJ_{1},\calJ_{2}}^{\calJ}(f)(\la).
\end{gather*}

\bdem
 Toutes les assertion, sauf la dernière, 
se démontrent de même façon 
que dans le lemme \ref{lem:LamPros} et la dernière découle 
du fait que $\upeta_{\calI} = (-1)^{\#\calI} \upeta_{\calI^{\sharp}}$ 
pour tout $\calI \sbse I_{0}$. 
\edem
\elem

\brop\label{prop:zetaDefProps} 
Soit $\calJ \sbse I_{0}$.
Alors, l'intégrale 
\[
\vol(T_{I_{2} \smin I_{0}}(\rmF) \bsl T_{I_{2} \smin I_{0}}(\A))
\int\limits_{\mathclap{
U(\A, X_{\calJ})\bsl U(\A)}}
f(x^{-1}X_{\calJ}x)
e^{\la(H_{0}(x))}dx, \quad \la \in \all_{\acalJ,\C}^{*}
\]
converge absolument 
et uniformément sur tous les compacts de
\[
\Rel(\la) \in 
\upla_{\calJ^{\sharp}} + \all_{\calJ^{\sharp}}^{*}
\]
et admet un prolongement méromorphe 
à $\all_{\acalJ,\C}^{*}$, noté $\zeta_{\calJ}(f)$,
qui vérifie
\[
\zeta_{\calJ}(f) = 
\sum_{\calJ_{1} \sqcup \calJ_{3} \sbs \calJ}
(-1)^{\#(\calJ \smin \calJ_{13})}
\sum_{\calJ_{2} \sbs \calJ_{3}}
\La_{\calJ_{1},\calJ_{2},\calJ_{3}}^{\calJ}(f).
\]

\bdem 
Remarquons d'abord que l'intégrale considérée dans la proposition c'est juste
\[
\int\limits_{\mathclap{
T_{I_{2} \smin I_{0}}(\rmF)T_{I_{0} \smin \acalJ}(\A)\bsl U(\A)}}
f_{x}(X_{\calJ})
e^{\la(H_{0}(x))}
dx.
\]
Pour tout $\calJ_{1} \sbs \calJ$ et tout $x \in U(\A)$,
en utilisant l'égalité (\ref{eq:PoissonInlcExcl1}) du lemme \ref{lem:PoissonInlcExcl}, on a
\begin{equation}\label{eq:anothPoiss}
\sum_{\eta \in T_{\acalJ}(\rmF)}f_{\eta x}(X_{\calJ})
= \sum_{\calJ_{3} \sbs \calJ \smin \calJ_{1}}
(-1)^{\#(\calJ \smin \calJ_{13})}
\sum_{\calJ_{2} \sbs \calJ_{3}}
\sum_{\eta \in T_{|\calJ_{12}|}(\rmF)}
\hf_{\eta x}^{\calJ_{3}}(X_{\calJ_{1} \cup \calJ_{2}^{\sharp}}).
\end{equation}
Remarquons aussi que pour tout $\calI \sbse I_{0}$ tel que $\acalI = \acalJ$ et
tout $\eta \in T_{\acalJ}(\rmF)$ on a 
$\ind_{\calI}(H_{0}(\eta x)) = \ind_{\calI}(H_{0}(x))$. 
En utilisant ceci, l'égalité (\ref{eq:anothPoiss}) ci-dessus, 
ainsi que l'égalité (\ref{eq:decOf1}) on trouve:
\begin{multline}\label{eq:zetaTateProof}
\sum_{\eta \in T_{\acalJ}(\rmF)}f_{\eta x}(X_{\calJ})= 
\sum_{\calJ_{1} \sbs \calJ}
\sum_{\eta \in T_{\acalJ}(\rmF)}f_{\eta x}(X_{\calJ})
\ind_{\calJ_{1} \cup (\calJ\smin \calJ_{1})^{\sharp}}(H_{0}(\eta x))
= \\
\sum_{\calJ_{1} \sqcup \calJ_{3} \sbs \calJ}
(-1)^{\#(\calJ \smin \calJ_{13})}
\sum_{\calJ_{2} \sbs \calJ_{3}}
\sum_{\eta \in T_{|\calJ_{12}|}(\rmF)}
\hf_{\eta x}^{\calJ_{3}}(X_{\calJ_{1} \cup \calJ_{2}^{\sharp}})
\ind_{\calJ_{1} \cup (\calJ \smin \calJ_{1})^{\sharp}}(H_{0}(\eta x)).
\end{multline}
Soit $\la \in \all_{\acalJ,\C}^{*}$. Pour tout $\eta \in T_{\acalJ}(\rmF)$ 
on a $\la(H_{0}(\eta x)) = \la(H_{0}(x))$. 
On multiplie alors l'égalité (\ref{eq:zetaTateProof}) ci-dessus par $e^{\la(H_{0}(x))}$ 
et 
en utilisant le lemme \ref{lem:LamPros} on voit qu'on 
peut intégrer le terme correspondant aux $\calJ_{1} \sqcup \calJ_{3} \sbs \calJ$
et $\calJ_{2} \sbs \calJ_{3}$ sur $T_{(I_{2} \smin I_{0}) \cup \acalJ}(\rmF)
T_{I_{0} \smin \acalJ}(\A)\bsl U(\A)$
pour $\la \in \all_{\acalJ,\C}^{*}$ tels que $\Rel(\la)$ appartient à:
\[
\all_{|\calJ_{1} \cup \calJ_{2}|}^{*} \times 
(\upla_{\calJ_{3 \smin 2}^{\sharp}} + 
\all_{(\calJ \smin \calJ_{12})^{\sharp}}^{*}) \sps 
\upla_{\calJ^{\sharp}} + \all_{\calJ^{\sharp}}^{*}.
\]
Tout est alors intégrable sur $\Rel(\la) \in \upla_{\calJ^{\sharp}} + \all_{\calJ^{\sharp}}^{*}$ et on obtient 
le résultat voulu en invoquant de nouveau le lemme \ref{lem:LamPros}. 
 \edem
 \erop
 
 \blem\label{lem:sumZetasUpsilon}
 Pour tout $I' \sbs I_{0}$
  on a l'égalité de fonctions méromorphes
sur $ \all_{I',\C}^{*}$
 \[
\sum_{\acalJ = I'}\zeta_{\calJ}(f)(\la) =  
\sum_{\acalJ = I'}
\sum_{\calJ_{1} \sqcup \calJ_{2} \sbs \calJ}
(-1)^{\#(\calJ \smin \calJ_{12})}
\brLa_{\calJ_{1}, \calJ_{2}}^{\calJ}(f)(\la).
 \] 
 En particulier, en vertu du lemme 
 \ref{lem:UpsPros}, l'expression ci-dessus est holomorphe 
 en $ \la = 0$.

 \bdem 
En utilisant la proposition 
\ref{prop:zetaDefProps} 
on a que pour 
tout $\la \in \all_{I',\C}^{*}$ 
la somme $\sum_{\acalJ  = I'}\zeta_{\calJ}(f)(\la)$ 
égale
\begin{multline*}
\sum_{\acalJ  = I'}
\sum_{\calJ_{1} \sqcup \calJ_{3} \sbs \calJ}
(-1)^{\#(\calJ \smin \calJ_{13})}
\sum_{\calJ_{2} \sbs \calJ_{3}}
\La_{\calJ_{1},\calJ_{2},\calJ_{3}}^{\calJ}(f)(\la) = 
\sum_{\calJ_{1} \sqcup \calJ_{3} \sbse I'}
(-1)^{\#(I' \smin |\calJ_{13}|)} \\
\sum_{\calJ_{2} \sbs \calJ_{3}}
\upeta_{\calJ_{3 \smin 2}^{\sharp}}
(\la_{|\calJ_{3 \smin 2}|} + 
\upla_{\calJ_{3 \smin 2}})
\Upsilon_{\calJ_{1},\calJ_{2},\calJ_{3}}(f)(\la_{|\calJ_{12}|})
\sum_{\acalJc = I' \smin |\calJ_{13}|}
\upeta_{\calJ_{4}}(\la_{\acalJc})
\end{multline*}
où l'on a utilisé le lemme \ref{lem:LamPros} 
pour faire apparaître les fonctions $\upeta$.

On prétend que 
\begin{equation}\label{eq:arthurSumOnly}
\sum_{\acalJc = I' \smin |\calJ_{13}|}
\upeta_{\calJ_{4}}(\la_{\acalJc}) \equiv 0, \quad 
\text{ si }I' \smin |\calJ_{13}| \neq \varnothing.
\end{equation}
En effet, il résulte du lemme \ref{lem:upeta} que si 
$\calI_{1} \sqcup \calI_{2} \sbse I'$ alors 
$\upeta_{\calI_{1} \sqcup \calI_{2}} = (-1)^{\#\calI_{2}}
\upeta_{\calI_{1} \sqcup \calI_{2}^{\sharp}}$. 
On voit que (\ref{eq:arthurSumOnly}) revient 
à montrer que $\sum_{S \sbs I' \smin |\calJ_{13}|}(-1)^{\#S} = 0$ 
si $I' \smin|\calJ_{13}| \neq \varnothing$
ce qui est 
juste l'identité (\ref{basicidentity}). 

On obtient donc
 \[
\sum_{\acalJ = I'}\zeta_{\calJ}(f) =  
\sum_{\acalJ = I'}
\sum_{\calJ_{1} \sqcup \calJ_{2} \sbs \calJ}
\La_{\calJ_{1}, \calJ_{2}}^{\calJ}(f)
 \] 
et le résultat suit du lemme \ref{lem:UpsPros}. 
\edem
 \elem

%



\end{document}